\documentclass[draft, 12pt]{article}
\usepackage[centertags]{amsmath}
\usepackage{amsfonts}
\usepackage{amssymb}
\usepackage{latexsym}
\usepackage{amsthm}
\usepackage{newlfont}
\usepackage{graphicx}
\usepackage{listings}
\usepackage{booktabs}
\usepackage{abstract}
\usepackage{xcolor}

\usepackage[normalem]{ulem}
\date{}

\newlength{\defbaselineskip}
\newcommand{\setlinespacing}[1]%
           {\setlength{\baselineskip}{#1 \defbaselineskip}}

\newcommand{\bR}{{\mathbb{R}}}

\newcommand{\N}{{\mathbb{N}}}

\newcommand{\actaqed}{\hfill $\actabox$}
{\medskip\noindent \textit{Proof of #1. }}%
{\actaqed \medskip}

\def\D{{\mathcal D}}

\def\cB{{\mathcal B}}
\def\cC{{\mathcal C}}

\def \Tr{\mathcal T}

\def \cP{\mathcal P}

\def \U{\mathcal U}
\def \E{\mathcal E}
\def \cE{\mathcal E}

\def \cM{\mathcal M}
\def \cN{\mathcal N}
\def \cF{\mathcal F}
\def \cS{\mathcal S}
\def\R{{\mathbb R}}
\def\Z{\mathbb Z}

\def \T{\mathbb T}
\def\bP{\mathbb P}
\def\bE{\mathbb E}
\def \bbE{\mathbb E}
\def\bbC{\mathbb C}
\def \<{\langle}
\def\>{\rangle}

\def \La{\Lambda}

\def \e{\varepsilon}
\def \de{\delta}

\def \ff{\varphi}
\def\la{\lambda}
\def \ro{\varrho}

\def \sp{\operatorname{span}}
\def \rank{\operatorname{rank}}

\def\bx{\mathbf x}
\def\by{\mathbf y}
\def\bz{\mathbf z}
\def\bk{\mathbf k}
\def\bv{\mathbf v}
\def\bu{\mathbf u}
\def\bw{\mathbf w}

\def\bs{\mathbf s}
\def\bb{\mathbf b}
\def\bfe{\mathbf e}
\def\bN{\mathbf N}
\def\bW{\mathbf W}
\def\bH{\mathbf H}

\def\bF{\mathbf F}

\def\btt{\mathbf t}

\def\U{{\mathcal U}}

\def\ga{\gamma}

\newtheorem{Theorem}{Theorem}[section]
\newtheorem{Lemma}{Lemma}[section]

\newtheorem{Definition}{Definition}[section]
\newtheorem{Proposition}{Proposition}[section]
\newtheorem{Remark}{Remark}[section]

\newtheorem{Corollary}{Corollary}[section]

\numberwithin{equation}{section}

\newcommand{\be}{\begin{equation}}
\newcommand{\ee}{\end{equation}}

\newcounter{om}
\newcommand{\om}{\par\refstepcounter{om}
\textbf{M.\arabic{om}.} }

\begin{document}

\title{Sampling discretization and related problems}

\author{B. Kashin\thanks{Steklov Institute of Mathematics, Moscow Center for Fundamental and Applied Mathematics,},\, E. Kosov\thanks{Lomonosov Moscow State University, National Research University Higher School of Economics,},\, I. Limonova\thanks{Lomonosov Moscow State University, Moscow Center for Fundamental and Applied Mathematics,}, \,  and V. Temlyakov\thanks{University of South Carolina, Steklov Institute of Mathematics, Lomonosov Moscow State University, and Moscow Center for Fundamental and Applied Mathematics.} } \maketitle

\begin{abstract}
	{This survey addresses sampling discretization and its connections with other areas of mathematics. {The survey concentrates on sampling discretization of norms of elements of finite-dimensional subspaces.} We present here known results on sampling discretization of both integral norms and  the uniform norm beginning with classical results and ending with very recent achievements.   We also show how  sampling discretization  connects to spectral properties and operator norms of submatrices, embedding of finite-dimensional subspaces, moments of marginals of high-dimensional distributions, and learning theory.  Along with the corresponding results, important techniques for proving those results are discussed as well.   }
\end{abstract}


 \newpage

 \tableofcontents

\section{Introduction}
\label{Int}

 Discretization is an essential step 
  in making a continuous problem computationally feasible. The construction of good sets of points in a multidimensional domain is a fundamental problem of pure and computational mathematics.

A classical example of a discretization problem is to estimate metric entropy (covering numbers, entropy numbers).
Bounds for the $\varepsilon$-entropy of function classes are inherently important  and also  important in connection with other fundamental problems.
   Numerical integration is another prominent example of a discretization problem.  We now briefly illustrate difficulties in finding good points for covering and for numerical integration on the example  of discretizing the $d$-dimensional unit cube $[0,1]^d$.
There are different ways of doing that. For example, we can interpret $[0,1]^d$ as a compact set of $\R^d$ and use the idea of {\it covering numbers} (metric entropy). With this approach, for instance
in the case of $\ell_\infty$ norm, we can find optimal coverings. For a given $n\in\N$ the regular grid with coordinates at the centers of intervals $[(k-1)/n,k/n]$, $k=1,\dots,n$, provides a net for
an optimal $\ell_\infty$ covering with the number of points $N=n^d$. Very often the unit cube
$[0,1]^d$ plays the role of a domain, where smooth functions of $d$ variables are defined and we are interested in discretizing some continuous operations with these functions. Numerical integration is such an example.
It turns out that the regular grids mentioned above are very far from being a good discretization of $[0,1]^d$ for numerical integration {of functions with mixed smoothness}.

A problem of optimal recovery is one more example of a discretization problem. This problem turns out to be very difficult for the mixed smoothness classes of multivariate functions, which are important in applications. It is not solved even in the case of the $L_2$ norm.

A standard approach to solving a continuous problem numerically -- the Galerkin method --
suggests searching for an approximate solution from a given finite-dimensional subspace. A standard way to measure an error of approximation is an appropriate $L_q$ norm, $1\le q\le\infty$. Thus, the problem of   discretization of the $L_q$ norms of functions from a given finite-dimensional subspace arises in a very natural way.
The first results in this direction were obtained  in the 1930s by Bernstein, by Marcinkiewicz  and
by Marcinkiewicz--Zygmund  for discretization of the $L_q$ norms of the univariate trigonometric polynomials. Even though this problem is very important in applications, its systematic study has begun only  recently.

The problem of arranging points in a multidimensional domain, in particular the problem of optimal spherical codes, is also a fundamental problem in coding theory.  Finding optimal spherical codes is equivalent to building
large incoherent dictionaries in $\R^d$ in compressed sensing. Compressed sensing is an actively developing area of numerical analysis, which is very important in processing large data sets, in particular in signal and image processing {(see, for instance, \cite{FR} and \cite{VTbook}, Ch. 5)}.

The goals of this paper are to survey the results, and connect together  ideas, methods, and results from different areas of research {related to problems of discretization and recovery in the case of finite-dimensional subspaces. We note that this survey does not cover some classical areas where discretization problems are studied. For instance, we do not discuss discretization in spaces of analytic functions. The reader can find such results in \cite{Seip} and \cite{OU}.} Many results on discretization {discussed in this survey} are obtained with the aid of probabilistic technique. The authors of those results use a variant of the concentration of measure phenomenon, the chaining technique, and bounds of the entropy numbers. We note that the fundamental idea of chaining goes back to the 1930s, when it was suggested by A.N. Kolmogorov in 1934 for finding sufficient conditions for almost everywhere continuity of trajectories of random processes. The first publication of his proof (with his approval) appeared in the paper by E.~Slutsky
\cite{Sluts1} (see also \cite{Sluts2}).
 Similar probabilistic techniques were used in
other areas of research: Spectral properties and operator norms of submatrices, embedding of finite-dimensional subspaces into $\ell_q^m$, moments of marginals of high-dimensional distributions, and learning theory.

We now proceed to precise settings and to a more detailed discussion.

{\bf Discretization.}
Let $\Omega$ be a compact subset of $\R^d$ with the probability measure $\mu$. By {\it sampling discretization in a broad sense} we understand replacement of the measure $\mu$ by a discrete measure $\mu_m$ supported on a set of $m$ points from $\Omega$, which keeps some properties of  $\mu$.   We would like {a} discrete replacement $\mu_m$ to represent (approximate) the original measure $\mu$ well. A specific (narrow) sampling discretization setting is determined by a criterion on the error of approximation of $\mu$ by $\mu_m$. There are classical ways to measure how close two measures are, for example, the Hellinger distance and the Kullback-Leibler information ({see, for instance, \cite{KuLe} and the books \cite{Poll} and} \cite{VTbook}, Section 4.2.3), which are important in probability theory and statistics. Our criteria of good approximation of $\mu$ by $\mu_m$ are motivated by applications in approximation theory and numerical analysis. We begin with
the classical problem of numerical integration, which provides an example of such criterion.  Numerical integration seeks good ways of approximating an integral
$$
\int_\Omega f(\bx)d\mu(\bx)
$$
by an expression of the form
\be\label{1.1ni}
\La_m(f,\xi) :=\sum_{j=1}^m\la_jf(\xi^j),\qquad \xi=\{\xi^j\}_{j=1}^m \subset \Omega,\quad j=1,\dots,m.
\ee
It is clear that we must assume that $f$ is integrable and defined at the points
 $\xi^1,\dots,\xi^m$. Expression~\eqref{1.1ni} is called a {\it cubature formula} $(\xi,\La)$ (if $\Omega \subset \R^d$, $d\ge 2$) or a {\it quadrature formula} $(\xi,\La)$ (if $\Omega \subset \R$) with nodes $\xi=\{\xi^j\}_{j=1}^m$ and weights $\La:=\{\la_j\}_{j=1}^m \subset \R$ or $\subset\bbC$. We do not impose any {\it a priori} restrictions on nodes and weights. Some nodes may coincide and both positive and negative weights are allowed. The above setting means that in the integration problem we replace the measure
$\mu$ by the discrete measure (signed measure) $\mu_m$ such that $\mu_m(\xi^j)=\la_j$, $j=1,\dots,m$.

Some classes of cubature formulas are of special interest. For instance,  cubature formulas, which have equal weights $1/m$, are important in applications. We use a special notation for these cubature formulas
$$
Q_m(f,\xi) :=\frac{1}{m}\sum_{j=1}^mf(\xi^j).
$$
In this notation $Q$ stands for Quasi-Monte Carlo -- a standard in some areas of research name for the above cubature formulas.

Typically, one is interested in {\it good cubature formulas} for a given function class. The term {\it good} can be understood in different ways. Cubature formulas providing
exact numerical integration for functions from a given class can be considered {\it best}. If a cubature formula is not exact on a given class, then we need to introduce a concept of error.  Following the standard approach, for a function class $\bW$ we introduce the   errors of the cubature formulas $\La_m(\cdot,\xi)$ and $Q_m(\cdot,\xi)$ by
$$
\La_m(\bW,\xi):= \sup_{f\in \bW} \left|\int_\Omega fd\mu -\La_m(f,\xi)\right|,
$$
$$
Q_m(\bW,\xi):= \sup_{f\in \bW} \left|\int_\Omega fd\mu -Q_m(f,\xi)\right|.
$$
The quantities $\La_m(\bW,\xi)$ and $Q_m(\bW,\xi)$  are classical characteristics of the quality of given cubature formulas. This setting is called {\it the worst case setting} in
the Information Based Complexity (see, e.g., \cite{Wo} and {\cite{TWW}}). Notice that the above characteristics provide
an absolute error independent of an individual function from the class. A typical class $\bW$ in the numerical integration problem is a smoothness class of functions continuous on $\Omega$, for instance, it might be the unit ball of a Sobolev, Nikol'skii, or Besov space.  Thus, for a given class $\bW$ the error $\La_m(\bW,\xi)$ shows how close measures $\mu$ and $\mu_m$ are, where $\mu_m$ is the discrete measure, concentrated on the set $\xi=\{\xi^j\}_{j=1}^m$ with weights $\La:=\{\la_j\}_{j=1}^m$, i.e. $\mu_m(\xi^j)=\la_j$, $j=1,\dots,m$. In particular, if $\bW$ is a unit ball of a Banach space embedded in $\cC(\Omega)$, then $\La_m(\bW,\xi)$ is the norm of the difference $\mu-\mu_m$ in the dual norm.

Recently, in a number of papers (see \cite{VT158}, \cite{VT159}, \cite{VT160}, \cite{DPTT}, \cite{VT171}, \cite{DPSTT1}, \cite{DPSTT2}, \cite{Kos}) a systematic study of a variant of the numerical integration problem has begun. It is the problem of discretization of the $L_q$ norms of elements from a given function class $\bF$. By $L_q$, $1\le q< \infty$, norm we understand
$$
\|f\|_q:=\|f\|_{L_q(\Omega,\mu)} := \left(\int_\Omega |f|^qd\mu\right)^{1/q}.
$$
By $L_\infty$ norm we understand the uniform norm of continuous functions
$$
\|f\|_\infty := \max_{\bx\in\Omega} |f(\bx)|
$$
and with some abuse of notation we occasionally write $L_\infty(\Omega)$ for the space $\cC(\Omega)$ of continuous functions on $\Omega$.

By discretization of the $L_q$ norm we understand a replacement of the measure $\mu$ by
a discrete measure $\mu_m$ with support on a set $\xi =\{\xi^j\}_{j=1}^m \subset \Omega$   in such a way that the error $|\|f\|^q_{L_q(\Omega,\mu)} - \|f\|^q_{L_q(\Omega,\mu_m)}|$ is small for $f\in\bF$. This means that integration with respect to $\mu$ is replaced by an appropriate cubature formula. Thus, integration is replaced by evaluation of a function $f$ at a
finite set $\xi$ of points. This is why this way of discretization is called {\it sampling discretization}.  The reader can find a corresponding discussion in a recent survey \cite{DPTT}.
{The first results in sampling discretization were obtained by Bernstein \cite{Bern1} and \cite{Bern2} (see also \cite{Z}, Ch.10, Theorem (7.28)) for discretization of the uniform norm of the univariate trigonometric polynomials in 1931-1932. In this case $\bF$ is the unit $L_\infty$-ball of  a subspace of the univariate trigonometric polynomials of a given degree.}
The first results in sampling discretization of the integral norms $L_q$, $1\le q<\infty$, were obtained by Marcinkiewicz ($1<q<\infty$) and
by Marcinkiewicz--Zygmund ($1\le q\le \infty$) (see \cite{Z}, Ch.10, Theorems (7.5) and (7.28)) for discretization of the $L_q$ norms of the univariate trigonometric polynomials in 1937. Therefore, we also call sampling discretization results  {\it Bernstein-type theorems} in the case of discretization of the uniform norm of functions from a finite-dimensional subspace and {\it Marcinkiewicz-type theorems} in the case of integral norms of functions from a finite-dimensional subspace (see \cite{VT158}, \cite{VT159}, \cite{DPTT}). In the literature this kind of results is also known under the name {\it Marcinkiewicz--Zygmund inequalities} (see, for instance, \cite{DW}, \cite{Gro} and references therein). Recently, a substantial progress in sampling discretization has been made in \cite{VT158}, \cite{VT159}, \cite{VT168}, \cite{DPTT},  \cite{DPSTT1}, \cite{DPSTT2}, \cite{Kos}, and \cite{LT}.

It will be convenient for us to use the following notation. For  given sets $\xi:=\{\xi^j\}_{j=1}^m$ of sampling points and $\La:=\{\la_j\}_{j=1}^m\subset \R$ of weights,  and for a class $\bF\subset L_q(\Omega,\mu)$, $1\le q<\infty$, of continuous on $\Omega$ functions we define
$$
er(\bF,\xi,L_q):=  \sup_{f\in \bF} \left|\|f\|_q^q - \frac{1}{m}\sum_{j=1}^m |f(\xi^j)|^q\right|
$$
and
$$
er(\bF,\xi,\La,L_q):=  \sup_{f\in \bF} \left|\|f\|_q^q - \sum_{j=1}^m \la_j |f(\xi^j)|^q\right|.
$$
It is clear that the above errors of discretization on the class $\bF$ coincide with the numerical integration errors on the class $\bW=\bF^q:= \{|f|^q\,:\, f\in \bF\}$
$$
er(\bF,\xi,L_q) = Q_m(\bF^q,\xi),\qquad er(\bF,\xi,\La,L_q)=\La_m(\bF^q,\xi).
$$
A special case, when $\bF$ is a unit ball of a given $N$-dimensional subspace $X_N$ of $L_q(\Omega)$, plays a very important role in sampling discretization. Denote
$X_N^q := \{f\,:\, f\in X_N,\,\, \|f\|_q \le 1\}$. If $er(\bF,\xi,L_q)=\e$, then it is guaranteed that
for any $f\in \bF$ we have the following error bound, which does not depend on $f$ and only depends on the function class $\bF$
\be\label{ar}
\left|\|f\|_q^q - \frac{1}{m}\sum_{j=1}^m |f(\xi^j)|^q\right| \le \e.
\ee
In the case $\bF=X_N^q$, $\e<1$, we can write (\ref{ar}) for $f$ such that $\|f\|_q=1$ and
obtain the following relative error bounds
\be\label{rr}
(1-\e)\|f\|_q^q \le \frac{1}{m}\sum_{j=1}^m |f(\xi^j)|^q \le (1+\e)\|f\|_q^q\quad\text{for all}\quad f\in X_N.
\ee
Inequalities of the form (\ref{rr}) are known under the name of the Marcinkiewicz-Zygmund inequalities. We point out that inequalities (\ref{rr}) consist of two inequalities -- the left inequality, which provides the lower bound for the discrete norm, and the right inequality, which provides the upper bound for the discrete norm. These one-sided inequalities are also of interest and they are discussed here as well. 

{In this survey we concentrate on results on sampling discretization of the integral norms of elements of finite-dimensional subspaces. There are some known results 
on the behavior of the quantities $er(\bF,\xi,L_q)$ and $er(\bF,\xi,\La,L_q)$ for $\bF$ being 
an infinitely dimensional function class (see, for instance,  \cite{VT171}). We do not discuss these results here.
We stress that in the above setting we discretize the $\|f\|_q^q$. There are results on optimal estimation of the $\|f\|$ under assumption that $f\in \bF$ (see \cite{Was}). At a first glance the problems of estimation of $\|f\|$ and, say, estimation of $\|f\|^2$, are very close. A simple inequality $|a^2-b^2| \le 2M|a-b|$ for numbers satisfying $|a|\le M$ and $|b|\le M$ shows that normally we can get an upper bound for estimation of $\|f\|^2$ in terms of the error of estimation of $\|f\|$. However, it turns out that the problems of optimal errors in sampling discretization of $\|f\|$ and $\|f\|^2$ are different as it is explained in \cite{VT171}. 
 }

The main goal of the sampling discretization theory is to find good sampling points for a given function class $\bF$. Naturally, this goal includes two steps. The first step consists of finding good bounds for the following optimal quantities
$$
er_m(\bF,L_q):= \inf_{\xi^1,\dots,\xi^m} er(\bF,\xi,L_q),
$$
$$
er_m^o(\bF,L_q):= \inf_{\substack{\xi^1,\dots, \xi^m\\\la_1,\dots, \la_m}} er(\bF,\xi,\La,L_q).
$$
{We point out that the superscript $o$ in $er_m^o(\bF,L_q)$ stands for {\it optimal} over both the point sets and the weights. This notation is similar to the corresponding notations in numerical integration and discrepancy theory (see \cite{VTbookMA}, p. 316)}.
The second step consists of finding sampling points, which provide discretization errors close to
the optimal ones.  {The tradeoff between accuracy and complexity is also an interesting question. In this survey we mostly present results on the first step. }

{In the above setting we are interested in the behavior of the error characteristics, which are optimal 
either over the point sets ($er_m$) or over both the point sets and the weights ($er_m^o$). 
The setting, when the point set $\xi:=\{\xi^1,\dots,\xi^m\}$ is given (deterministically or randomly) 
and we want either to control  the error $er(\bF,\xi,L_q)$ or to prove inequalities (\ref{rr}) is also interesting and important in applications (see, for instance, \cite{BaGr}, \cite{Gr}, \cite{O-CS}, \cite{PST}, \cite{SmZ}, and references therein) but not considered in the survey. This setting falls into {\it sampling theory} and the corresponding Marcinkiewicz-Zygmund inequalities are known under the name 
{\it sampling theorem}.
}

In Section~\ref{A} we discuss results on
sampling discretization of the $L_q$, $1\le q\le \infty$, norm of elements of $N$-dimensional subspaces $X_N$ of
the space $\cC(\Omega)$  of functions continuous on $\Omega$. We call problems of that type the {\it Marcinkiewicz discretization problems}.  These problems concern the upper bounds on the errors
$er_m(\bF,L_q)$ and $er_m^o(\bF,L_q)$ where $\bF=X_N^q$ is the $L_q$ unit ball of the subspace $X_N$.

In Section \ref{A} we present two types of results on sampling discretization. Some results involve sampling discretization for specific
$N$-dimensional subspaces, for instance, the trigonometric polynomials with frequencies from
parallelepipeds or hyperbolic crosses, algebraic polynomials. In this case techniques from harmonic analysis, approximation theory, and number theory are used. Other results involve sampling discretization in general $N$-di\-men\-sional subspaces, which may satisfy some natural conditions such as the Nikol'skii-type inequality or inequalities on the entropy numbers of the $L_q$ unit ball of the $X_N$ in the uniform norm. In this case techniques from functional analysis and probability theory are used.

{\bf Submatrices.} Section \ref{M} is devoted to results on spectral properties, operator norms of submatrices, and related quantities of a given matrix, which are directly connected to the sampling discretization in finite-dimensional subspaces. {The point is that the problem of discretizing functional systems can usually be reduced to the corresponding matrix problem. For many problems a given subspace can be replaced by the corresponding subspace consisting of piecewise constant functions with general metric properties to be preserved. An advantage of the matrix setting is that its analysis allows us to apply the iteration technique (see Section \ref{M} and Subsection \ref{Em} for a detailed discussion).} 

Let us consider sampling
discretization of the $L_q$ norm in the space 
$X_N~=~\sp\{\psi_1,\dots,\psi_N\}$ (see (\ref{rr}) above and (\ref{A1}) below). We want to find a set of $m$ points $\xi=\{\xi^j\}_{j=1}^m$ such that
\be\label{sub1}
C_1\|f\|_q^q \leq \frac{1}{m} \sum_{j=1}^m |f(\xi^j)|^q \leq C_2\|f\|_q^q, \qquad \forall f\in X_N.
\ee
 Then for $f=\sum_{i=1}^N b_i\psi_i$ and for the matrix
$S:=[\psi_i(\xi^j)]$, $j=1,\dots, m$, 
${i=1,\dots, N}$, with $m$ rows (enumerated by $j$) and $N$ columns (enumerated by $i$) inequalities (\ref{sub1}) mean that
\be\label{sub2}
C_1\|f\|_q^q \leq   \|S\bb\|_{\ell_{q,m}}^q\leq C_2\|f\|_q^q, \qquad {\forall} \bb=(b_1,\dots, b_N),
\ee
where for $\by=(y^1,\dots, y^m)\in \bbC^m$
$$
\|\by\|_{\ell_{q,m}} :=\left( \frac{1}{m} \sum_{j=1}^m |y^j|^q\right)^{1/q}.
$$
{\sloppy

}
Assume now that all $\psi_i$ are piecewise constant and  the set $\Omega$ can be split $\Omega=\sqcup_{j=1}^M \omega_j$, $\mu(\omega_j)=1/M$, $1\leq j\leq M$, in such a way that
$$
\psi_i(\bx)= c_i^j, \qquad \bx\in\omega_j.
$$
Denote $F:=[c_i^j]$, $j\!=\!1,\dots,M$, $i=1,\dots, N$. It is clear that ${\|F\bb\|_{\ell_{q,M}}\!=\!\|f\|_q}$. Therefore,
inequalities (\ref{sub2}) mean that we are looking for a submatrix $S$ of the full matrix $F$, which represents it well{, namely:
$$
C_1\|F\bb\|_q^q \leq   \|S\bb\|_{\ell_{q,m}}^q\leq C_2\|F\bb\|_q^q, \qquad \forall\bb=(b_1,\dots, b_N)\in \bbC^N.
$$
}
For applications,  the most important case is discretization of the $L_2$ norm in subspaces of $L_2(\Omega)$. In this case we consider an orthonormal basis $\{\psi_i\}_{i=1}^N$ of a given subspace and study the spectral properties of matrices of the form
$$
[\psi_i(\xi^j)], \qquad \xi=\{\xi^j\}\subset\Omega.
$$

For what follows, we need to introduce the $(p,q)$-norm of a matrix: For an $M\times N$ matrix $A$ and $1\leq p, q\leq \infty$ define
\begin{gather}\label{pq_norm_def}
\|A\|_{(p,q)}=\sup\limits_{\|\bx\|_{\ell_p^N}\leq 1}\|A\bx\|_{\ell_q^M}.
\end{gather}
For brevity in the case $p=q=2$, we denote $\|A\|_{(2,2)}$ simply as $\|A\|$. We present known estimates for the $(p,q)$--norms of submatrices (and related  quantities) of a given matrix in Section \ref{M}.

In a number of applications it is not required to have two-sided inequalities \eqref{sub1}, it is sufficient to obtain one--sided estimates of the form
\be\label{low_est}
\|f\|_q^q\geq C\sum\limits_{j=1}^m \lambda_j|f(\xi^j)|^q, \qquad \forall f\in X_N,
\ee
or
\be\label{upp_est}
\|f\|_q^q\leq C\sum\limits_{j=1}^m \lambda_j|f(\xi^j)|^q, \qquad \forall f\in X_N.
\ee

One   example when we need to know one--sided estimates of the form \eqref{low_est} is the Large Sieve method in number theory. It is based on upper estimates of the $(2,2)$--norms of submatrices of the form
$$
[e^{ik\xi_j}], \qquad 1\leq k\leq n, \qquad \xi=\{\xi_j\}_{j=1}^m\subset(0, 2\pi),
$$
in terms of properties of the set $\xi$ (see \cite[Ch.2]{Mon} for details).
Another example of application of (only) one-sided estimates \eqref{low_est} is Marcinkiewicz's proof of convergence in the  $L_p$, $1\leq p<\infty$, of interpolation polynomials for continuous functions (see \cite{M36}). Inequalities of the form \eqref{upp_est} are useful for estimating conditional numbers of the corresponding matrices when applying Galerkin's method.

The first estimates of the norms of submatrices of a general matrix were established by Kashin in \cite{Ka80b} and \cite{Ka80} in connection with A.~N.~Kol\-mo\-gorov's problem (still unsolved) on permutations of orthonormal systems. It was proved in \cite{Ka80b} (see {\bf M.\ref{Ka80b}} below) that every matrix $A$ with $\|A\|=1$ and big enough $M/N$ contains an $N\times N$ submatrix with small $(2,2)$--norm. At the first step of the proof the estimates of $(2,1)$--norms of random $2N\times N$ submatrices of $A$ were established. At the second step of the proof the Grothendieck factorization theorem on operators acting from $\ell_2^N$ to $\ell_1^{2N}$ was used. Sharp estimates in the mentioned above theorem on submatrices were obtained by Lunin \cite{Lun}. This direction was further developed in the series of
works by Bourgain and Tzafriri (see \cite{BT87}, \cite{BT89}), where the remarkable restricted invertibility theorem {\bf M.\ref{restr}, M.\ref{restr_}} (see the survey \cite{NY} about it) was proved: Let $A$ be an $N\times N$ matrix with ones on the diagonal, then there is a submatrix of $A$ whose rank is proportional to that of $A$ and which is well invertible in the sense that the norm of the inverse matrix depends only on the norm of $A$ and does not depend on the dimension $N$ (see {\bf M.\ref{restr_}}, {\bf M.\ref{restr_p}}).  A probabilistic selection based on a decoupling principle and a factorization argument were used here.  The discussed result has important applications in harmonic analysis and in the study of hilbertian systems (see \cite{BT89}). The method from the paper \cite{BT91} of averaging of the norms of submatrices of a given size (raised to the power of $\log N$) gave a partial solution of the old Kadison-Singer problem from quantum mechanics \cite{KS}, which has been finally solved by Marcus, Spielman, and Srivastava   in \cite{MSS}.

A new efficient method for finding both upper and lower bounds for the norms of submatrices was proposed by Batson, Spielman, and Srivastava in the paper \cite{BSS}. They considered positive semidefinite matrices $A$ of rank $N$ written as a sum of outer products
$$
A=\sum\limits_{i\leq M}v_iv_i^T,
$$
where the number of terms $M$ may be much larger than the rank $N$. Such a representation arises from any $N\times M$ matrix $B$ with columns $b_i$ by considering $A=BB^T=\sum_{i\leq M}b_ib_i^T$. It is shown in \cite{BSS} that for any $N\times M$ matrix $B$ and $\e>0$, there is a nonnegative diagonal $M\times M$ matrix $S$ with at most $\lceil N/\e^2\rceil$ nonzero entries such that $(1-\e)^2BB^T \preceq BSB^T\preceq (1+\e)^2 BB^T$. We comment on this result at the end of  Subsection \ref{greedy}.
 {The} recently developed method of interlacing polynomials by Marcus, Spielman, and Srivastava led to a number of important results and allowed, in particular, to obtain a solution to the Kadison-Singer problem in \cite{MSS} (see {\bf M.\ref{Weav}} for another formulation):
\begin{Theorem}\label{MSS} Let a system of vectors $\bv_1,\dots,\bv_M$ from $\bbC^N$ have the following properties: for all $\bw\in \bbC^N$
    \begin{equation}\label{A5c}
    \sum_{j=1}^M |\<\bw,\bv_j\>|^2 = \|\bw\|_2^2
    \end{equation}
    and for some $\e>0$
    \begin{equation*}
    \|\bv_j\|_2 \le \e,\qquad j=1,\dots,M.
    \end{equation*}
    Then there is a partition of $ \{1,2,\dots,M\}$ into two sets $S_1$ and $S_2$ such that for all $\bw\in\bbC^N$ and for each $i=1,2$
    \begin{equation*}
      \sum_{j\in S_i} |\<\bw,\bv_j\>|^2 \leq \frac{(1+\e\sqrt{2})^2}{2}\|\bw\|_2^2.
    \end{equation*}
    \end{Theorem}
    Note that condition \eqref{A5c} means that $\sum_{j\leq M}\bv_j\bar \bv_j^T$ is the $N\times N$ identity matrix.
 {The main application of Theorem \ref{MSS} in view of this survey is that it leads to the final in the sense of order Marcinkiewicz--type discretization theorem in $L_2$ (see {\bf D.15} and {\bf T.4}). Theorem \ref{MSS} has found many other applications. }For example, it is used for the estimation of the effectiveness of computational methods (\cite{Osw}), see {\bf M.\ref{Osw_}} below. In graph theory it is related to sparsification and used to find a small subgraph of a graph with similar spectral properties (see, e.g., \cite{Sr13} and \cite{Spi}).  The proof of existence of exponential frames in $L_2(S)$ for every unbounded set $S$ of finite measure in $\mathbb{R}$ (see \cite{NOU}) is also based on Theorem~\ref{MSS}. We refer the reader to the survey \cite{Bow}  on a recent progress on continuous frames inspired by the solution to the Kadison-Singer problem. The method of \cite{BSS} was used to construct one feature selection technique for linear Support Vector Machine, see \cite{Dr}.

The above discussed topic  is related to the results on the number of contact points of a convex body with its John ellipsoid (see, e.g. \cite{Rud97}, \cite{Versh01}, and \cite{Sr12} for algorithmic results). We have listed here only some possible applications of the results on the spectral and related properties of submatrices.

The scope of applications of results on the estimates of the operator norms and other properties of submatrices of a given matrix, considered in Section \ref{M}, is not limited to discretization problems. In turn, this topic belongs to a wider range of problems on the restriction of operators to coordinate subspaces. The problem here is formulated as follows. Let $X$, $Y$ be the Banach spaces, $T:X\rightarrow Y$ be a linear operator, $\{\varphi_j\}_{j\in\Omega}$ where $\Omega=\mathbb N$ (in the case of infinite-dimensional  $X$)
or $\Omega=\{1,\dots, N\}$, $N\in\mathbb N$ (for finite-dimensional  $X$) be a basis of $X$. We are looking for a subspace $L\subset X$ of the form $L=\sp\limits_{j\in \Gamma}\{\varphi_j\}$, $\Gamma\subset \Omega$, such that the restriction of the operator $T$ on the subspace $L$ has additional properties. The dual formulation of the problem uses a basis $\{f_j\}$ of $Y$. In that case we want to find $\Gamma\subset{\mathbb{N}}$ such that the operator
{$T|_{L_1}$, $L_1=\sp\limits_{j\in \Gamma}\{f_j\}$},
has additional properties.  As an example let us formulate the Banach theorem (see, e.g. \cite{KaczSt}): Let $\{\psi_j\}_{j=1}^{\infty}$ be an orthonormal system with the following property: for some $p>2$ and $D>0$ $\|\psi_j\|_p\leq D$, $j\in\mathbb N$. Then there exists an infinite subsystem $\{\psi_{\nu}\}_{\nu\in \Gamma}$ so that for all $\{a_{\nu}\}_{\nu\in\Gamma}\subset\bR$ the following inequality holds:

\begin{gather}\label{p_2}
\|\sum\limits_{\nu\in\Gamma}a_{\nu}\psi_{\nu}\|_{L_p}\leq C_{D, p}(\sum\limits_{\nu\in\Gamma}a_{\nu}^2)^{1/2}.
\end{gather}
In this case $X=\ell_2$, $Y=L_2$, $T:\ell_2\rightarrow L_2$, $T(\{a_j\}_{j=1}^{\infty})=\sum\limits_{j=1}^{\infty}a_j\psi_j$. Banach's theorem means that the operator $T$ acts between $\ell_2$ and $L_p$ on the coordinate subspace $\sp\{e_{\nu},\  \nu\in\Gamma\}$  of $\ell_2$. In the finite-dimensional case, the largest dimension of the corresponding good subspace is usually of interest. We recall, for example, a deep result by Bourgain \cite{BoLam}: Under the additional assumption $\|\psi_j\|_{L_{\infty}}\leq D$, $1\leq j\leq N$, we can choose a subsystem $\{\psi_j\}_{j\in\Gamma}$ from an orthonormal system $\{\psi_j\}_{j=1}^N$ with $|\Gamma|\geq N^{2/p}$ such that \eqref{p_2} holds. This theorem can be seen as a finite-dimensional analog of Banach's theorem above.

{\bf Connections with other areas.}
 In Section \ref{connect}
we discuss connections between the sampling discretization problem
and other areas of research and describe how
ideas and techniques from these areas can be applied
in the study of sampling discretization.

{\bf Moments of random vectors.}
 In Subsection \ref{MM}
we discuss the probabilistic approach to sampling discretization of the $L_q$ norms of
elements of an $N$-dimensional subspace $X_N$.
For $f\in X_N$ the quantity
$$
er(f,\xi,L_q):= \left|\|f\|_q^q - \frac{1}{m}\sum_{j=1}^m |f(\xi^j)|^q\right|,\quad \xi = (\xi^1,\dots,\xi^m) ,
$$
can be seen as a random variable on $\Omega^m$ with the measure $\mu^m$.
Let a system of functions $\{u_i(x)\}_{i=1}^N$, $x\in \Omega$, form a basis of $X_N$.
Consider the vector function $\bu(x) := (u_1(x),\dots,u_N(x))$,
which can be interpreted as a random vector on the probability space $(\Omega, \mu)$
and consider $m$ independent copies $\bu^j:=\bu(\xi^j)$
of this random vector.
We associate with $f=\sum_{i=1}^N y_i(f)u_i$ a vector of its coefficients $\by(f) = (y_1(f),\dots,y_N(f))$.
Then a class $\bF\subset X_N$ of functions can be described by the class
$K(\bF) := \{\by(f)\,:\, f\in \bF\} \subset \R^N$
of their coefficients
and vice versa for $K\subset \R^N$ define
$\bF(K):= \{f\,:\, \sum_{i=1}^N y_iu_i,\, \by\in K\}$.

For a set $K\subset \R^N$ we consider the following random variable
$$
V_q(K) :=
\sup_{\by \in K} \left| \frac{1}{m}\sum_{j=1}^m |\<\by,\bu^j\>|^q - \bbE|\<\by,\bu\>|^q\right|,
$$
where $\<\by,\bx\>=y_1x_1+\dots+y_Nx_N$ for $\bx,\by\in\R^N$ and $\bbE$ denotes the expectation of a random variable.
Clearly,
$$
V_q(K) = er(\bF(K),\xi,L_q) := \sup_{f\in \bF(K)}er(f,\xi,L_q) .
$$
In the sampling discretization we are interested in dependence of bounds for
$\inf\limits_{\xi} er(\bF(K),\xi,L_q)$ on the number $m$ of sampling points.
If now one obtains bounds for the expectation
$\bbE V_q(K)$ or for the probability of the event $\{V_q(K)>\varepsilon\}$,
then clearly as a corollary one gets a result about sampling discretization.
Thus, the following general problem arises:
How many independent copies of a random vector $\bu$
are needed to get $V_q(K)\le\varepsilon$ with high probability.

 Such a probabilistic reformulation
of the initial problem proved to be useful in many areas.
The first results on approximation of the moments of random vectors
were actually motivated
by the study of problems of complexity in computing volumes in high dimensions
(see \cite{B}, \cite{KLS}), by the study of restrictions of operators to coordinate subspaces
(see \cite{Rud}, where the problem proposed by B. Kashin and L. Tzafriri
in \cite{KaTz} was discussed), and by the study of the embedding
of finite-dimensional subspaces $X_N$ of $L_q$ into $\ell_q^m$ (see \cite{Ta2} and \cite{Ta3}).

{\bf Embedding of finite-dimensional subspaces.}
In Subsection \ref{Em} we discuss the problem of a good (almost isometrical) embedding
of an $N$-di\-men\-sional subspace $X_N$ of $L_q[0, 1]$
into an $N$-dimensional subspace of $\ell_q^m$ with the smallest possible $m$.
 This type of problem was originally motivated
by the study of the bounds for the dimension of an almost spherical section in the
Dvoretzky theorem (see \cite{FLM} and \cite{Ka77}, or Chapter 5 in \cite{A-AGM}).
In \cite{JS82} and \cite{Sche85} embeddings of general spaces into $\ell_1^m$ were studied.
Finally, in \cite{Sche3}, the described above problem appeared.
In this problem researchers do not restrict themselves to embedding by sampling
and are allowed to use any linear embedding.
Nevertheless, almost all approaches to this problem 
go through sampling (the empirical method).
In Subsection~\ref{Em} we describe several methods for solving
 the problem of embedding, which can be applied in sampling discretization.

{\bf Sparse approximation.} In Subsection \ref{greedy}
we give some remarks, which illustrate  {a} connection between sparse approximation and sampling discre\-ti\-zation of the $L_2$ norm. The properties of the subspace $X_N$ can be described in terms of a system $\U_N:=\{u_i\}_{i=1}^N$ of functions such that $X_N:=\sp\{u_i\}_{i=1}^N$.
In the case $X_N \subset L_2$ we assume that
the system $\U_N$ is orthonormal on $\Omega$ with respect to  {the} measure $\mu$. In the case of real functions   we associate with $x\in\Omega$ the matrix $G(x) := [u_i(x)u_j(x)]_{i,j=1}^N$. Clearly, $G(x)$ is a symmetric positive semi-definite matrix of rank $\leq 1$.
It is easy to see that for a set of points $\xi^k\in \Omega$, $k=1,\dots,m$, and $f=\sum_{i=1}^N b_iu_i$ we have
\be\label{Gb}
 \sum_{k=1}^m\la_k f(\xi^k)^2 - \int_\Omega f(x)^2 d\mu = {\mathbf b}^T\left(\sum_{k=1}^m \la_k G(\xi^k)-I\right){\mathbf b},
\ee
where ${\mathbf b} = (b_1,\dots,b_N)^T$ is the column vector and $I$ is the identity matrix. Therefore,
the sampling discretization problem is closely connected with a problem of approximation (representation) of the identity matrix $I$ by an $m$-term approximant with respect to the system $\{G(x)\}_{x\in\Omega}$. We can obtain some results on this approximation problem by greedy approximation methods, which provide a constructive way of finding good sampling points and by applying results from the theory of random matrices. We discuss these methods in more detail in Subsection \ref{greedy}.

{\bf Learning theory.} In Subsection \ref{learn} we demonstrate that the sampling discretization problem of the $L_2$ norm is closely related to
 supervised learning theory.
This is a vast area of research with a wide range of different settings. In Subsection \ref{learn} we only discuss a development of a setting from \cite{CS} (see \cite{VTbook}, Ch.4, for detailed discussion). In our further discussion we are interested in discretization of the $L_2$ norm of real functions from a given function class $\bF$. It is well known in learning theory that performance of the empirical risk minimization algorithm (least squares) can be controlled by asymptotic characteristics (entropy, Kolmogorov widths) measured in the uniform norm. Typically, analysis there is based on Bernstein's concentration measure inequality and chaining techniques.

Let $\Omega \subset\R^d$, $Y\subset \R$ be Borel sets, $\rho$ be a Borel probability measure on a Borel set $Z\subset \Omega \times Y$.  For $f:\Omega\to Y$ define {\it the error}
$$
\cE(f)  :=\int_Z(f(\bx)-y)^2 d\rho.
$$

Let   $\rho_\Omega$  be the marginal probability  measure of $\rho$ on $\Omega$, i.e.  for any Borel set $S\subset \Omega$ one has $\rho_\Omega (S) = \rho(S\times Y) $. Define
$$
f_\rho(\bx) := \bbE(y|\bx)
$$
to be a conditional expectation of $y$.
The function $f_\rho$ is known in statistics as the {\it regression function} of $\rho$. It is clear that if $f_\rho\in L_2(\rho_\Omega)$, then it minimizes the error $\cE(f)$ over all $f\in L_2(\rho_\Omega)$ i.e. $\cE(f_\rho)\le \cE(f)$. Thus, in the sense of error $\cE(\cdot)$ the regression function $f_\rho$ is the best to describe the relation between inputs $\bx\in \Omega$ and outputs $y\in Y$. The goal is to find an estimator $f_\bz$, on the base of given data $\bz:=((\bx^1,y_1),\dots,(\bx^m,y_m))$ that approximates $f_\rho$  well with high probability. We assume that $(\bx^i,y_i)$, $i=1,\dots,m$, are independent and distributed according to $\rho$.   We measure the error between $f_\bz$ and $f_\rho$ in the $L_2(\rho_\Omega)$ norm.

We note that a standard setting in the distribution-free theory of regression (see \cite{GKKW}) involves the expectation $\bbE(\|f_\rho- f_\bz\|_{L_2(\rho_\Omega)}^2)$ as a measure of quality of an estimator. An important new feature of  the setting in learning theory formulated in \cite{CS} (see \cite{VTbook} for detailed discussion) is the following.
 We begin with a class ${\mathcal M}$ of admissible measures $\rho$. Usually, we impose restrictions on $\rho$ in the form of restrictions on the regression function $f_\rho\in \Theta$, where $\Theta$ is a given class  of priors. Then the first step is to find an optimal estimator for  $\Theta$.  In regression theory we typically evaluate the performance of an estimator $f_\bz$  by studying its convergence in expectation, i.e. the rate of decay of the quantity $\bbE(\|f_\rho-f_\bz\|^2_{L_2(\rho_\Omega)})$ as the sample size $m$ increases. Here the expectation is taken with respect to the product measure $\rho^m$ defined on $Z^m$. We note that ${\mathcal E}(f_\bz)-{\mathcal E}(f_\rho) = \|f_\bz-f_\rho\|_{L_2(\rho_\Omega)}^2$.  A more accurate and more refined way of evaluating the performance of $f_\bz$ was proposed in \cite{CS}:   to
  systematically study the probability distribution function
$$
\rho^m\{\bz:\|f_\rho-f_\bz\|_{L_2(\rho_\Omega)}\ge\eta\}
$$
instead of the expectation.

We define the {\it empirical error} of $f$ as
$$
\E_\bz(f):= \frac{1}{m}\sum_{i=1}^m(f(\bx^i)-y_i)^2.
$$
Let $f\in L_2(\rho_\Omega)$. The {\it defect function} of $f$ is
$$
L_\bz(f) := L_{\bz,\rho}(f) := \E(f)-\E_\bz(f);\quad \bz=(z_1,\dots,z_m),\quad z_i=(\bx^i,y_i).
$$
We are interested in estimating the supremum of $|L_\bz(f)|$ over functions $f$ coming from a given class.

Settings for a supervised learning problem and a discretization problem are different.
In the supervised learning problem we are given a sample $\bz$ and we want to approximately recover the regression function $f_\rho$. It is important that we do not know
$\rho$. We only assume that we know that $f_\rho\in \Theta$. In the discretization of the $L_q$, $1\le q<\infty$, norm we assume that $f\in \bF$ and the probability measure $\mu$ is known. We want to find a discretization set $\xi=\{\bx^j\}_{j=1}^m$, which is good for the whole class $\bF$. However, the technique, based on the defect function, for solving the supervised learning problem can be applied to  the discretization problem.   Let us consider a given function class $\bF$ of real functions, defined on $\Omega$.
Suppose that the probability measure $\rho$ is such that $\rho_\Omega=\mu$ and for all $\bx \in \Omega$
we have $y=0$. In other words, we assume that $Y=\{0\}$. Then for the defect function we have
$$
L_\bz(f) =\int_\Omega f^2d\mu -\frac{1}{m}\sum_{j=1}^m f(\bx^j)^2,\quad |L_\bz(f)|= er(f,(\bx^1,\dots,\bx^m),L_2) 
$$
and
$$
\rho^m\{\bz:\sup_{f\in \bF} |L_\bz(f)|\ge \eta\} = \mu^m\{\bw  \in \Omega^m:\sup_{f\in \bF} er(f,\bw,L_2)\ge \eta\}.
$$
In Subsection \ref{learn} we illustrate this connection by highlighting some known results from learning theory and demonstrating how those results imply results on sampling discretization.

{The technique used there is a probabilistic technique. We derive sampling discretization results from 
the estimates of $\rho^m\{\bz:\sup_{f\in \bF}|L_\bz(f)|\le\e\}$. We stress that in order to prove existence of a point set that is good for the whole class we study the probability distribution of
 the supremum of the random quantity of interest over the class. There is other approach of using 
{\it randomized algorithms} (see, for instance, \cite{H}). In this approach we first take the expectation of the random quantity of interest for each individual function in the class and after that we take the supremum over the whole class. This approach is similar to the Monte Carlo approach for numerical integration discussed in Subsection \ref{learn}. However, as it is pointed out in 
Subsection \ref{learn}, it does not help to obtain results for the whole class, which contains infinitely many functions. It is an interesting and important setting but it is not directly related to our setting of looking for a good point set for the whole class. 
}

{\bf The structure of the paper.}
We now comment on the further structure of the paper. Section \ref{A} contains results on sampling discretization of integral norms $L_q$, $1\le q<\infty$, and the uniform norm $L_\infty$ of elements of finite-dimensional subspaces. In Subsection \ref{Atr} we present results on discretization of norms of trigonometric polynomials, in Subsection \ref{Ag} we
discuss general subspaces, and in Subsections \ref{Alg} and \ref{cap} -- algebraic polynomials. Mostly, these subsections contain known results, which we formulate in
paragraphs enumerated by letter {\bf D} with a number. Letter {\bf D} stands for
Discretization. In Subsection \ref{rec} we present a brief discussion of recent results on
sampling recovery proved with the aid of sampling discretization results. We enumerate the corresponding results by letter {\bf R} with a number. Letter {\bf R} stands for Recovery.
In Subsection \ref{tech} we complement the discussion on techniques from Section \ref{Int}
by more detailed discussion of important techniques used in sampling discretization.
We enumerate the corresponding results by letter {\bf T} with a number. Letter {\bf T} stands for Technique. Finally, in Subsection~\ref{OP} we formulate some open problems.
We note that the material of Section \ref{A} complements the earlier survey on sampling discretization \cite{DPTT}.

  In Section \ref{M} we give  a survey on results about the operator $(p,q)$--norms of submatrices which are connected to discretization. In addition to the classical case of the $(2,2)$-norm (spectral norm) we pay attention to the case of arbitrary $1\leq p,q\leq\infty$, which is not sufficiently covered in the literature.  

Section \ref{connect} is devoted to a detailed discussion of connections between settings and results from sampling discretization with settings and results from
other areas of mathematics. We show that the corresponding connections are very close.
It is a new and important contribution of this survey.
In Subsection \ref{MM} we demonstrate that the sampling discretization problem is very close to the corresponding problem on approximation
of moments of marginals of high-dimensional distributions.
We also show how several results about high-dimensional random vectors
can be applied to  sampling discretization.
In Subsection \ref{MM} we split the material into three large pieces enumerated by
{\bf MM.1} -- {\bf MM.3} and formulate results in terms of theorems and corollaries.
 In Subsection \ref{Em} we discuss the problem of good embedding of an $N$-dimensional
subspace $X_N$ of $L_q$ into $\ell_q^m$ with the smallest possible~$m$.
In this subsection we present the best known results
concerning the described problem ({\bf Em.1})
and
then we give some ideas, that are used in the proofs of these results
and can also be useful in the sampling discretization problem
({\bf Em.2} -- {\bf Em.6}).
In Subsection \ref{greedy} we give a brief comment on the connection between sampling discretization of the $L_2$ norm and the $m$-term approximation of the identity matrix with respect to a special dictionary, determined by the subspace $X_N$, of matrices of rank one.
In Subsection \ref{learn} we demonstrate how known results from learning theory can be used for establishing upper bounds for the sampling discretization error $er_m(\bF,L_2)$ for certain function classes $\bF$. We formulate results of this subsection in a form of theorems and corollaries.

Occasionally, for the reader's convenience we write $a_m\ll b_m$ instead of $a_m\le Cb_m$, where $C $ is a positive constant independent of $m$. In case $a_m\ll b_m$ and $b_m\ll a_m$ we write $a_m\asymp b_m$.

\section{Discretization}
\label{A}

  We give some general remarks before we proceed to technical details. We are interested in discretization of the $L_q$ norms, $1\le q\le \infty$, of elements of finite-dimensional subspaces. If a subspace of interest is a subspace of $\cC(\Omega)$, then all the functions from this subspace are defined at each point of $\Omega$. This is the case in many concrete situations, for instance, in the case of trigonometric and algebraic polynomials. However, we encounter a problem when we want to consider
arbitrary subspaces of $L_q$, $1\le q<\infty$. In this case by a function $f\in L_q(\Omega,\mu)$ we
understand a specific function (not an equivalency class), which is defined almost everywhere with respect to $\mu$ on $\Omega$. In other words, for $f\in L_q(\Omega,\mu)$ there exists a set $E(f)\subset \Omega$ such that $\mu(E(f))=0$ and $f(x)$ is defined for all $x\in \Omega \setminus E(f)$.
We say that a subspace $X_N\subset L_q(\Omega,\mu)$ is an $N$-dimensional subspace (more precisely, has dimension $\le N$) if there are $N$ functions $u_i\in X_N$, $i=1,\dots,N$, such that
$X_N=\sp\{u_1,\dots,u_N\}$. In this case, for a subspace $X_N$ there exists a set $E(X_N)\subset \Omega$ such that $\mu(E(X_N))=0$ and each $f\in X_N$ is defined for all $x\in \Omega \setminus E(X_N)$. 

{\bf The Marcinkiewicz discretization problem.} Let $\Omega$ be a compact subset of $\R^d$ with the probability measure $\mu$. We say that a linear subspace $X_N$ (index $N$ here, represents the dimension of $X_N$) of $L_q(\Omega,\mu)$, $1\le q < \infty$, admits the Marcinkiewicz-type discretization theorem with parameters $q$  and $m\in \N$ and positive constants $C_1\le C_2$ if there exists a set
$$
\Big\{\xi^j \in \Omega: j=1,\dots,m\Big\}
$$
 such that for any $f\in X_N$ we have
\be\label{A1}
C_1\|f\|_q^q \le \frac{1}{m} \sum_{j=1}^m |f(\xi^j)|^q \le C_2\|f\|_q^q.
\ee

{\bf The Bernstein discretization problem.} In the case $q=\infty$ we define $L_\infty$ as the space of continuous functions on $\Omega$  and ask for
\begin{equation*}
C_1\|f\|_\infty \le \max_{1\le j\le m} |f(\xi^j)| \le  \|f\|_\infty.
\end{equation*}

We will also use the following notation to express the above properties: The $\cM(m,q)$ (more precisely the $\cM(m,q,C_1,C_2)$) theorem holds for a subspace $X_N$, written $X_N \in \cM(m,q)$ (more precisely $X_N \in \cM(m,q,C_1,C_2)$). In the case $q=\infty$ we write $\cM(m,\infty,C_1)$.

{\bf The Marcinkiewicz discretization problem with weights.}  We say that a linear subspace $X_N$ of the $L_q(\Omega,\mu)$, $1\le q < \infty$, admits the weighted Marcinkiewicz-type discretization theorem with parameters $m$ and $q$ and positive constants $C_1\le C_2$ if there exist a set of knots $\{\xi^\nu \in \Omega\}$ and a set of weights $\{\la_\nu \in\bR\}$, $\nu=1,\dots,m$, such that for any $f\in X_N$ we have
\be\label{A3}
C_1\|f\|_q^q \le  \sum_{\nu=1}^m \la_\nu |f(\xi^\nu)|^q \le C_2\|f\|_q^q.
\ee
Then we also say that the $\cM^w(m,q)$ (more precisely the $\cM^w(m,q,C_1,C_2)$) theorem holds for  a subspace $X_N$ and write  $X_N\!\in\!\cM^w(m,q)$ (more precisely, $X_N\!\in\!\cM^w(m,q,C_1,C_2)$).
Obviously, $X_N\in\cM(m,q)$ implies that $X_N~\in~\cM^w(m,q)$.

{\bf The Marcinkiewicz discretization problem with $\e$.} For $1\le q<\infty$ we write $X_N\in \cM(m,q,\e)$ if (\ref{A1}) holds with $C_1=1-\e$ and $C_2=1+\e$.  Respectively,
we write $X_N\in \cM^w(m,q,\e)$ if (\ref{A3}) holds with $C_1=1-\e$ and $C_2=1+\e$.

In the case, when we need to specify either the set $\Omega$ or the measure $\mu$ we write in the notation $L_q(\Omega,\mu)$ instead of $q$, for instance, $\cM^w(m,L_q(\Omega,\mu))$ instead of $\cM^w(m,q)$.
We note that the most powerful results are 
when the $L_q$ norm of $f\in X_N$ is discretized exactly by the formula with equal weights $1/m$, i.e. the $\cM(m,q,0)$ theorems.

The optimal in the sense of order sampling discretization results are those, which guarantee $X_N \in \cM(CN,q,C_1,C_2)$. This is closely related to the following concept of {\it quasi-matrix systems} introduced in \cite{KaC2} in the case of $p=\infty$.

\begin{Definition}\label{ADqm} Let $p\in [1,\infty]$. A system of functions $\Phi=\{\varphi_j(x)\}_{j=1}^\infty$, $x\in\Omega$,
is called quasi-matrix if there exist positive constants $K_i$, $i=1,2,3$, such that for any $N\in\N$
there is a finite set $\{\xi^j\}_{j=1}^m\subset \Omega$ with $m\le K_1N$ with the property: For any
$f=\sum_{j=1}^Na_j \varphi_j$   the following inequalities hold in the case $1\le p<\infty$
$$
K_2\|f\|_p \le \left(\frac{1}{m}\sum_{j=1}^m|f(\xi^j)|^p\right)^{1/p} \le K_3\|f\|_p
$$
and for $p=\infty$
$$
K_2\|f\|_\infty \le \max_{1\le j\le m}|f(\xi^j)|.
$$
\end{Definition}

\subsection{Trigonometric polynomials}
\label{Atr}

In this subsection we specify $\Omega := \T^d := [0,2\pi)^d$ and $\mu$ to be the normalized Lebesgue measure $d\mu = (2\pi)^{-d}d\bx$.
By $Q$ we denote a finite subset of $\Z^d$, and $|Q|$ stands for the number of elements in $Q$. Let
$$
\Tr(Q):= \left\{f: f=\sum_{\bk\in Q}c_\bk e^{i\langle\bk,\bx\rangle},\  \  c_{\bk}\in\mathbb{C}\right\}.
$$

{\bf Polynomials with frequencies from parallelepipeds.} Consider $d~$-dimensional parallelepipeds
$$
\Pi(\mathbf N,d) :=\bigl \{\mathbf a\in \Z^d  : |a_j|\le N_j,\
j = 1,\dots,d \bigr\} ,\quad \bN=(N_1,\dots,N_d),
$$
where $N_j$ are nonnegative integers and the corresponding subspaces of the trigonometric polynomials
$$
\Tr(\bN,d) := \Tr(\Pi(\bN,d)).
$$
Then $\dim \Tr(\bN,d) =   \prod_{j=1}^d (2N_j  + 1)$.

{\bf D.1.} The following classical result on exact discretization of the $L_2$ norm is well known (see, for instance, \cite{VTbookMA}, p.7)
\begin{equation*}
\Tr(\bN,d) \in \cM(\dim \Tr(\bN,d),2,0).
\end{equation*}
 
{\bf D.2.} The following result was obtained by  Bernstein (see, for instance,  \cite{Z}, Ch.10, Theorem (7.28)) in the case $d=1$
\begin{equation*}
\Tr(N,1) \in \cM((1+\delta)\dim \Tr(N,1),\infty,C_1(\delta)) .
\end{equation*}

{\bf D.3.} The following result was obtained by  Marcinkiewicz in the case $d=1$, $1<q<\infty$, and by
Marcinkiewicz--Zygmund in the case $d=1$, $1\le q\le \infty$ (see \cite{Z}, Ch.10, \S7). For the multivariate analogs see \cite{VTbookMA}, p.102, Theorem 3.3.15.
\begin{equation*}
\Tr(\bN,d) \in \cM(C(d)\dim \Tr(\bN,d),q,C_1(d),C_2(d)),\quad 1\le q\le\infty.
\end{equation*}
The reader can find a generalization of the Marcinkiewicz--Zygmund results to the case of weighted $L_q$ spaces, namely, to the case of $L_q(\T,\mu)$, $d\mu = w(x)dx$, where $w$ is the Muckenhoupt weight, in \cite{KaLiz89}.

 {\bf Hyperbolic cross polynomials.} For $\bs\in\Z^d_+$
define
$$
\rho (\bs) := \{\bk \in \Z^d : [2^{s_j-1}] \le |k_j| < 2^{s_j}, \quad j=1,\dots,d\}
$$
where $[x]$ denotes the integer part of $x$. We define the step hyperbolic cross
$Q_n$ as follows
$$
Q_n := \bigcup_{\bs:\|\bs\|_1\le n} \rho(\bs)
$$
and the corresponding set of the hyperbolic cross polynomials as
$$
\Tr(Q_n) := \{f: f=\sum_{\bk\in Q_n} c_\bk e^{i\langle\bk,\bx\rangle}\}.
$$
In addition to the step hyperbolic cross $Q_n$ we also consider a more general step hyperbolic
cross $Q_n^\ga$, where $\ga = (\ga_1,\dots,\ga_d)$ has the form 
$$1=\ga_1=\dots=\ga_\nu< \ga_{\nu+1} \le \dots\le \ga_d$$
 with $\nu\in \N$, $\nu \le d$:
$$
Q_n^\ga := \bigcup_{\bs:(\ga,\bs)\le n} \rho(\bs),\quad (\ga,\bs):=\ga_1s_1+\dots+\ga_ds_d.
$$
It is clear that in the case $\ga = {\mathbf 1} := (1,\dots,1)$ we have $Q_n^{\mathbf 1} = Q_n$.
Note that $|Q_n^\ga| \asymp 2^nn^{\nu-1}$.

 {\bf D.4.} We begin with the case $q=1$. In this case we have (see \cite{VT180}): There are three positive constants
$C_i=C_i(\ga)$, $i=1,2,3$, such that we have
\begin{equation*}
\Tr(Q_n^\ga) \in \cM(m,1,C_1,C_2)\quad \text{provided}\quad m\ge C_3|Q_n^\ga|n^3 .
\end{equation*}
Weaker versions of the above result for $\ga = {\mathbf 1} $   were obtained in \cite{Bel} (with $n^3$ replaced by
$n^4$) and in \cite{VT158} (with $n^3$ replaced by $n^{7/2}$).

{\bf D.5.} We continue with the case $q=2$. In this case the problem is solved in the sense of order (see \cite{VT158}): There are three absolute positive constants
$C_i $, $i=1,2,3$, such that we have
\begin{equation*}
\Tr(Q_n^\ga) \in \cM(m,2,C_1,C_2)\quad \text{provided}\quad m\ge C_3|Q_n^\ga|  .
\end{equation*}
See {\bf D.9} and {\bf D.15} for further generalizations.

 {\bf D.6.} Consider the case $1<q<\infty$, $q\neq 2$. In this case we have (see \cite{VT180}): There are three positive constants
$C_i=C_i(q,\ga)$, $i=1,2,3$, such that we have
\begin{equation*}
\Tr(Q_n^\ga) \in \cM(m,q,C_1,C_2)\quad \text{provided}\quad m\ge C_3|Q_n^\ga|n^{w(\nu,q) },
\end{equation*}
 where
$$
 w(\nu,q) = 2, \quad q\in (1,2); \qquad
w(\nu,q) = (\nu-1)(q-2) + \min(q,3),\quad q\in (2,\infty).
$$
A weaker version of the above result for $\ga = {\mathbf 1} := (1,\dots,1)$ with $n^{w(\nu,q) }$ replaced by $n^{\max((d-1)(q-2),0)+4}$ was obtained in \cite{Bel}.

{\bf D.7.} Consider the case $q=\infty$. It turns out that in the case of the hyperbolic cross trigonometric polynomials sampling discretization results involving the $L_\infty$ norm are very different from those involving the $L_q$ norm, $1\le q<\infty$. The following statement was proved in \cite{KT3} -- \cite{KaTe03} for $d=2$: There exists a positive absolute constant $c$ with the property: If $m$ is such that
\begin{equation*}
\Tr(Q_n) \in \cM(m,\infty,C_1), \quad \text{then}\quad m\ge C(C_1) |Q_n|^{1+c} .
\end{equation*}

{\bf D.8.} The upper bound for $m$ is given in \cite{DPTT}. Define
$$\alpha_d:=\sum_{j=1}^d \frac 1j\qquad\mbox{ and}\qquad
\beta_d:=d-\alpha_d.
$$
Then there are two positive constants
$C_i=C_i(d)$, $i=1,2$, such that we have
\begin{equation*}
\Tr(Q_n) \in \cM(m,\infty,C_1)\quad \text{provided}\quad m\ge C_22^{n\alpha_d}n^{\beta_d}.
\end{equation*}

{\bf General trigonometric polynomials.} We now proceed to the case of general trigonometric subspaces of dimension $N$, in other words
subspaces $\Tr(Q)$ with $|Q|=N$.

{\bf D.9.} We begin with the case $q=2$. The following result was obtained in \cite{VT158}. There are three absolute positive constants
$C_i $, $i=1,2,3$, such that for any $Q\subset \Z^d$ we have
\begin{equation*}
\Tr(Q) \in \cM(m,2,C_1,C_2)\quad \text{provided}\quad m\ge C_3|Q|  .
\end{equation*}
See  {\bf D.15} for a further generalization.

{\bf D.10.} Consider the case $q\in [1,2)$. In this case we have (see \cite{VT180}): There are three positive constants
$C_i=C_i(q)$, $i=1,2,3$, such that we have for any $Q\in \Z^d$
$$
\Tr(Q) \in \cM(m,q,C_1,C_2)\quad \text{provided}\quad m\ge C_3|Q|(\log (2|Q|))^{w(q) },
$$
 where $w(1) = 3$ and $w(q) = 2$ for  $q\in (1,2).$
 See {\bf D.16} for more general results.

{\bf D.11.} We do not have results similar to {\bf D.10} in the case $q\in (2,\infty)$. We mention here a corollary of {\bf D.9} for the case of even $q$. Let $q=2s$, $s\in \N$. The case $s=1$ is covered by {\bf D.9}. We consider the case $s\ge 2$. For a given $Q\subset \Z^d$ consider the set
$$
Q^s:= \{\bk\,:\, \bk = \bk^1+\bk^2+\dots+\bk^s,\, \bk^j\in Q,\, j=1,\dots,s\}.
$$
Then, $|Q^s| \le |Q|^s$ and for any $f\in \Tr(Q)$ we have $f^s\in \Tr(Q^s)$. Therefore, discretization
of the $L_{2s}$ norm of elements of $\Tr(Q)$ follows from discretization of the $L_2$ norm of elements
of $\Tr(Q^s)$. Applying {\bf D.9} to $\Tr(Q^s)$ we obtain the following result: There are three absolute positive constants
$C_i $, $i=1,2,3$, such that for any $Q\subset \Z^d$ we have
\be\label{A11a}
\Tr(Q) \in \cM(m,2s,C_1,C_2)\quad \text{provided}\quad m\ge C_3|Q|^s.
\ee
 We also point out (see \cite[Corollary 4.5]{Kos}) that  in the case $2< q<\infty$ one has
\be\label{A11a'}
\Tr(Q) \in \cM(m, q , \varepsilon)\quad \text{provided}\quad m\ge C(q,\e)|Q|^{\frac{q}{2}}\log|Q|.
\ee
See  {\bf D.17} for more general results.

We point out that bound (\ref{A11a}) is sharp and bound (\ref{A11a'}) is sharp up to the
logarithmic factor $\log |Q|$. This follows from the corresponding necessary conditions discussed
in {\bf D.20. A lower bound}.

{\bf D.12.} In the case $q=\infty$ the discretization results are very different. We refer the reader to {\bf D.7} above and to the survey paper \cite{DPTT} for a discussion.   The following result is from \cite{DPTT} (see Theorem 6.7 there).
Let $\Lambda_n = \{k_j\}_{j=1}^n$ be a lacunary sequence: $k_1=1$, $k_{j+1} \ge bk_j$, $b>1$, $j=1,\dots,n-1$. Assume that a finite set $\xi=\{\xi^\nu\}_{\nu=1}^m\subset \mathbb T$ has
the following property
\begin{equation}\label{A11b}
\forall f\in \Tr(\Lambda_n) \qquad \|f\|_\infty \le K\max_{\nu}|f(\xi^\nu)|.
\end{equation}
Then
$$
m \ge (n/e)e^{Cn/K^2}
$$
with a constant $C>0$ which may only depend on $b$.

Thus, for successful discretization of the uniform norm of general trigonometric polynomials we may need the number of points to be exponential in the dimension of the subspace. Namely, for good discretization of the uniform norm of functions from $\Tr(\Lambda_n)$ (see (\ref{A11b})) we need exponentially many points. It turns out that the exponential number of points is always enough for discretization of the uniform norm (see \cite{KKT} and {\bf D.21} below).

\subsection{General subspaces}
\label{Ag}

We begin with the definition of entropy numbers.
  Let $X$ be a Banach space and let $B_X$ denote the closed unit ball of $X$ with center at $0$. Denote by $B_X(y,r)$ a ball with center $y$ and radius $r$: $B_X(y,r):=\{x\in X:\|x-y\|\le r\}$. For a compact set $A$ and a positive number $\e$ we define the covering number $N_\e(A)$
 as follows
$$
N_\e(A) := N_\e(A,X)
:=\min \{n : \exists y^1,\dots,y^n, y^j\in A :A\subseteq \cup_{j=1}^n B_X(y^j,\e)\}.
$$
It is convenient to consider the entropy numbers $\e_k(A,X)$:
$$
\e_k(A,X)  :=\inf \{\e : \exists y^1,\dots ,y^{2^k} \in A : A \subseteq \cup_{j=1}
^{2^k} B_X(y^j,\e)\}.
$$
In our definition of $N_\e(A)$ and $\e_k(A,X)$ we require $y^j\in A$. In a standard definition of $N_\e(A)$ and $\e_k(A,X)$ this restriction is not imposed.
However, it is well known that these entropy number characteristics may differ at most by a factor $2$. Namely, the following statement holds (see, for instance, \cite{VTbookMA}, p.322, Theorem 7.1.1 and \cite{VTbook}, p.208). Let $\e_k^*(A,X)$ be the entropy numbers without the restriction $y^1,\dots ,y^{2^k} \in A$. Then $\e_k^*(A,X)\le\e_k(A,X)\le 2\e_k^*(A,X)$. As above, we use the following notation for the unit $L_q$-ball of $X_N$
$$
X_N^q:= \{f\in X_N\subset L_q\,:\, \|f\|_q \le 1\}.
$$

{\bf Conditional results with entropy.} We formulate here some results on sampling discretization, which are proved under conditions on the entropy numbers of the unit balls of finite-dimensional subspaces of interest. In a spirit, these results are connected with known results in learning theory (see, for instance, \cite{VTbook}, Ch.4).

{\bf D.13.} In the general case $q\in [1,\infty)$  the following statement holds (see \cite{DPSTT1}):
 There exists a positive constant $C(q)$ such that for any $N$-dimensional subspace $X_N$ satisfying the condition
\be\label{A12}
\e_k(X^q_N,L_\infty) \le  B (N/k)^{1/q}, \quad 1\leq k\le N,
\ee
where $B\ge 1$, we have
$$
X_N \in \cM(m,q,\e) \quad \text{provided}\quad m \ge C(q,\e)NB^{q}(\log_2(2BN))^2.
$$
Note that in the case $q=1$, $\e=1/2$, the above statement was proved in \cite{VT159}.

{\bf D.14.} In the case $q\in[1,\infty)$ we have one more conditional result. In the above paragraph {\bf D.13} the condition is formulated in terms of the entropy numbers in the uniform norm $L_\infty$. Very recently, a new idea in this direction was developed in \cite{Kos}.  The corresponding   theorem with the conditions imposed on the entropy numbers in a weaker metric than the uniform norm was proved in \cite{Kos}. We now formulate that result. Let $Y_s:= \{y_j\}_{j=1}^s \subset \Omega$ be a set of sample points from the domain $\Omega$. Introduce a semi-norm
$$
\|f\|_{Y_s}:=\|f\|_{L_{\infty}(Y_s)} := \max_{1\le j\le s} |f(y_j)|.
$$
Clearly, for any $Y_s$ we have $\|f\|_{Y_s} \le \|f\|_\infty$.
The following result is from \cite{Kos} (see Corollary 3.8 there): There exists a number $C_1(q)>0$ such that for $m$ and $B\ge 1$ satisfying
$$
m\ge C_1(q) NB^q (\log N)^{w(q)},\quad w(1):=2, \quad w(q) := \max(q,2)-1, 1<q<\infty,
$$
and for a subspace $X_N$ satisfying the condition: For any set $Y_m\subset \Omega$
\be\label{A13}
\e_k(X^q_N,L_{\infty}(Y_m)) \le  B (N/k)^{1/q}, \quad 1\leq k\le N
\ee
we have
$$
X_N \in \cM(m,q),\quad 1\le q<\infty.
$$

 The entropy assumption with the discretized uniform norm
$\|f\|_{Y_s}$ was used to obtain upper bounds for the expectation of the supremum of a random process
$\sum\limits_{j=1}^m\varepsilon_j|f(y_j)|^q$, where $\varepsilon_j$ are 
i.i.d symmetric Bernoulli random variables. 
It would be interesting to understand whether one can relax the above entropy condition.
For example, in \cite{KaTz-lower} the following norm of a function was introduced
\begin{equation}\label{KT-norm}
\|f\|_{m, \infty} = \int_\Omega\ldots\int_\Omega\max\limits_{1\le j\le m}|f(y_j)|\, \mu(dy_1)\ldots\mu(dy_m).
\end{equation}
This norm provided lower bounds for the expectation 
of the supremum of a random process in a similar setting (see also \cite{Grig1} and \cite{Grig2}
for the development of ideas from \cite{KaTz-lower}). 
One may wonder: Is it possible to replace the entropy condition (\ref{A13})
for each subset $Y_s$
with a single entropy condition similar to (\ref{A13})
with respect to norm (\ref{KT-norm}) 
to get a good Marcinkiewicz-type discretization theorem?

{\bf Conditional results with the Nikol'skii inequality.}
 The following condition is widely used in discretization.

{\bf Condition E.} The orthonormal system $\{u_i(x)\}_{i=1}^N$ defined on $\Omega$ satisfies Condition E with a constant $t$ if for all $x\in \Omega$
$$
 \sum_{i=1}^N |u_i(x)|^2 \le Nt^2.
$$

It is well known that Condition E is equivalent to the following Nikol'skii inequality (see, for instance, \cite{LT} for an explanation and \cite{DP} for a detailed discussion). We say that $X_N$ satisfies the Nikol'skii inequality for the pair $(2,\infty)$ if there exists a constant $t$ such that
\begin{equation}\label{A14}
\|f\|_\infty \leq t N^{\frac12}\|f\|_2,\   \ \forall f\in X_N.
\end{equation}
It is known that (\ref{A14}) implies Condition E for any orthonormal basis of $X_N$ and that Condition E for an orthonormal basis of $X_N$ implies (\ref{A14}). {Also, Condition E is closely related to the Christoffel function of the subspace $X_N$ (see {\bf D.26} below). }

{\bf D.15.} Consider the case $q=2$. In this case the sampling discretization problem for subspaces satisfying Condition E is solved in the sense of order (see \cite{LT}): There are three absolute positive constants $C_i $, $i=1,2,3$, such that for any $N$-dimensional subspace $X_N \subset L_\infty(\Omega)$   satisfying (\ref{A14})   we have
\be\label{A15}
X_N \in \cM(m,2,C_1,C_2t^2)\quad \text{provided}\quad m\ge C_3t^2N  .
\ee

The above result has history beginning with the Rudelson's result from \cite{Rud}. In the paper \cite{Rud} it is formulated in terms of submatrices of an orthogonal matrix  (see {\bf M.\ref{Rud_matrices}} below). We reformulate it in our notation. Let $\Omega_M~=~\{x^j\}_{j=1}^M$ be a discrete set with the probability measure $\mu(x^j)=1/M$, $j=1,\dots,M$. Suppose that $X_N\subset L_2(\Omega_M,\mu)$ satisfies (\ref{A14}). Then for $\e \in (0,1)$ we have
\be\label{A16}
X_N \in \cM(m,2,\e)\quad \text{provided}\quad  m\ge C\frac{t^2}{\e^2}N\log\frac{Nt^2}{\e^2}.
\ee
In the followup paper \cite{Rud2} more general results were obtained. In particular, these more general results imply $X_N \in \cM(m,2,\e)$ 
provided that $m~\ge ~C(t/\e)^2N\log N$ for a general set $\Omega$
(see a discussion after Corollary~\ref{CRud} below).
In \cite{VT159}  it was demonstrated how the Bernstein-type concentration of measure inequalities for random matrices can be used for proving an analog of the above Rudelson's result for a general $\Omega$. The proof in \cite{VT159} is based on a different idea than Rudelson's proof.

Both (\ref{A16}) and its version for a general $\Omega$ provide sufficient conditions on $m$   for existence of a good set of cardinality $m$ for sampling discretization. These sufficient conditions are close to the necessary condition, which is $m\ge N$, but still have an extra $\log N$ factor in the bound for $m$.  {The} result in (\ref{A15}) gives a sufficient condition on $m$ without an extra $\log N$ factor.   The first result {of that type} was obtained in \cite{VT158} where  (\ref{A15}) was proved for special domains $\Omega_M$ and under a stronger condition than Condition E: $\sum_{i=1}^N |u_i(x^j)|^2 = N$, $j=1,\dots,M$.

{\bf D.16.} Consider the case $q\in [1,2)$. The following result was obtained in \cite{DPSTT2}:
Let $X_N$ be an $N$-dimensional subspace of $L_\infty(\Omega)$ satisfying (\ref{A14}) with  $\log t\le \alpha \log N$.
Then we have
\be\label{A17}
X_N \in \cM(m,q,\e)\quad \text{provided}\quad m\ge C(q,\alpha,\e)tN(\log N)^3.
\ee
In the case $q\in (1,2)$ result (\ref{A17}) was improved in
 \cite{Kos} by replacing $(\log N)^3$ by $(\log N)^2$. It is interesting to note that in results for $q\in [1,2]$ the condition on $X_N$ -- the Nikol'skii inequality for the pair $(2,\infty)$ -- does not depend on $q$.

Let $p\in [1,\infty)$ and $X_N\subset L_\infty(\Omega)$. The inequality
\begin{equation*}
\|f\|_\infty \leq M\|f\|_p,\   \ \forall f\in X_N
\end{equation*}
is called the Nikol'skii inequality for the pair $(p,\infty)$ with the constant $M$.
Denote
$$
M_p(X_N):= \sup_{f\in X_N;f\neq 0} \|f\|_\infty/\|f\|_p.
$$

{\bf D.17.} The understood case $q\in(2,\infty)$ is somewhat different from the case  $q\in[1,2)$. In the former we use the Nikol'skii inequality for the pair $(q,\infty)$. The following results are known. It was proved in \cite{Kos}
  that the Nikol'skii inequality for the pair $(q,\infty)$ with the constant $M_q(X_N)$ of the order $N^{1/q}$ implies
$$
X_N \in \cM(m,q) \quad \text{provided} \quad m \gg N(\log N)^q.
$$ 
It is pointed out in \cite{VT180} that the Nikol'skii inequality for the pair $(q,\infty)$ with the constant $M_q(X_N)$ of the order $N^{1/q}$
combined with the assumption that $X_N \in \cM(s,\infty)$ with $s\le aN^c$ imply
$$
X_N \in \cM(m,q) \quad \text{provided} \quad m \gg N(\log N)^3.
$$
We also mention  (see \cite[Corollary 4.5]{Kos}) that, for subspaces $X_N$
satisfying Condition E,
one has
\begin{equation*}
X_N \in \cM(m, q)\quad \text{provided}\quad m\gg CN^{\frac{q}{2}}\log N.
\end{equation*}
For further results see Corollaries \ref{C-ma-1}--\ref{C-GR}.

{\bf Comment.} We discussed above conditional results for sampling discretization -- {\bf D.~13-14} under the entropy condition and {\bf D.~15-17} under the Nikol'skii inequality condition. These two
conditions are related. We will return to this connection at the end of Subsection \ref{tech}.

 {\bf Unconditional results.}
In the case of weighted discretization we can prove Marcinkiewicz-type discretization theorems for an arbitrary finite-di\-men\-si\-onal subspace.

{\bf D.18.} In the case $q=2$ we begin with he following result from \cite{LT}: There are three absolute positive constants $C_i $, $i=1,2,3$, such that for any (real or complex) $N$-dimensional subspace $X_N \subset L_2(\Omega,\mu)$   we have
\be\label{A19}
X_N \in \cM^w(m,2,C_1,C_2)\quad \text{provided}\quad m\ge C_3N  .
\ee
Moreover, the weights can be made positive.

The first result in that direction was obtained in  \cite{BSS} (see Theorem 3.1 there). In the case of a special domain $\Omega_M$ and real subspaces, the authors proved (\ref{A19}) with $C_3=b$, where $b>1$ is any number, $C_1=1$, and $C_2=\left(\frac{\sqrt{b}+1}{\sqrt{b}-1}\right)^2$.
Further, it was observed  in \cite[Theorem 2.13]{DPTT} that  this last result from \cite{BSS} with  a general probability space $(\Omega, \mu)$ in place of the discrete space  $(\Omega_M, \mu)$ remains true (with other constant $C_2$) if $X_N\subset L_4(\Omega,\mu)$.  It was proved in \cite{DPSTT2} (see Theorem 6.3 there) that the additional assumption ${X_N\subset L_4(\Omega,\mu)}$ can be dropped.

{\bf D.19.} Consider the case $q\in [1,2)$. The following result was proved in \cite{DPSTT2}: For given
$q\in [1,2)$ and $\e\in (0,1)$ there exists $C(q,\e)$ such that for any $N$-dimensional subspace of $L_q(\Omega,\mu)$ we have
\begin{equation*}
X_N \in \cM^w(m,q,\e)\quad \text{provided}\quad m\ge C(q,\e)N(\log N)^3  .
\end{equation*}

{\bf D.20.~An upper bound.} We do not have results similar to paragraph {\bf D.19} in the case $q\in (2,\infty)$, but we have a corollary of {\bf D.18} for the case of even $q$. Let $q=2s$, $s\in \N$. The case $s=1$ is covered by {\bf D.18}. We consider the case $s\ge 2$. For a given $N$-dimensional subspace $X_N \subset L_q(\Omega,\mu)$ consider
$$
(X_N)^s:= \{f\,:\, f = f_1\times f_2\times\dots\times f_s,\, f_j\in X_N,\, j=1,\dots,s\}.
$$
It is well known and easy to see that $(X_N)^s \subset L_2(\Omega,\mu)$ and
$\dim (X_N)^s \le N^s$. For any $f\in X_N$ we have $f^s\in (X_N)^s$. Therefore, discretization
of the $L_{2s}$ norm of elements of $X_N$ follows from discretization of the $L_2$ norm of elements
of $(X_N)^s$. Applying {\bf D.18} to $(X_N)^s$ we obtain the following result. There are three absolute positive constants $C_i $, $i=1,2,3$, such that for any (real or complex) $N$-dimensional subspace $X_N \subset L_{2s}(\Omega,\mu)$   we have
\begin{equation*}
X_N \in \cM^w(m,2s,C_1,C_2)\quad \text{provided}\quad m\ge C_3N^s  .
\end{equation*}
Moreover, the weights can be made positive.

  Actually, in the case $2<q<\infty$
   one has (see the discussion after Corollary~4.5 in {\bf MM.3})
\begin{equation*}
X_N \in \cM^w(m, q , \varepsilon)\quad \text{provided}\quad m\ge CN^{\frac{q}{2}}\log N.
\end{equation*}

{\bf D.20.~A lower bound.}
We point out that, when the $L_q$ norm with $q\in (2,\infty)$ on a subspace $X_N$ is equivalent to the $L_2$
norm, it is necessary to have at least $cN^{q/2}$ points
for discretization with positive weights (see e.g. \cite{BDGJN}).
Indeed, assume that $\|f\|_q\le M\|f\|_2$ $\forall f\in X_N$ and for some
points $\xi^1,\ldots, \xi^m$ and for some positive weights $\lambda_1, \ldots, \lambda_m$
one has
$$
c\|f\|_q^q\le \sum\limits_{\nu=1}^{m}\lambda_\nu|f(\xi^\nu)|^q\le C\|f\|_q^q,
$$
where $C>c>0$ are some fixed constants.
Let $u_1,\ldots, u_N$ be an orthonormal basis in $X_N$.
Then for each choice of signs $\varepsilon=(\varepsilon_1,\ldots, \varepsilon_N)$
for the function $f_\varepsilon:=\sum\limits_{k=1}^{N}\varepsilon_ku_k$
one has
$$
c N^{q/2}=c\|f_\varepsilon\|_2^{q}\le c\|f_\varepsilon\|_q^q
\le \sum\limits_{\nu=1}^{m}\lambda_\nu|f_\varepsilon(\xi^\nu)|^q=
\sum\limits_{\nu=1}^{m}\lambda_\nu\Bigl|\sum\limits_{k=1}^N\varepsilon_ku_k(\xi^\nu)\Bigr|^q.
$$
We now assume that signs are chosen randomly and take the average:
$$
c N^{q/2}
\le
\sum\limits_{\nu=1}^{m}\lambda_\nu\mathbb{E}_\varepsilon\Bigl|\sum\limits_{k=1}^N\varepsilon_ku_k(\xi^\nu)\Bigr|^q
\le B_q\sum\limits_{\nu=1}^{m}\lambda_\nu\Bigl(\sum\limits_{k=1}^N|u_k(\xi^\nu)|^2\Bigr)^{q/2},
$$
where we have applied the Khintchine inequality.
We note that
$$
\lambda_\nu\Bigl(\sum\limits_{k=1}^N|u_k(\xi^\nu)|^2\Bigr)^{q/2} =
\sup\limits_{f\in X_N, \|f\|_2\le 1}\lambda_\nu |f(\xi^\nu)|^q\le CM^q,
$$
since
$$
\lambda_\nu |f(\xi^\nu)|^q\le \sum\limits_{\nu=1}^{m}\lambda_\nu|f(\xi^\nu)|^q\le
C\|f\|_q^q\le CM^q\|f\|_2^q.
$$
Thus, $c N^{q/2}\le B_qCM^qm$ and $m\ge \frac{c}{CB_qM^q} N^{q/2}$.

A classical example of a subspace with equivalent $L_q$, $1\leq q<\infty$, and $L_2$ norms is
a subspace $\Tr(\Lambda_n)$, where $\Lambda_n = \{k_j\}_{j=1}^n$ is a lacunary sequence: $k_1=1$, $k_{j+1} \ge bk_j$, $b>1$, $j=1,\dots,n-1$.

As another example of a space with equivalent $L_q$ and $L_2$ norms, consider the space
$X_N:=\{\langle\cdot, \by\rangle\colon \by\in \mathbb{R}^N\}\subset L^q(Q_N)$,  where $Q_N:=[-\frac{1}{2}, \frac{1}{2}]^N$ -- is the unit cube endowed with
the standard Lebesgue measure.
It is known {(see Theorem 2.4.6 in \cite{BGVV})} that for each $q$ there is a constant $c_q$,
dependent only on $q$ and, in particular, independent of the dimension $N$,
such that
$\|f\|_q\le c_q\|f\|_2$ $\forall f\in X_N$.
Moreover, $X_N$ satisfies the $(2, \infty)$ Nikol'skii-type inequality assumption
with constant $3$ i.e. $\|f\|_\infty\le 3N^{1/2}\|f\|_2$ $\forall f\in X_N$.

{\bf D.21.}  {The} recent paper  \cite{KKT} addresses sampling discretization of the uniform norm. In particular, it is proved there that for discretization of the uniform norm of elements of any $N$-dimensional subspace of $L_\infty(\Omega)= \mathcal C(\Omega)$ it is sufficient to use $e^{CN}$
sample points. The following statement is proved in   \cite{KKT}.
Let $X_N$ be an $N$-dimensional subspace of $L_\infty(\Omega)$.
There exists a set $\xi=\{\xi^\nu\}_{\nu=1}^m$ of $m\le 9^N$ points such that for any $f\in X_N$ we
have
\begin{equation*} 
\|f\|_\infty \le 2\max_\nu |f(\xi^\nu)|\quad \text{which means}\quad X_N \in \cM(9^N,\infty,1/2).
\end{equation*}

{We now mention some results, where{, in contrast to our setting,} the constant $C_1$ in the Bernstein discretization problem is allowed to depend on $N$.
The following result was obtained in \cite{NoLN} (see Proposition 1.2.3 there).
 Let $X_N$ be an $N$-dimensional subspace of $\cC(\Omega)$ and let $\e>0$. 
There exists a set $\xi=\{\xi^\nu\}_{\nu=1}^N$ of $N$ points such that for any $f\in X_N$ we 
have
$$
\|f\|_\infty \le (N+\e)\max_\nu |f(\xi^\nu)|.
$$
}

{The following conditional result, which connects the upper bound in the discretization theorem for the uniform norm with the Nikol'skii-type inequality between $\cC$ and $L_2$ norms, was proved in \cite{DPTT}.
 Let  $\Omega := [0,1]^d$. Assume that a real $N$-dimensional subspace $Y_N\subset \cC(\Omega)$ satisfies the Nikol'skii-type inequality: For any $f\in Y_N$
 $$
 \|f\|_\infty \le H(N)\|f\|_2,\quad \|f\|_2 := \left(\int_\Omega |f(\bx)|^2d\mu\right)^{1/2},
 $$
 where $\mu$ is the Lebesgue measure on $\Omega$.
 Then for any $a>1$ there exists a set $\xi(m)=\{\xi^j\}_{j=1}^m\subset \Omega$ with the property:  $m \le a N$ and
 for any $f\in Y_N$ we have  
$$
 \|f\|_\infty \le C(a)H(N)\max_{1\le j\le m} |f(\xi^j)|, 
$$
where $C(a)$ is a  positive constant.
}

\subsection{Real algebraic polynomials}
\label{Alg}

In this subsection we discuss the case $\Omega =[-1,1]^d$ with a probability measure $\mu$ such that $d\mu=w(\bx)d\bx$. A subspace $X_N$ will be a subspace of algebraic polynomials.

{\bf D.22.} We begin with direct corollaries of results on the trigonometric polynomials. For $\bx\in [-1,1]^d$ and $\btt\in [0,\pi]^d$ consider the following change of variables $x_j=-\cos t_j$, $j=1,\dots,d$.
Also, consider the measure $\mu_c$ with density $w_c(\bx):= \pi^{-d} \prod_{j=1}^d (1-x_j^2)^{-1/2}$ -- the Chebyshev measure -- on $[-1,1]^d$. Then we have
$$
\int_{[-1,1]^d}|f(\bx)|^qw_c(\bx)d\bx = \pi^{-d}\int_{[0,\pi]^d}|f(-\cos t_1,\dots,-\cos t_d)|^qd\btt.
$$
Thus, result {\bf D.3} implies the following relation. 
{Let $\bN~=~(N_1,\dots,N_d)$ and $\cP(\bN,d)$ be
  the subspace of algebraic polynomials of degree $N_j$ in the variable $j$, $j=1,\dots,d$. }Then we have
 for $1\le q<  \infty$
\begin{equation*}
\cP(\bN,d) \in \cM(C(d)\dim \cP(\bN,d), L_q([-1,1]^d,\mu_c)),\quad d\mu_c = w_c(\bx)d\bx,
\end{equation*}
and
\begin{equation*}
\cP(\bN,d) \in \cM(C(d)\dim \cP(\bN,d), L_\infty([-1,1]^d)).
\end{equation*}

{\bf D.23.} Consider special algebraic polynomials -- Chebyshev polynomials
$$
T_k(x) := \frac{1}{2}\left((x+(x^2-1)^{1/2})^k + (x-(x^2-1)^{1/2})^k\right), \quad k=0,1,\dots.
$$
The system $\{T_k(x)\}_{k=0}^\infty$ forms an orthonormal basis in $L_2([-1,1],\mu_c)$. Note that
$T_k(-\cos t)=(-1)^k\cos kt$. Denote
$$
T_\bk(\bx) := \prod_{j=1}^d T_{k_j}(x_j),\quad \bk=(k_1,\dots,k_d),\quad \bx=(x_1,\dots,x_d).
$$
For a finite set $Q\subset \Z_+^d$ denote
$$
\cC(Q,d) := \left\{f\,:\, f(\bx) =\sum_{\bk\in Q} c_\bk T_\bk(\bx)\right\}.
$$
Then, result {\bf D.15} implies for any $Q\subset \Z_+^d$
\be\label{Al3}
\cC(Q,d) \in \cM(m,L_2([-1,1]^d,\mu_c),C_1,C_2)\quad \text{provided}\quad m\ge C_3|Q|  .
\ee
We now consider sets $Q\subset \Z^d_+$ with the property: if ${\bk= (k_1,\dots,k_d)\in Q}$, then the box
$\prod_{j=1}^d [0,k_j]$ also belongs to $Q$. These sets were used in the multivariate approximation of periodic functions in \cite{Tem4} (see also \cite{VTbookMA}, p.141). There, these sets were identified by saying: "$Q$ has the property $S$". Later, sets with this property
have become important in Stochastic PDEs, where these sets are called {\it lower sets}.

For a finite set $Q\subset \Z_+^d$ denote
$$
\cP(Q,d) := \left\{f\,:\, f(\bx) =\sum_{\bk\in Q} c_\bk \bx^\bk\right\},\quad \bx^\bk=\prod_{j=1}^d x_j^{k_j}.
$$
Note that for any lower set $Q$ we have $\cP(Q,d)=\cC(Q,d)$. Therefore, (\ref{Al3}) implies: For any lower set $Q$ we have
\begin{equation*}
\cP(Q,d) \in \cM(m,L_2([-1,1]^d,\mu_c),C_1,C_2)\quad \text{provided}\quad m\ge C_3|Q|  .
\end{equation*}

{\bf D.24.} We formulate a result from \cite{KalTz}. Consider the case $d=1$ and the measure $\mu(\gamma)$ with density $w(x,\gamma)= (1-x^2)^{\gamma -1/2}$, $-1/2<\gamma <\infty$. Then we have
\begin{equation*}
\cP(n-1,1) \in \cM^w(n,L_q([-1,1],\mu(\gamma)),\quad 1<q<\infty.
\end{equation*}
Note that in the case $\ga=1/2$ the measure $\mu(1/2)$ is the standard Lebesgue measure.
The reader can find some discretization results for this case in the book \cite{Shar}. In particular, it is known {(see \cite{Riv} and \cite{Shar1})} that if we restrict ourselves to uniformly distributed points $\xi =\{\xi^j\}_{j=0}^{2m}$, $\xi^j = (j-m)/m$, then for
sampling discretization of the $L_q([-1,1],\mu(1/2))$ norm of polynomials from $\cP(n,1)$ it is necessary and sufficient that $m$ grows as $n^2$ with $n \to \infty$.

{\bf D.25.} There are some recent papers on discretization in the uniform norm $L_\infty$ of subspaces of algebraic polynomials defined on a   set $\Omega \subset \R^d$. Denote by $\cP^d_n$ the set of real algebraic polynomials of total degree $\le\!n$. In other words, $\cP^d_n$ is the set of real polynomials from $\cP(\Delta_n,d)$, where ${\Delta_n := \{\bk\in \Z_+^d\,:\, \|\bk\|_1 := k_1+\dots+k_d \le n\}}$, {the cardinality
$|\cP^d_n|=\binom{n+d-1}{d-1}$, $|\cP^d_n|^{\frac{1}{d-1}}\asymp\frac{n+d-1}{d-1}.$} The following result was obtained in \cite{Kr1}: Let $d=2$ and let $\Omega \subset \R^2$ be a convex compact set with nonempty interior, then for $\e\in (0,1)$
\begin{equation*}
\cP^2_n \in \cM(4\cdot 10^5(n/\e)^2,L_\infty(\Omega),(1+\e)^{-1}).
\end{equation*}
The following result is from \cite{BBCL}: Let $\Omega \subset \R^d$, $d\in \N$, be a compact set with nonempty interior, then there exist positive  constants $C_0(d)$ and $C_1(d)$ such that for any $m\ge C_0(d)(n\log n)^d$ we have
\begin{equation*}
\cP^d_n \in \cM(m,L_\infty(\Omega),C_1(d)).
\end{equation*}

{\bf D.26.} We continue a discussion from {\bf D.25} and present a conditional result from \cite{Pr1}. Let
$\{P_j\}_{j=1}^N$ be an orthonormal basis of $\cP^d_n$ with respect to the normalized Lebesgue measure on $\Omega$. Define the Christoffel function for the $\cP^d_n$ as
$$
\lambda_n(\bx,\Omega):= \left(\sum_{j=1}^N P_j(\bx)^2\right)^{-1/2}.
$$
The following conditional result is proved in \cite{Pr1} with the aid of a lemma from \cite{BV}: Suppose
$\Omega$ is such that
$$
\lambda_n(\bx,\Omega)\le C(\Omega)\lambda_{2n}(\bx,\Omega),\quad \bx\in\Omega.
$$
Then
\begin{equation*}
\cP^d_n \in \cM(C(\Omega)\dim(\cP^d_n),\infty).
\end{equation*}

{\bf D.27.} We now formulate one recent result for the case $1\le q\le \infty$.  The following result was obtained in \cite{DaP}: Let $\Omega$ be a bounded $C^2$-domain in $\R^d$ with the normalized Lebesgue measure $\mu_L$ on it. Then for $1\le q<\infty$
\begin{equation*}
\cP^d_n \in \cM^w(C(\Omega)\dim(\cP^d_n),L_q(\Omega,\mu_L),1/2,3/2)
\end{equation*}
and for $q=\infty$
\begin{equation*}
\cP^d_n \in \cM(C(\Omega)\dim(\cP^d_n),L_\infty(\Omega)).
\end{equation*}

{\bf A comment on other settings.} The above discussion in this section shows that there is only one case, namely the case $q=2$, where the weighted sampling discretization problem is solved (in the sense of order). This result is presented in {\bf D.18}. We point out that there are many ingredients in the setting of the sampling discretization problem. Here they are: (1) $L_q(\Omega,\mu)$, which contains three ingredients: $1\le q\le \infty$, $\Omega$, and {the} measure $\mu$; (2) {the} subspace $X_N$; (3) restrictions on the weights $\la_1,\dots,\la_m$ in the discretization formula (\ref{A3}). Clearly, a result like {\bf D.18} with $1\le q<\infty$ and equal weights will solve the problem of discretization for all $1\le q<\infty$. Note, that results of {\bf D.7} show that there is no analog of {\bf D.18} for $q=\infty$. However, there is no known result at the level of generality of {\bf D.18} for $1\le q<\infty$ distinct from $2$. Therefore, in order to obtain results on discretization we impose some restrictions on the above mentioned ingredients of the problem. Typical restrictions are as follows. Instead of a general subspace $X_N$ we study a specific one, say, trigonometric or algebraic polynomials. Instead of a general domain $\Omega$ we study a specific one, say, a cube or a convex domain. Instead of a general measure $\mu$ we study a specific one, say, the Lebesgue or the Chebyshev measure.
We now comment on a special case, when $X_N$ is the $\cP^d_n$. We discuss different sets $\Omega$. Let us begin with the case $d=2$. Consider a unit circle $S^1$ as the $\Omega$. Then it is clear that the subspace $\cP^2_n$ restricted to $S^1$ becomes the subspace of the trigonometric polynomials (with the angle from the polar coordinates as the argument) of degree $n$. There is a number of papers, where $\cP^d_n$ is restricted to the sphere $S^{d-1}$ or to a part of it. We briefly mention here some of these papers and refer the reader to these papers for  further discussion.
 The authors of \cite{MT} proved the weighted Marcinkiewicz discretization theorems with general doubling measures for algebraic polynomials on intervals and for  trigonometric polynomials on circles, with the number of required sample points comparable to the dimension of the corresponding polynomial spaces.
 The paper \cite{Da1} extends the results of \cite{MT}  to the sphere and the ball. Weighted  Marcinkiewicz discretization theorems with doubling measures were established in this paper for spherical polynomials on the unit sphere and for algebraic polynomials on the unit ball. Further, the weighted Marcinkiewicz discretization theorems for spherical polynomials on spherical caps were established in \cite{DW}. In both \cite{Da1} and \cite{DW}, the number of sampling points required in the Marcinkiewicz discretization theorems is asymptotically equivalent to the dimension of the corresponding polynomial spaces.  

\subsection{Complex algebraic polynomials}
\label{cap}

 {The} discretization theorems of Bernstein and Marcinkiewicz--Zygmund (see {\bf D.2} and {\bf D.3}) were obtained as natural corollaries of the Bernstein and Zygmund inequalities for derivatives of the trigonometric polynomials. The corresponding inequalities for derivatives have been established and used in the proofs of the inverse theorems of approximation -- i.e. statements on smoothness properties of a function, which is approximated well by trigonometric polynomials. Applications of these inequalities
to discretization problems appeared two decades later, despite the fact that an approach to construction of economical nets, which provide good sampling discretization of a function, based on information of its derivative, is absolutely natural. We demonstrate this remark on the example of
a subspace $X_N \subset \cC(\Omega)$, $\dim X_N=N$, consisting of smooth functions defined on a compact subset $\Omega\subset \R^d$. Let $B_N:= X_N^\infty$ be the unit $L_\infty$-ball in $X_N$. Suppose that for the gradient $\nabla f$ we have the following bound for each $f\in B_N$
$$
|\nabla f(\bx)| \le M,\quad \bx\in \Omega,
$$
where, as usual, ${|\by|\!=\!\|\by\|_{\ell^d_2}\!=\!
 (|y_1|^2\!+\!\dots\!+\!|y_d|^2)^{1/2}}$ for $\by=(y_1,\dots,y_d) \in \R^d$.
 Then for any finite subset $\La\subset \Omega$, which forms an $(2M)^{-1}$-net
for the $\Omega$ in the metric $\ell^d_2$, and for each $f\in B_N$ the following inequality holds
$$
\max_{\xi\in\La} |f(\xi)| \le \|f\|_{\cC(\Omega)} \le 2\max_{\xi\in\La} |f(\xi)|.
$$
Indeed, let $\bz\in\Omega$ be one of the points, where the norm of $f\in B_N$ is achieved, say,
$|f(\bz)|= \|f\|_{\cC(\Omega)} =1$. Then for $\xi\in\La$ such that $|\bz-\xi|\le (2M)^{-1}$ we obtain
$$
|f(\xi)| \ge |f(\bz)|- |f(\bz)-f(\xi)| \ge 1- |\bz-\xi|\|\nabla f\|_{\cC(\Omega)} \ge 1/2.
$$
Typically, establishing discretization results for subspaces of smooth functions (first of all for
subspaces of polynomials) is reduced to the proof of the corresponding inequalities of the Bernstein-type.

{\bf D.28.} {We begin with results from  \cite{KaC1}  and its unpublished followup note \cite{HoC}} on discretization in the subspace $\varPi(d,N)$ of homogeneous polynomials on $d$ complex variables of degree $N$ restricted onto the complex sphere
$$
 \cS^d := \{\bz=(z_1,\dots,z_d) \in \bbC^d\,:\, |\bz| = (|z_1|^2+\dots+|z_d|^2)^{1/2}=1\}.
$$
It is proved in \cite{KaC1} that the absolute value of the tangential derivative of a function $|P|$, $P\in \varPi(d,N)$, is bounded by $CN^{1/2}\|P\|_{\cC(\cS^d)}$ and, therefore, for any two $|\bz^1|=|\bz^2|=1$
\be\label{AC1}
\left| |P(\bz^1)|-|P(\bz^2)|\right| \le CN^{1/2}d(\bz^1,\bz^2)\|P\|_{\cC(\cS^d)},
\ee
where $d(\bz^1,\bz^2)$ is the spherical distance between $\bz^1$ and $\bz^2$, and $C$ is an absolute constant. By the argument explained above (\ref{AC1}) implies (see \cite{KaC1}):
\be\label{AC2}
\varPi(d,N) \in \cM(C_dN^{d-1},\infty,1/2).
\ee
Note that $\dim \varPi(d,N) \asymp N^{d-1}$. Therefore, relation (\ref{AC2}) shows that alike subspaces of trigonometric polynomials subspaces $\varPi(d,N)$ allow good discretization with the number of sampling points of the same order as their dimensions. In the case of the uniform norm, the corresponding $\e$-nets of $\Omega$ with $\e\asymp N^{-1/2}$ can be taken as those  point sets. Relation (\ref{AC2}) was successfully used to construct uniformly bounded bases in the spaces of holomorphic functions in the $d$-dimensional ball $\cB^d$. Also, it was used to prove the existence
of nontrivial inner functions in $\cB^d$ (see \cite{BoC1}, \cite{BoC2}, and \cite{KaC1}). Note that
L.~H{\"o}rmander, \cite{HoC}, found the best constant in the inequality (\ref{AC1}):
$$
C_{min} = (1-1/N)^{(1-N)/2} \le e^{1/2},\quad N=1,2,\dots.
$$
Also, a class of surfaces for which an analog of (\ref{AC2}) holds was suggested in \cite{HoC}.

{\bf D.29.} An analog of (\ref{AC2}) for the $L_p$ norms, $1\le p\le \infty$, was obtained in \cite{WoC}.
Namely, it was proved in \cite{WoC} that for $1\le p\le \infty$
\be\label{AC3}
\varPi(d,N) \in \cM(C_d\dim(\varPi(d,N)),p,C_1(d),C_2(d)).
\ee
  The maximal $c_dN^{-1/2}$-distinguishable subset $\La\subset \Omega$ can be taken as a good discretization point set.

{\bf D.30.} In the case of functions on real variables, i.e. in the case of subspaces $\cP^d_n$ restricted to the unit sphere in $\R^d$, analogs of (\ref{AC3}) were considered later in \cite{MNWaC}. In this case a good discretization point set was constructed in a different way. The unit sphere is divided into disjoint subsets
of equal measure and of  diameters $\le c_dN^{-1}$. Then, a good discretization point set is formed by including one point from each subset of the subdivision. An interesting application of  the Marcinkiewicz--Zygmund inequalities from \cite{MNWaC} (the case $p=1$) was found in \cite{BRViC}. The authors of \cite{BRViC} solved
the problem of existence of a cubature formula with equal weights defined on the unit sphere of $\R^d$
with the minimal in the sense of order number of points ($\asymp N^{d-1}$), which is exact for algebraic polynomials of total degree $\le N$ (see  \cite{BRV15}   for further developments).

{\bf D.31.} We now briefly comment on some results on discretization on manifolds (see \cite{FM}, \cite{PeC}, \cite{BD} and references therein). Let $M$ be a compact Riemannian manifold. Consider a negative definite self adjoint elliptic differential operator of the second order defined on $L_2(M)$. Optimal in the sense of order results on sampling discretization of the $L_p$ norm, $1\le p\le \infty$, of elements of subspaces
generated by eigenfunctions corresponding to the eigenvalues of the operator not exceeding a given number $N$ are obtained in \cite{PeC}.  Similar results on sampling discretization were also established  on compact two-point homogeneous spaces in the paper \cite{BD}.

{\bf D.32.} The author of \cite{SkC} found sufficient conditions on a generating function $\varphi$, which guarantee that the Marcinkiewicz--Zygmund inequalities hold for subspaces spanned by the corresponding system of periodic wavelets for all $L_p$, $1\le p\le \infty$.

It is clear that many of these results in Section \ref{A} can be formulated as results on the
quasi-matrix property (see Definition \ref{ADqm} above) of the corresponding system. We note that the concept of quasi-matrix system was used
in \cite{KTeC} for establishing general lower bounds for the best $m$-term approximation in the $L_1$ norm.

\subsection{Sampling recovery}
\label{rec}

The problem of recovery (reconstruction) of an unknown function defined on a subset of  $\R^d$ from its samples at a finite number of points is a fundamental problem of pure and applied mathematics. The goal is to construct recovering operators (algorithms) which are good in the sense of accuracy, stability, and computational complexity. In this subsection we discuss
the issue of accuracy.  It turns out that, in the case of recovery in
the $L_2$ norm, the weighted least squares algorithms are reasonably good recovering methods. The corresponding analysis is based on recent deep results in discretization of the $L_2$ norms of
functions from finite-dimensional subspaces  (see {\bf D.18} above).

 Recall the setting
 of optimal recovery. For a fixed $m$ and a set of points  $\xi:=\{\xi^j\}_{j=1}^m\subset \Omega$, let $\Phi$ be a linear operator from $\bbC^m$ into $L_p(\Omega,\mu)$.
Denote for a class $\bF$
$$
\varrho_m(\bF,L_p) := \inf_{\text{linear}\, \Phi; \,\xi}\,\, \sup_{f\in \bF} \|f-\Phi(f(\xi^1),\dots,f(\xi^m))\|_p.
$$
This recovery procedure is linear;
the following modification of it is also of interest. We allow any mapping $\Phi : \bbC^m \to X_N \subset L_p(\Omega,\mu)$ where $X_N$ is a linear subspace of dimension $N\le m$ and define
$$
\varrho_m^*(\bF,L_p) := \inf_{\Phi; \xi; X_N, N\le m} \sup_{f\in \bF}\|f-\Phi(f(\xi^1),\dots,f(\xi^m))\|_p.
$$

In both of the above cases we build an approximant, which comes from a linear subspace of dimension at most $m$.
It is natural to compare quantities $\varrho_m(\bF,L_p)$ and $\varrho_m^*(\bF,L_p)$ with 
Kolmogorov widths. Let $\bF\subset L_p$ be a centrally symmetric compact. The quantities
$$
d_n (\bF, L_p) := \operatornamewithlimits{inf}_{\{u_i\}_{i=1}^n\subset L_p}
\sup_{f\in \bF}
\operatornamewithlimits{inf}_{c_i} \left \| f - \sum_{i=1}^{n}
c_i u_i \right\|_p, \quad n = 1, 2, \dots,
$$
are called the {\it Kolmogorov widths} of $\bF$ in $L_p$. In the definition of
 Kolmogorov widths we take for $f\in \bF$ an approximating element
from $U := \sp \{u_i \}_{i=1}^n$ the element of best
approximation. This means
that in general (i.e. if $p\neq 2$) this method of approximation is not linear.

We have the following obvious inequalities
\begin{equation*}
d_m (\bF, L_p)\le \varrho_m^*(\bF,L_p)\le \varrho_m(\bF,L_p).
\end{equation*}

 {The characteristics $\varrho_m$ and $\varrho_m^*$ are well studied for many particular classes of functions. For an exposition of known results we refer to the books 
\cite{TWW}, \cite{NoLN}, \cite{DTU}, \cite{VTbookMA}, \cite{NW1}--\cite{NW3} and references therein. The characteristics $\varrho_m^*$ and $\varrho_m$ are inspired by the concepts of 
the Kolmogorov width and the linear width. Probably, $\varrho_m^*$ were introduced in \cite{TWW} and $\varrho_m$ in \cite{VT51}.}

In {\bf R.1} and {\bf R.2} we consider the case $p=2$, i.e. recovery takes place in the Hilbert space $L_2$.

{\bf R.1.} The following general inequality was proved in \cite{VT183} with the aid of result {\bf D.18}:
   There exist two positive absolute constants $b$ and $B$ such that for any   compact subset $\Omega$  of $\R^d$, any probability measure $\mu$ on it, and any compact subset $\bF$ of $\cC(\Omega)$ we have
\be\label{R1}
\ro_{bn}(\bF,L_2(\Omega,\mu)) \le Bd_n(\bF,L_\infty).
\ee

Note, that it is pointed out in \cite{LT} (see Theorem 3.4 of v2) that in the case of real spaces the
constant $b>1$ can be taken arbitrarily close to $1$ and in the case of complex spaces -- the constant $b>2$ arbitrarily close to $2$ with the constant $B$ allowed to depend on $b$.
 We refer the reader to the very recent paper \cite{CoDo} {(also see references therein)} for a discussion of the
$L_2$ recovery when sampling points are drawn according to some random distribution.

{\bf R.2.}  {Probably, the first paper, which gave a general upper bound on the sampling numbers $\ro_{n}(\bF,L_2)$ in terms of the Kolmogorov numbers $d_k (\bF, L_2)$ is \cite{WW}.  Recently an outstanding progress has been done in the sampling recovery in the $L_2$ norm (see, for instance, \cite{CM}, \cite{KU}, \cite{KU2},  \cite{NSU}, \cite{CoDo}, \cite{TU1}, and  \cite{KUV}, \cite{KKP}). We give a very brief comments on those interesting results. }
For special sets $\bF$ (in the reproducing kernel Hilbert space setting) the following inequality is known (see \cite{NSU} and \cite{KU}):
\be\label{R2}
\ro_{n}(\bF,L_2) \le C\left(\frac{\log n}{n}\sum_{k\ge cn} d_k (\bF, L_2)^2\right)^{1/2}
\ee
with absolute constants $C,c>0$.  We refer the reader for further results in this direction to the   paper \cite{KU2}. {In particular, the close to optimal bounds on $\ro_{n}(\bF,L_2)$ in terms of $d_k (\bF, L_2)$ are obtained there for general sets $\bF$.}

The above results were proved with the aid of a classical type of algorithm -- weighted least squares. {A variant of this method -- the empirical risk minimization -- is a standard tool in learning theory and nonparametric statistics (see, for instance, \cite{GKKW} and \cite{VTbook}, Ch.4). 
The weighted least squares algorithm was successfully applied in the sampling recovery for the first time in \cite{CM}. }

Let $X_N$ be an $N$-dimensional subspace of the space of continuous functions $\cC(\Omega)$ and let $\bw:=(w_1,\dots,w_m)\in \R^m$ be a positive weight, i.e. $w_i>0$, $i=1,\dots,m$. Consider the following classical weighted least squares recovery operator (algorithm) (see, for instance, \cite{CM})
$$
 \ell 2\bw(\xi,X_N)(f):=\text{arg}\min_{u\in X_N} \|S(f-u,\xi)\|_{2,\bw},\quad \xi=\{\xi^j\}_{j=1}^m\subset \Omega,
$$
where
$$
\|S(g,\xi)\|_{2,\bw}:= \left(\sum_{\nu=1}^m w_\nu |g(\xi^\nu)|^2\right)^{1/2} .
$$

The above weighted least square operator is a linear operator. This linearity guarantees that its application
gives an upper bound for $\ro_{n}(\bF,L_2)$. Along with the operator
 $ \ell 2\bw(\xi,X_N)(f)$ we
can  also consider the operator
 $ \ell p\bw(\xi,X_N)(f)$, $1\le p<\infty$, and the operator
  $ \ell \infty(\xi,X_N)(f)$. These operators are not linear operators for $p\neq 2$. Therefore, their applications provide upper bounds for $\ro_{n}^*(\bF,L_p)$ but not for $\ro_{n}(\bF,L_p)$.

{\bf R.3.} Inequality (\ref{R1}) gives a  bound on sample recovery in $L_2$ with the right hand side
being a characteristic of the set $\bF$ measured in the uniform norm. We now discuss an analog of
(\ref{R1}) for the $L_p$ spaces with $p\in (2,\infty)$.  If instead of
 {\bf D.18} we use {\bf D.20.~An upper bound},  then the proof from \cite{VT183} allows us to establish an analog of (\ref{R1}). Then similarly to {\bf R.1} we obtain the following statement.
Let $p=2s$, $s\in \N$.   There exist two positive absolute constants $b$ and $B$ such that for any   compact subset $\Omega$  of $\R^d$, any probability measure $\mu$ on it, and any compact subset $\bF$ of $\cC(\Omega)$ we have
\begin{equation*}
\ro_{bn^{p/2}}^*(\bF,L_p(\Omega,\mu)) \le Bd_n(\bF,L_\infty).
\end{equation*}

{\bf R.4.} In the case $p=\infty$ we use the proof from \cite{VT183}
with {\bf D.18} replaced by {\bf D.21}. Then similarly to {\bf R.1}, {\bf R.3} we obtain the following statement.
There exists a positive absolute constant $B$ such that for any   compact subset $\Omega$  of $\R^d$ and any compact subset $\bF$ of $\cC(\Omega)$ we have
\begin{equation}\label{R1b}
\ro_{9^n}^*(\bF,L_\infty(\Omega)) \le Bd_n(\bF,L_\infty).
\end{equation}

{ It is interesting to compare the above inequality (\ref{R1b}) with the following result of E. Novak (see \cite{NoLN}, Proposition 1.2.5)
\be\label{Novak}
\ro_{n}^*(\bF,L_\infty(\Omega)) \le (n+1)d_n(\bF,L_\infty).
\ee
 Inequality (\ref{R1b}) complements inequality (\ref{Novak}) in the case, when the sequence $\{d_n(\bF,L_\infty)\}$ decays 
slowly, slower than $1/n$. 
}

{
 We refer the reader for some recent results on sampling recovery in the uniform norm to the   paper \cite{PoUl}.
}

Very recently the result {\bf R.1} was applied to estimate the error of optimal sampling recovery in $L_2$ of functions with small mixed smoothness (see \cite{TU1}). In the case of small smoothness we cannot use {\bf R.2} because the series in (\ref{R2}) diverges. This calls to mind an interesting phenomenon, which was discovered in \cite{TU1}. We illustrate that phenomenon on the example of
the $\bH^r_p$ classes. Let  $\btt =(t_1,\dots,t_d )$ and $\Delta_{\btt}^l f(\bx)$
be the mixed $l$-th difference with
step $t_j$ in the variable $x_j$, that is
$$
\Delta_{\btt}^l f(\bx) :=\Delta_{t_d,d}^l\dots\Delta_{t_1,1}^l
f(x_1,\dots ,x_d ) .
$$
Let $e$ be a subset of natural numbers in $[1,d ]$. We denote
$$
\Delta_{\btt}^l (e) =\prod_{j\in e}\Delta_{t_j,j}^l,\qquad
\Delta_{\btt}^l (\varnothing) = I .
$$
We define the class $\bH_{p,l}^r B$, $l > r$, as the set of
$f\in L_p$ such that for any $e$
\begin{equation*}
\bigl\|\Delta_{\btt}^l(e)f(\bx)\bigr\|_p\le B
\prod_{j\in e} |t_j |^r .
\end{equation*}
In the case $B=1$ we omit it. It is known (see, for instance, \cite{VTbookMA}, p.137) that the classes $\bH^r_{p,l}$ with different $l$ are equivalent. So, for convenience we fix one $l= [r]+1$ and omit $l$ from the notation.

There are constructive methods for the sampling recovery based on sparse grids (Smolyak point sets $SG(n)$) (see \cite{DTU}, Chapter 5 and \cite{VTbookMA}, Section 6.9). For instance, these methods give the following upper bounds for the sampling recovery for all $r>1/p$, $2\le p\le \infty$ (see  \cite{VTbookMA}, p.307, Theorem 6.9.2).
\be\label{C5}
\ro_m(\bH^{r}_p,L_2) \le C(r,p,d) m^{-r} (\log m)^{(d-1)(1+r)} .
\ee
It is also known that in many cases  sparse grid methods are optimal sampling recovery methods for classes of functions with mixed smoothness (\cite{VTbookMA}, Section 6.9.3). However, results of
\cite{TU1} show that it is not the case for the small smoothness situation. It is proved in \cite{TU1} with the aid of (\ref{R1}): Let $2<p\le\infty$ and $1/p<r<1/2$ then
\be\label{C3h}
\ro_m(\bH^r_p,L_2) \le C(r,p,d) m^{-r} (\log m)^{d-1+r}
\ee
and in the case $r=1/2$
\be\label{C4h}
\ro_m(\bH^{1/2}_p,L_2) \le C(r,p,d) m^{-1/2} (\log m)^{d-1/2}(\log \log m)^{3/2}.
\ee
It follows from known results (see, for instance, \cite{DU}, \cite{VTbookMA}, Section 6.4, and a discussion in \cite{TU1}) that the upper bound (\ref{C5}) cannot be improved if we use the $(n,n-1)$-nets, in particular the sparse grids (see \cite{VTbookMA}, Section 6.4 for their definitions), as  sample sets for recovery. Comparing above mentioned lower bounds for (\ref{C5}) and inequalities
(\ref{C3h}) and (\ref{C4h}), we conclude that, in the range of parameters  $1/p<r\le 1/2$, $2< p<~\infty$, ~$d>2$,
there exists a weighted least squares algorithm, which provides better (albeit, nonconstructive) upper bounds for sampling recovery than algorithms based on sparse grids or even based on a wider class of point sets -- the $(n,n-1)$-nets.

\subsection{Some comments on techniques}
\label{tech}

It turns out that the technique used in sampling discretization   is rather rich and connected with deep techniques from different areas of mathematics -- functional analysis, harmonic analysis, probability, and others.
In   brief, here are a few of these connections.

{\bf T.1.} {The} result in {\bf D.1} is based on a simple well known identity:
$$
\sum_{l=0}^{2n} e^{ik2\pi l/(2n+1)} = \begin{cases} 2n+1,& k=0, \\
0, & 0<|k| \le 2n. \end{cases}
$$

{\bf T.2.} Result {\bf D.3} is based on the Riesz-Thorin interpolation theorem (see, for instance, \cite{VTbookMA}, Sections 1.3.3 and 3.3.4).

{\bf T.3.} Sampling discretization for the hyperbolic cross polynomials $\Tr(Q_n)$ turns out to be a very difficult and interesting problem. It may serve as a testing field for different methods. These results were obtained with the aid of probabilistic techniques, in particular: a variant of the Bernstein concentration measure inequality
from \cite{BLM} (see Lemma \ref{L2.1} below), the chaining technique from \cite{KoTe} (see also \cite{VTbook}, Ch.4), and bounds of the entropy numbers.  The reader can find further results on the chaining technique in \cite{Ta1}.

We formulate a lemma, which is an analog of the well known concentration measure inequalities (see, for instance \cite{VTbook}, Ch.4). Lemma \ref{L2.1} is from \cite{BLM}.
\begin{Lemma}\label{L2.1} Let $\{g_j\}_{j=1}^m$ be independent random variables with $\bE g_j=0$, $j=1,\dots,m$, which satisfy
$$
\|g_j\|_1\le 2,\qquad \|g_j\|_\infty \le M,\qquad j=1,\dots,m.
$$
Then for any $\eta \in (0,1)$ we have the following bound on the probability
$$
\bP\left\{\left|\sum_{j=1}^m g_j\right|\ge m\eta\right\} < 2\exp\left(-\frac{m\eta^2}{8M}\right).
$$
\end{Lemma}

This lemma,\! in  combination with different chaining techniques and bounds on the entropy numbers of the  unit $L_q$-ball in the uniform norm, was used in {\bf D.4}, {\bf D.6}, {\bf D.10}, and {\bf D.13}.

{\bf T.4.} For all $q\in [1,\infty)$ distinct from $2$, we obtain the sampling discretization of the $L_q$ norm results close to the ones  optimal  in the sense of order. The necessary condition for
successful discretization in an $N$-dimensional subspace is $m\ge N$. The results in {\bf D.16} and {\bf D.17}  
provide good discretization for $m$ of order $N(\log N)^a$ with $a>0$ depending on the subspace.
It seems like the probability technique used does not allow us to get rid of the logarithmic factor in the bound for $m$. It turns out that in the case $q=2$ we can get rid of the logarithmic factor and obtain
optimal (in the sense of order) results. This technique is based on results by A.~Marcus, D.A.~Spielman, and
N.~Srivastava from \cite{MSS} (see Corollary 1.5 with $r=2$ there) obtained for solving the Kadison-Singer problem (see Theorem \ref{MSS} from Section \ref{Int}).   We now discuss this technique.

{The following  simple technical result (see \cite{DPTT}, Proposition $2.1$)  illustrates reduction of a continuous problem to the discrete one.
\begin{Theorem}\label{Pr2} Let $Y_N:=\sp\{u_1(x),\dots,u_N(x)\}$ with $\{u_i(x)\}_{i=1}^N$ being a real (or complex) orthonormal on $\Omega$ with respect to a probability measure $\mu$ basis for $Y_N$. Assume that $\|u_i\|_4:=\|u_i\|_{L_4(\Omega,\mu)} <\infty$ for all $i=1,\dots,N$.   Then for any $\de>0$ there exists
    a set $\Omega_M=\{x^j\}_{j=1}^M\subset \Omega$ such that for any $f\in Y_N$
    \begin{equation*}
   | \|f\|_{L_2(\Omega)}^2 - \|f\|_{L_2(\Omega_M)}^2| \le \de \|f\|_{L_2(\Omega)}^2,
    \end{equation*}
    where
    $$
    \|f\|_{L_2(\Omega_M)}^2 := \frac{1}{M}\sum_{j=1}^M |f(x^j)|^2.
    $$
\end{Theorem}
}
{ We now discuss some lemmas, that are used for solving the Marcinkiewicz-\-type discretization problem in $L_2$ for discrete sets.} 
{The following Lemma \ref{HO_CorE2} (see \cite{HO}, Corollary E.2) by N.J. Harvey and N. Olver was derived from Theorem \ref{MSS}.
\begin{Lemma}\label{HO_CorE2} Let $\bv_1, \dots, \bv_M\in\mathbb{R}^N$ be such that for all $\bw\in\mathbb{C}^N$
$$
\alpha\|\bw\|_2^2\leq \sum\limits_{j=1}^M |\langle \bw, \bv_j\rangle|^2\leq \beta\|\bw\|_2^2
$$ and $\|v_j\|^2=\delta:=N/M$ for all $j$. Suppose that $\alpha\in[1/2, 1]$ and $\beta\in[1,2]$. Then there exists $J\subset\{1, 2,\dots, M\}$ satisfying
$$
(\alpha -C\sqrt{\delta})\|\bw\|_2^2\leq 2 \sum\limits_{j\in J} |\langle \bw, \bv_j\rangle|^2\leq (\beta+C\sqrt{\delta})\|\bw\|_2^2,
$$
where $C$ is an absolute constant $C\geq 1$. 
\end{Lemma}
The following Lemma \ref{NOUL2} (see Lemma 2  in \cite{NOU}  and also see \cite{OU}, Lemma 10.22, p.105) was proved with the aid of Lemma \ref{HO_CorE2} by using the iteration technique which is discussed in Subsection \ref{Em}.}

\begin{Lemma}\label{NOUL2} Let a system of vectors $\bv_1,\dots,\bv_M$ from $\bbC^N$ satisfy \eqref{A5c} for all $\bw\in \bbC^N$
and
\begin{equation*}
\|\bv_j\|_2^2 = N/M,\qquad j=1,\dots,M.
\end{equation*}
Then there is a subset $J\subset \{1,2,\dots,M\}$ such that for all $\bw\in\bbC^N$
\begin{equation*}
c_0 \|\bw\|_2^2 \le \frac{M}{N} \sum_{j\in J} |\<\bw,\bv_j\>|^2 \le C_0\|\bw\|_2^2,
\end{equation*}
where $c_0$ and $C_0$ are some absolute positive constants.
\end{Lemma}

A simple Remark \ref{R4.1} is from \cite{VT158}.

\begin{Remark}\label{R4.1} For the cardinality of the subset $J$ from Lemma \ref{NOUL2} we have
$$
c_0 N \le |J| \le C_0 N.
$$
\end{Remark}

{We note that the iteration technique was already used in the paper \cite{HO} in order to obtain the following result (see Claim E.3 in \cite{HO}), which is similar to Lemma \ref{NOUL2}). Here we give a slight reformulation of it, {which follows from the proof of Claim  E.3 in \cite{HO}}.
\begin{Lemma}
    Let a system of vectors $\bv_1,\dots,\bv_M$ from $\mathbb{R}^N$ satisfy \eqref{A5c} for all $\bw\in \mathbb{R}^N$
    and
    \begin{equation*}
    \|\bv_j\|_2^2 = N/M,\qquad j=1,\dots,M.
    \end{equation*}
    For any $\varepsilon\leq1/9$ there is a subset $J\subset \{1,2,\dots,M\}$ 
    with $|J|\leq C_1N/\varepsilon^2$,
    where $C_1$ is a positive absolute constant, such that for all $\bw\in\mathbb{R}^N$
  $$
  (1-\varepsilon) \|\bw\|_2^2 \le \frac{M}{|J|} \sum_{j\in J} |\<\bw,\bv_j\>|^2 \leq (1+\varepsilon)\|\bw\|_2^2.
$$
\end{Lemma}
}
The following Lemma \ref{Lim}, which is a generalization of Lemma \ref{NOUL2}, and
Corollary \ref{weighted} were proved in \cite{LT}.

\begin{Lemma}\label{Lim} Let a system of vectors $\bv_1,\dots,\bv_M$ from $\bbC^N$ satisfy \eqref{A5c} for all $\bw\in \bbC^N$
and
\begin{equation*}
\|\bv_j\|_2^2 \le \theta N/M, \qquad  \theta \le M/N,\qquad j=1,\dots,M.
\end{equation*}
Then there is a subset $J\subset \{1,2,\dots,M\}$ such that for all $\bw\in\bbC^N$
\begin{equation*}
c_0\theta \|\bw\|_2^2 \le \frac{M}{N} \sum_{j\in J} |\<\bw,\bv_j\>|^2 \le C_0\theta\|\bw\|_2^2,\quad |J|\le C_1\theta N,
\end{equation*}
where $c_0$, $C_0$, and $C_1$ are some absolute positive constants.
\end{Lemma}

 We note that the proof of Lemma \ref{Lim} in \cite{LT} gives a slightly stronger result than Lemma \ref{Lim} -- the tight frame condition (\ref{A5c}) can be replaced by a frame condition
$$
A\|\bw\|_2^2 \le\sum_{j=1}^M |\<\bw,\bv_j\>|^2 \le B \|\bw\|_2^2.
$$
This stronger version of Lemma \ref{Lim} was formulated in the followup paper
\cite{NSU}, where it was used for sampling recovery. Also, \cite{NSU} gives a proof with
specified constants $c_0$, $C_0$, and $C_1$.

\begin{Corollary}\label{weighted} Let a system of vectors $\bv_1,\dots,\bv_M$ from $\bbC^N$ satisfy \eqref{A5c} for all $\bw\in \bbC^N$.
Then there exists a set of weights $\lambda_j\ge 0$, $j=1,\dots, M$, such that $|\{j: \lambda_j\neq 0\}| \le  2C_1N $ and for all $\bw\in\bbC^N$ we have
\begin{equation*}
c_0\|\bw\|_2^2 \le \sum_{j=1}^M \lambda_j|\<\bw, \bv_j\>|^2 \le C_0\|\bw\|_2^2.
\end{equation*}
 where $c_0$, $C_0$, and $C_1$ are absolute positive constants from Lemma \ref{Lim}.
\end{Corollary}

Lemma \ref{NOUL2} was used in {\bf D.5} and {\bf D.9}.  Lemma \ref{Lim} was used in {\bf D.15}
and Corollary \ref{weighted} was used in {\bf D.18}. We point out that results based on Lemma \ref{Lim}
and Corollary \ref{weighted} found their applications in the optimal recovery problems (see \cite{NSU} and \cite{VT183}).

{\bf T.5.} In \cite{Kos}  Talagrand's generic chaining and its development
due to van Handel (see \cite{VH}) have been applied.
A key step in the proofs by chaining is to bound certain entropy numbers.
In \cite{Kos} this was accomplished by the following theorem (see \cite{Ta1}, p.552, Lemma 16.5.4). For a Banach space $X$ we define the modulus of smoothness
$$
\rho(u):= \rho(X,u) := \sup_{\|x\|=\|y\|=1}\left(\frac{1}{2}(\|x+uy\|+\|x-uy\|)-1\right).
$$
Let $\D_n=\{g_j\}_{j=1}^n$ be a system of  elements of cardinality $|\D_n|=n$ in a Banach space $X$. We equip the linear space
 $W_n:=[\D_n]:= \sp\{\D_n\}$   with the norm
 \begin{equation*}
 \|f\|_A := \|f\|_{A_1(\D_n)} := \inf\left\{ \sum_{j=1}^n |c_j|\, : \, f = \sum_{j=1}^n c_jg_j\right\}.
 \end{equation*}
 Denote by $W_{n,A}$ the $W_n$ equipped with the norm $\|\cdot\|_A$.
 We are interested in the dual norm to the norm $\|\cdot\|_A$, which we denote $\|\cdot\|_U$:
 $$
 \|F\|_U:= \|F\|_{U(\D_n)}:=   \sup_{f\in W_n; \|f\|_A\le 1} |F(f)|.
 $$
 Denote $W_{n,U}^*$ the $W_n^*$ equipped with the norm $\|\cdot\|_U$. Note that $\|\cdot\|_U$ is a semi-norm on the dual to $X$, space $X^*$.

 \begin{Theorem}\label{IT2c} Let $X$ be $q$-smooth: $\rho(X,u) \le \gamma u^q$, $1<q\le 2$ and let $\D_n$ be a  normalized system in $X$ of cardinality $|\D_n|=n$.  Then for the unit ball $B(X^*)$ of $X^*$ we have
  \be\label{I5c}
 \e_k(B(X^*),\|\cdot\|_{U(\D_n)}) \le C(X) \left(\frac{\log n}{k}\right)^{1-1/q},\quad k=1,\dots,n.
 \ee
 \end{Theorem}
We note that  Talagrand's inequality (\ref{I5c}) was slightly improved in \cite{VT181}, where a different technique, based on a combination of the deep duality result for the entropy numbers proved in \cite{BPST} and results from greedy approximation, was used: The factor  $\log n$
was replaced by  $\log (2n/k)$.  The following result from \cite{Kos},
 was proved with the aid of Theorem \ref{IT2c}.

\begin{Lemma}\label{BL1} Let $q\in (2,\infty)$. Assume that for any $f\in X_N$ we have
\be\label{B6'c}
\|f\|_\infty \le M\|f\|_q
\ee
with some constant $M$. Then for $k\in [1,N]$ we have for any finite subset $Y_s\subset \Omega$ of $s$ points
\begin{equation*}
\e_k(X_N^q, L_\infty(Y_s)) \le C(q)(\log s)^{1/q} (MN^{-1/q})(N/k)^{1/q} .
\end{equation*}
\end{Lemma}

This technique was used in {\bf D.14} and {\bf D.17}.

We now comment on the relation between the Nikol'skii inequality condition (\ref{B6'c}) and
the entropy condition (\ref{A12}). On one hand it is easy to see (see \cite{DPSTT2}) that the entropy condition (\ref{A12}) for $k=1$ implies the following Nikol'skii inequality
$$
\|f\|_\infty \le 4BN^{1/q}\|f\|_q \quad \text{for all}\quad f\in X_N.
$$
On the other hand, we note that the Nikol'skii-type
inequality condition \eqref{B6'c} with $M=B N^{1/q}$ implies
 the entropy condition \eqref{A12} with $k=1$. Thus, the Nikol'skii-type
inequality condition is equivalent to the  the entropy condition (\ref{A12}) for $k=1$.
Moreover, Lemma \ref{BL1} shows that the Nikol'skii inequality (\ref{B6'c}) implies
the entropy condition in the $L_\infty(Y_s)$ norm for the whole range $k\in [1, N]$. Further, the following Lemma \ref{BL2} from
\cite{VT180} shows that the Nikol'skii inequality (\ref{B6'c}), combined with an extra mild condition, implies  the entropy condition in the $L_\infty$ norm.

\begin{Lemma}\label{BL2} Let $q\in (2,\infty)$. Assume that for any $f\in X_N$ we have
\begin{equation*}
\|f\|_\infty \le BN^{1/q}\|f\|_q
\end{equation*}
with some constant $B$. Also, assume that $X_N \in \cM(s,\infty,C_1)$ with $s\le aN^c$.
Then for $k\in [1,N]$ we have
\begin{equation*}
\e_k(X_N^q, L_\infty) \le C(q,a,c,C_1)(\log N)^{1/q}B(N/k)^{1/q} .
\end{equation*}
\end{Lemma}

\subsection{Some further comments and open problems}
\label{OP}

In Subsections \ref{Atr}--\ref{cap} we provided some sufficient conditions for the Mar\-cin\-kiewicz-type theorem
to hold, but optimal conditions would be more useful. We now formulate the corresponding optimization problem the formulation of which is the same in both the Marcinkiewicz
$\cM(m,L_q(\Omega,\mu),C_1,C_2)$ theorem and in its weighted analog $\cM^w(m,L_q(\Omega,\mu),C_1,C_2)$. For an $N$-dimensional subspace $X_N\subset L_q(\Omega,\mu)$, $1\le q<\infty$, define
$$
m^w(X_N,L_q(\Omega,\mu),C_1,C_2) := \inf\{m\,:\, X_N\in \cM^w(m,q,C_1,C_2)\}
$$
and
\begin{equation*}
sd^w(N,q,C_1,C_2) := \sup_{X_N,\Omega,\mu} m^w(X_N,L_q(\Omega,\mu),C_1,C_2).
\end{equation*}
Here $sd$ stands for the {\it sampling discretization}. In the above definition of $sd^w(N,q,C_1,C_2)$ and in the similar definition of $sd(N,q,C_1,C_2)$ we take supremum over all $N$-dimensional subspaces
of $L_q(\Omega,\mu)$. In the case when $\Omega$ and $\mu$ are fixed and  $N$-dimensional subspaces of $L_q(\Omega,\mu)$ come
from a special collection $\cF$ we write
\begin{equation*}
sd^w(\cF,N,L_q(\Omega,\mu),C_1,C_2) := \sup_{X_N\in\cF} m^w(X_N,L_q(\Omega,\mu),C_1,C_2).
\end{equation*}

{We are interested in behavior of characteriscics $sd^w(N,q,C_1,C_2)$ and their analogs when $N$ goes to infinity and the constants $C_1$ and $C_2$ are fixed (do not depend on $N$). Certainly, this behavior may depend on a specific pair of constants $C_1$ and $C_2$. In order to focus attention on dependence on $N$ we talk about a generic characteristic 
$sd^w(N,q)$, which we do not define formally.
}

\begin{Definition}\label{AD1} We say that the sequence of  sampling weighted discretization numbers $sd^w(N,q)$ grows as $\ff(N)$ with $N\to \infty$ if there exist three positive numbers $C_i^*(q)$, $i=1,2,3$, which may depend on $q$ such that
$$
sd^w(N,q,C_1^*(q),C_2^*(q)) \le C_3^*(q)\ff(N),\quad N\in \N,
$$
and for any $0<C_1\le C_2<\infty$ there is a constant $C(C_1,C_2)>0$ such that we have
$$
sd^w(N,q,C_1,C_2) \ge C(C_1,C_2)\ff(N),\quad N\in \N .
$$
\end{Definition}
We use similar definitions for the sequences
$$
\{sd(N,q)\},\quad \{sd(\cF,N,L_q(\Omega,\mu))\},\quad \text{and}\quad  \{sd^w(\cF,N,L_q(\Omega,\mu))\}.
$$
It is clear that for all characteristics mentioned above we have $\ff(N)\ge N$.

{\bf Main Problem.} The main problem of sampling discretization is to find the growth rate of the sequences $\{sd(N,q)\}$, $\{sd^w(N,q)\}$ for all $1~\le~q<~\infty$, the growth rate of the sequence $\{sd(N,\infty)\}$, and the growth rate of the sequences $\{sd(\cF,N,L_q(\Omega,\mu))\}$ and
 $\{sd^w(\cF,N,L_q(\Omega,\mu))\}$ for $1~\le~ q~<~\infty$, $\{sd(\cF,N,L_\infty(\Omega,\mu))\}$ for special collections $\cF$.

We point out that the Main Problem is solved in only a few cases. Results of {\bf D.18} claim that
the sequence $\{sd^w(N,2)\}$ grows as $N$. Results of {\bf D.15} claim that for any $\Omega$ and $\mu$ the sequence $\{sd(\cN(2,\infty,tN^{1/2}),N,L_2(\Omega,\mu))\}$ grows as $N$. Here
$$
\cN(p,\infty,M) := \{X_N\,:\, \|f\|_\infty \le M\|f\|_p,\quad \forall f\in X_N\subset L_p(\Omega,\mu)\}.
$$

Another important collection of subspaces is the collection $\cE(N,q,B)$ of $N$-dimensional subspaces
satisfying the entropy condition (\ref{A12}).

{\bf Open problem 1.} Find the growth rate of the sequence \newline $\{sd(N,q)\}$, $q\in [1,\infty]$.

{\bf Open problem 2.} Find the growth rate of the sequence \newline $\{sd^w(N,q)\}$, $q\in [1,\infty)\setminus\{2\}$.

{\bf Open problem 3.} Find the growth rate of the sequence  \newline$\{sd(\cN(2,\infty,tN^{1/2}),N,L_q(\Omega,\mu))\}$, $q\in [1,\infty)\setminus\{2\}$.

 {\bf Open problem 4.} Find the growth rate of the sequence  \newline$\{sd(\cN(q,\infty,tN^{1/q}),N,L_q(\Omega,\mu))\}$, $q\in [1,\infty)\setminus\{2\}$.

{\bf Open problem 5.} Find  the growth rate of the sequence  \newline$\{sd(\cE(N,q,B),N,L_q(\Omega,\mu))\}$, $q\in [1,\infty)\setminus\{2\}$.

The problem of the growth rate of the $sd$-sequences is open even in the case of trigonometric polynomials.
Denote for $\Omega=[0,2\pi]^d$ and $\mu$ the normalized Lebesgue measure on $\Omega$
$$
\Tr r(N,d) := \{ \Tr(Q)\,:\, Q\subset \Z^d,|Q|=N\}.
$$

{\bf Open problem 6.} Find the growth rate of the sequence \newline $\{sd(\Tr r(N,d),N,L_q(\Omega,\mu))\}$, $q\in [1,\infty)\setminus\{2\}$.

Results {\bf D.7}, {\bf D.12}, and {\bf D.21} show that the case $q=\infty$ is very different from the case $1\le q<\infty$.

{\bf Open problem 7.} Find the growth rate of the sequence \newline $\{sd(\Tr r(N,d),N,L_\infty(\Omega,\mu))\}$.

Note, that results {\bf D.12}, and {\bf D.21} give
$$
e^{c_1N} \le sd(\Tr r(N,d),N,L_\infty(\Omega,\mu),1/2) \le e^{c_2N} ,
$$
where $c_2\ge c_1>0$ are absolute constants.

\section{Spectral properties and operator norms of submatrices}
\label{M}
Theorem \ref{MSS} mentioned in Introduction is actually a theorem on partitioning a matrix into two submatrices with extremely small $(2,2)$-norms (see \eqref{pq_norm_def} for the definition). It can be reformulated as follows:
\begin{Theorem}\label{reform_matrix}
Let $A$ be an $M\times N$ matrix with columns forming an orthonormal system in $\ell_2^M$  and rows $\bv_j$, $j=1,\dots, M$, such that $\|\bv_j\|_2\leq \e$, for some $\e>0$. Then there exists a partition of the set $\langle M\rangle:=\{1,2,\dots,M\}$ into two parts
\begin{gather}\label{partition}
    \langle M\rangle=\Omega_1\sqcup\Omega_2
\end{gather}
such that
\be\label{MSS_norm_est}
\|A(\Omega_k)\|^2\leq \frac{(1+\sqrt{2}\e)^2}{2}, \qquad k=1,2,
\ee
where for $\omega\subset \langle M\rangle$ we define $A(\omega)$ as a submatrix of $A$ obtained from $A$ by extracting the rows numbered in $\omega$.
\end{Theorem}
For a matrix satisfying conditions of Theorem \ref{reform_matrix} inequality \eqref{MSS_norm_est} is equivalent to the following {\it pointwise estimate}: For  $k=1,2$
\be\label{MSS_pnt_est}
\|A(\Omega_k)\bx\|_2^2\leq \frac{(1+\sqrt{2}\e)^2}{2} \|A\bx\|_2^2, \qquad \forall \bx\in\bR^N.
\ee
As we saw in Subsection \ref{tech}, the best in the sense of order Marcinkiewicz-type discretization theorem for finite-dimensional subspaces in $L_2$ can be derived from Theorem \ref{reform_matrix}. There is a chance that similar results for other norms of matrices can improve known Marcinkiewicz-type discretization theorems.

In this section we present known results on partitioning a matrix into two submatrices with extremely small $(p,q)$--norms (as in \eqref{MSS_norm_est}) and with good pointwise estimate of the $\ell_q$ norm of the image of any vector (as in \eqref{MSS_pnt_est}) for different $1\leq p, q\leq \infty$.
 More precisely, {\it extremely small $(p,q)$--norm} here is understood as the property of the norm raised to the power $q$ to be extremely close to the half. Let us comment on these two types of problems. 

Let $A$ be an $M\times N$ matrix. We assume that $M\geq N$ and consider $A$ as an operator acting between $\bR^N$ and $\bR^M$. We start with the problem of partitioning of a matrix $A$ into two submatrices with good pointwise estimates for $\ell_q$ norm, $1\leq q\leq \infty$, of the image of any vector.
We say that partition \eqref{partition}
 has a good pointwise estimate if the following inequality holds for $k=1,2$:
\begin{gather}\label{pointwise_def}
\left\|A(\Omega_k)\bx\right\|_q^q\leq\left(\frac{1}{2}+\varphi(A)\right)\left\|A\bx\right\|_q^q, \qquad \forall \bx\in\bR^N
\end{gather}
with some parameter $0<\varphi(A)<1/2$ (depending only on $A$). Quality of the partition is measured by how small $\varphi(A)$ is.
For a partition of the form \eqref{partition} we obviously have
\begin{gather*}
  \left\|A(\Omega_1)\bx\right\|_q^q+\left\|A(\Omega_2)\bx\right\|_q^q= \left\|A\bx\right\|_q^q, \qquad \forall \bx\in\bR^N.
\end{gather*}
So, we cannot have the factor in the right-hand side of \eqref{pointwise_def} less than $1/2$.
For this reason, if inequality (\ref{pointwise_def}) or its analog for the norm is satisfied, then we say that the submatrix $A(\Omega_k)$ has {\it extremely small}  norm.

Thus if \eqref{pointwise_def} holds, then for $k=1,2$ we get
\begin{gather*}
  \left(\frac{1}{2}-\varphi(A)\right)\left\|A\bx\right\|_q^q \leq \left\|A(\Omega_k)\bx\right\|_q^q\leq\left(\frac{1}{2}+\varphi(A)\right)\left\|A\bx\right\|_q^q, \qquad \forall \bx\in\bR^N.
\end{gather*}
Actually, under some conditions on $A$ we can guarantee that $\varphi(A)$ is small enough (see, e.g. Theorem \ref{MSS} and {\bf M.\ref{K17}}). In this case we can find a submatrix $B$ of $A$ with {the number of rows not exceeding $M/2$} such that the quantities $2^{1/q}\left\|B\bx\right\|_q$ and $\left\|A\bx\right\|_q$ are almost the same for all $\bx\in\bR^N$.  In some cases we can continue and find a partition with good pointwise estimate for the matrix $B$ and so on {(see Subsection \ref{Em} for details)}. This process was applied in a number of papers (see, e.g., \cite{Lun}, \cite{NOU}, \cite{LT}, \cite{NSU}). The idea of such an iteration process was proposed by Lunin in the remarkable paper \cite{Lun}, where he obtained the first sharp estimates (see {\bf M.\ref{Ka80b}} for details).

Let $1\leq p\leq \infty$. By taking the supremum over all $\bx\in\bR^N$ with $\|\bx\|_p\leq 1$ in \eqref{pointwise_def}, we obtain the necessary condition for \eqref{pointwise_def} to hold, namely, the following inequality:
\begin{gather*}
\left\|A(\Omega_k)\right\|_{(p,q)}^q\leq\left(\frac{1}{2}+\varphi(A)\right)\left\|A\right\|_{(p,q)}^q, \qquad k=1,2.
\end{gather*}
This type of estimate is a different  type of result compared to the pointwise estimate.

 Another setting considered below concerns square matrices. Here problems are formulated not only on extracting rows of a matrix but simultaneous extraction of both rows and columns is allowed  (see, e.g. results  {\bf M.\ref{KTomega}} and {\bf M.\ref{And}} below). For $\omega\subset\langle N\rangle$ let $R_{\omega}^{N}$ be an $N\times N$ matrix with entries $r_i^j$ such that $r_i^i=1$ for $i\in \omega$ and other entries equal to $0$.
Note that for an $M\times N$ matrix $A$ and $\omega\subset\langle M\rangle $ the matrices $R_{\omega}^{M}A$ and $A(\omega)$ have the same non-zero rows.
Let $A$ be an $N\times N$ matrix. The goal is to find $\omega\subset\langle N\rangle$ such that $\|R_{\omega}^N A R_{\omega}^N\|$ is small enough. Due to a decoupling argument (see \cite[Proposition 1.9]{BT87}), this can be achieved if we have good estimates for the average of the quantity $\|R_{\omega}^N A R_{\omega'}^N\|$ over all $\omega,\omega'\subset\langle N\rangle$.

\subsection{Euclidean case}
{We begin our discussion with the  case that is most important for applications: The matrix $A$ is considered to be an operator acting between $\ell_2^N$ and $\ell_2^M$.}

\om \label{Ka80b} The first result on operator norms of submatrices is due to Kashin \cite{Ka80b}: For any $\e>0$ there exists $C(\e)>0$ such that for any $M\times N$ matrix $A$ with $\|A\|=1$ and
$M\geq C(\e)N$ there exists an $N\times N$ submatrix $A(\omega)$ with $\|A(\omega)\|\leq\e$. Later, in \cite{Lun} Lunin established the following sharp estimate of $C(\e)$: $C(\e)\leq B\e^{-2}$. Lunin's result is actually an estimate of the form \eqref{low_est} for arbitrary subspaces in $L_2(\Omega)$.

{In \cite{BT87} Bourgain and Tzafriri obtained results {\bf M.}\ref{restr} and {\bf M.}\ref{restr_}, which are widely known as the restricted invertibility theorems.}
\om \label{restr} For some absolute positive  constants $C_1$ and $C_2$ and an $N\times N$ matrix $A$ with $\|A(\bfe_i)\|_2=1$, $i=1,\dots, N$, ($\{\bfe_i\}_{i=1}^N$ is the standard basis of $\bR^N$) there exists a set $\omega\subset \langle N\rangle$ of cardinality $|\omega|\geq C_1N/\|A\|^2$, such that
$$
\|AR_{\omega}^N\bx\|_2\geq C_2\|R_{\omega}^N\bx\|_2, \qquad \forall \bx\in\mathbb{R}^N.
$$

\om \label{restr_} For every $C>0$ and $0<\e<1$, there is a constant $d=d(C, \e)>0$ such that, whenever $N\geq 1/d$ and $A$ is an $N\times N$ matrix with $\|A\|\leq C$ and $\langle A\bfe_i, \bfe_i\rangle=1$, $i=1,\dots, N$, there exists  $\omega\subset \langle N\rangle$ of cardinality $|\omega|\geq dN$ so that $R_{\omega}^NAR_{\omega}^N$ is invertible and
$$
\|(R_{\omega}^NAR_{\omega}^N)^{-1}\|<(1-\e)^{-1}.
$$

\om \label{KTomega} Here is a result from \cite{KaTz} on square matrices: Let $A$ be an $N\times N$ matrix. There exists a constant $C$ such that if $1/N<\e<1$, $1\leq r\leq N$, then
\begin{gather}\label{KTomega12}
\min_{|\omega_1|,|\omega_2|\geq \e N/4}
\|R_{\omega_1}^NAR_{\omega_2}^N\|\leq C \left(a_r(A)\e^{1/2}+\e+\left(\frac{\e r}{N}\right)^{1/2}\right)\|A\|,
\end{gather}
where $a_r(\omega)=\min \{\|A-B\|,  B \text{ is an } N\times N \text{matrix}, \rank B\leq r\}$.
The proof relies on averaging the quantity $\|R_{\omega_1}^NAR_{\omega_2}^N\|_{(\infty, 1)}$ and using  Grothendieck factorization.  By a decoupling argument the authors obtained a similar result with $\omega_1=\omega_2$ in \eqref{KTomega12} (see Corollary 6 in \cite{KaTz}) for $A$ of $\rank A=r$ with zeros on the diagonal.

\om \label{Rud_matrices} Here is the result from \cite{Rud} due to Rudelson (see also {\bf D.15, MM.1}): Let $t\geq 1$ and $A=[a_i^j]$, $j=1,\dots, M$, $i=1,\dots, N$, $M>N$, be a matrix with orthonormal columns. Suppose that for all $j$
$$
\sqrt\frac{M}{N}\cdot \left(\sum\limits_{i=1}^N (a_i^j)^2\right)^{1/2}\leq t.
$$
Then for every $\e>0$ there exists a set $\omega\subset \langle M\rangle$ so that
 
\be\label{log_est}
|\omega|\leq C \frac{t^2}{\e^2} N\log \frac{Nt^2}{\e^2}
\ee
and for all $\bx\in\bR^N$
$$
(1-\e)\|\bx\|\leq \sqrt{\frac{M}{|\omega|}} \|R_{\omega}^MA\bx\|\leq (1+\e)\|\bx\|.
$$
An example from \cite{KaTz} shows that estimate \eqref{log_est} is the best one can obtain by the random selection method. The obstacle here is that  in contrast to the upper estimate (see \cite{BT87}) the lower estimate does not hold in average.

The Kadison-Singer problem has a number of equivalent formulations. It was solved in \cite{MSS} by the new method of interlacing polynomials through proving Weaver's conjecture (see the stronger formulation of it in Theorems~\ref{MSS}, \ref{reform_matrix}, {\bf M.\ref{Weav}}).

\om \label{And} The following quantitative version of Anderson's paving conjecture \cite{And}, which is equivalent to the Kadison-Singer problem, was also proved in \cite{MSS}: For any $\e>0$, every zero-diagonal complex self-adjoint $N\times N$ matrix $T$ can be $(r, \e)$-paved with $r=(6/\e)^4$. This means that there are coordinate projections
 $P_1, \dots, P_r$ such that
 $\sum\limits_{i=1}^r P_i=I_{N}$ (identical matrix) and $\|P_iTP_i\|_{(2,2)}\leq \e\|T\|_{(2,2)}$ for all $i=1,\dots, r$.

\om \label{Weav}  We give here another formulation of Theorem \ref{reform_matrix} due to N.~Srivastava (see \cite{Sr}): For any vectors $\bv_1, \dots, \bv_M\in\mathbb{R}^N$ satisfying:
\be\label{cond_eucl}
|\langle\bv_{i_0}, \bx\rangle|^2\leq\varepsilon^2\sum_{i=1}^M|\langle\bv_i, \bx\rangle|^2,\quad \forall \bx\in\mathbb R^N, \quad \forall i_0\in\left<M\right>,
\ee
there is a partition of the form \eqref{partition} such that for $k=1,2$
\be\label{MSS_obtained}
\sum\limits_{i\in \Omega_k} |\langle\bv_i, \bx\rangle|^q\leq \left(\frac{1}{2}+O(\varepsilon)\right) \sum_{i=1}^M|\langle\bv_i, \bx\rangle|^q, \quad q=2, \quad \forall \bx\in\mathbb R^N.
\ee
In other words,  every quadratic form of this kind in which no term has too much influence can be divided into two quadratic forms which approximate it.

\om \label{Osw_} We present the result from \cite{Osw} which is applied to accelerate classical numerical methods for solving linear systems:  There exists an absolute constant $C$ such that for any $N\times N$ hermitian matrix $A$ there exists a permutation $\sigma$ of the set $\langle N\rangle$ for which
$$
\|L_{\sigma}\|\leq C\|A\|,
$$
where $L_{\sigma}$ is a strictly lower triangular part of $P_{\sigma}AP_{\sigma}^T$ ($P_{\sigma}$ here denotes the
associated $N\times N$ row permutation matrix). Anderson's paving conjecture was used to prove this.

\subsection{$p=2$, $q=1$}
Results for the case $p=2$, $q=1$ are interesting by themselves and they allow us to obtain important results for  {the} Euclidean case by using the Grothendieck factorization theorem.
\om \label{sh_2_1} Lunin established the following sharp estimate in \cite{Lun}: There is a constant $C$ such that for any $M\times N$ matrix there exists $\omega\subset\langle M\rangle$ with $|\omega|=N$ such that the following inequality holds:
$$
\|A(\omega)\|_{(2,1)}\leq C\frac{N}{M}\|A\|_{(2,1)}.
$$
Due to a factorization argument he deduced from this result the sharp estimate in {\bf M.\ref{Ka80b}}.
\om \label{KT1} The following result from \cite{KaTz} on averaging slightly modifies the estimate from \cite{BT91}: Let $\{\delta_i\}_{i=1}^M$ be a sequence of $\{0,1\}$--valued independent random variables of mean $0<\delta<1$ over some probability space $(\Omega, \mu)$ and for $\omega\in\Omega$ define {$\Omega_{\omega}=\{1\leq i\leq M: \delta_i(\omega)=1\}$.} Then there exists a constant $D$ such that for $M\geq D$ and $1\leq N\leq M$ we have
$$
\int\limits_{\Omega}\|A(\Omega_{\omega})\|_{(2,1)} d\mu\leq D(\delta
M^{1/2}+(\delta N)^{1/2})\|A\|_{(2,2)}.
$$

\om{\label{K17}} Let $I_N$ be the matrix of the identity operator in $\bR ^N$ and let $\varphi(\e):= 5\e^{1/2}\log^{1/4}(3/\e)$, $0<\e<1$. The following result from \cite{Ka17} is obtained by a modified approach proposed by Lunin in \cite{Lun} (see also \cite{Ka15}):
For some constant $\e_0>0$ and each $M\times N$ matrix $A$ such that $A^*A=I_N$ and the rows $\bv_j$, $j\in\langle M\rangle$, of $A$ satisfy $\|\bv_j\|_2\leq \e$ with $\e\leq \e_0$, there exists a partition of the form \eqref{partition} such that both of the norms $\|A(\Omega_1)\|_{(2,1)}$ and $\|A(\Omega_2)\|_{(2,1)}$ do not exceed $(1/2+\varphi(\e))\|A\|_{(2,1)}$. This result is similar to Theorem \ref{reform_matrix} for $(2,2)$-norms but it does not allow us to make a good iteration process.

\om \label{alw_2_1} We mention here a slight generalization of the result from \cite{KL19} (see also \cite{Lim16}) on the possibility of reducing $(2,1)$-norms of the matrix via partitioning into two submatrices: Let $1\leq q<2$,
$\|A\|_{(2,q)}=1$, $\|\bv_i\|_2\leq\e$, $i\in \langle M\rangle$. Then there exists a partition of the form \eqref{partition} with $||\Omega_1|-|\Omega_2||\leq 1$ such that
$$
\|A(\Omega_k)\|_{(2,q)}\leq \frac{1}{2^{\frac{1}{q}-\frac{1}{2}}}+2\e, \qquad k=1,2.
$$

\subsection{General case}
{In the case $(p,q)\notin \{(2,1), (2,2)\}$ no final (in the sense of order) results are known. The result {\bf M.\ref{sh_2_1}} is sharp, but there still remain important open questions in that problem, in particular, the question of obtaining two-sided inequalities is open.}
\om \label{restr_p} We start with formulating the restricted invertibility theorem for the $(p,p)$-norm from \cite{BT87}: For every $0<\e<1$, $1\leq p\leq \infty$ and $B>0$, there exists a constant $c=c(\e, p, B)$ such that, whenever $N\geq 1/c$ and $A$ is an $N\times N$ matrix with $\|A\|_{(p,p)}\leq B$ and ones on the diagonal, then one can find a subset $\sigma$ of $\{1,2,\dots, N\}$ of cardinality $|\sigma|\geq cN$ so that $R_{\sigma}^NAR_{\sigma}^N$ is invertible and
$$
\|(R_{\sigma}^NAR_{\sigma}^N)^{-1}\|_{(p,p)}\leq (1-\e)^{-1}.
$$

For the rest of this section we consider the problem of partitioning a matrix into two submatrices with small norms and with good pointwise estimate. We start with the case of $(1,q)$--norm, $1\leq q<\infty$. From convexity of the function $\|A(\bx)\|_q$ it is easy to see that $\|A\|_{(1,q)}$ is the maximal $\ell_q$-norm of the columns of $A$. Certainly, this quantity is much easier to deal with than in the case $1<p<\infty$, so here we need to impose a significantly weaker condition on the matrix $A$ to have a partition of a matrix with almost equal norms raised to the power $q$. Let
\begin{gather}\label{1_q_cond}
|a^j_i|\leq \e\|\bw_i\|_q,\qquad i\in\langle N\rangle, \qquad j\in\langle M\rangle,
\end{gather}
where $\bw_i$, $i=1,\dots, N$, are the columns of $A$.
\om We formulate positive results (for items a) and c) see \cite{Lim20}, for  item b) see \cite{Gl} and \cite{Lim20}): Assume that, for an $M\times N$ matrix $A$, inequality \eqref{1_q_cond} holds for some $0<\e<1$. Then there exists a partition of the form \eqref{partition} such that, for $k=1,2$ the following inequalities hold:

$\text{a) }\left\| A(\Omega_k)\right\|_{(1,q)}\leq
\left(\frac{1}{2}+\frac{3}{2}\varepsilon^{q/3}\ln^{1/3}{(4N)}\right)^{1/q}
\left\| A\right\|_{(1,q)},$

$\text{b) }\left\| A(\Omega_k)\right\|_{(1,q)}\leq
\left(\frac{1}{2}+\frac{1}{2}\varepsilon^{q}\sqrt{M}(1+\log(\frac{N}{M}+1))^{1/2}\right)^{1/q}
\left\| A\right\|_{(1,q)},$

$\text{c) } \left\| A(\Omega_k)\right\|_{(1,q)}\leq\left(\frac{1+N\varepsilon^q}{2}\right)^{1/q}\left\| A\right\|_{(1,q)}.$

Note that in the proof of c), the generalized ham sandwich theorem (see, e.g. \cite{Gray}) was applied.

\om \label{p_q_negative} The negative result is also from \cite{Lim20}: For $N=2^{2k-1}$, there exists a $2k\times N$ matrix $A$ for which condition \eqref{1_q_cond} with $\e^q\log_2 {2N}>2$ is satisfied, and the following relation holds for any partition of the form \eqref{partition}:
$$
\max\{\|A(\Omega_1)\|_{(1,q)}, \|A(\Omega_2)\|_{(1,q)}\}=\|A\|_{(1,q)}.
$$

Let $X$ be a normed space. We can define $\|A\|_{(X, q)}$ in the same way as in \eqref{pq_norm_def}.
\om For $q=\infty$ and any matrix $A$ there is no partition into two submatrices with smaller $(X, q)$--norms, since the matrix $A$ has a row $\bv_{\sup}$ such that $\|A\|_{(X, \infty)}=\sup_{\|\bx\|_X\leq 1}\langle \bx, \bv_{\sup}\rangle$, and the norm of the submatrix containing the row $\bv_{\sup}$ is equal to the norm of the matrix $A$.

The following condition is a direct analog of condition \eqref{cond_eucl} on a matrix in the case of arbitrary $q$, $1\leq q<\infty$:
\begin{equation}\label{cond}
  \forall \bx\in\mathbb R^N\ \  \forall i_0\in\left<M\right> \ \ |\langle\bv_{i_0}, \bx\rangle|\leq\varepsilon\left(\sum_{i=1}^M|\langle\bv_i, \bx\rangle|^q\right)^{1/q}.
\end{equation}
Note that for $1\leq q_1\leq q_2$ condition \eqref{cond} with $q=q_2$ is more restrictive than with $q=q_1$.   It may happen that another condition on $A$ will be natural in the sense that it will lead to sharp estimates. We note that condition \eqref{cond} is the $(q,\infty)$ Nikol'skii-type inequality.
Srivastava asked in \cite{Sr} whether there is an analog of the result {\bf M.\ref{Weav}} under condition \eqref{cond} for $q=1$.
\om \label{pointwise_p_q} A partially positive answer to Srivastava's question was
 obtained in \cite{KL19} and then it was generalized to the case $1\leq q<\infty$ in \cite{Lim20}.
Very recently (see \cite{Kos21}) it was strengthened as follows:
 Assume that $A$ satisfies condition \eqref{cond} for some $\e$ and $1 \leq q <\infty$ with $0<\e\leq N^{-1/q}$. Then there exists a partition of the form \eqref{partition} such that for any $\bx\in\mathbb{R}^N$ and $k=1,2$:
$$
\|A(\Omega_k)\bx\|_q^q\leq\gamma\|A\bx\|_q^q, \quad \gamma=\frac{1}{2}+C(q)\e^{q/2}N^{1/2}.
$$
In terms of the quantity $\varphi(A)$, introduced in the beginning of Section \ref{M}, this means that $\varphi(A) \leq C(q)\e^{q/2}N^{1/2}$. Necessity of the condition $0<\e\leq N^{-1/q}$, which was imposed in both   \cite{KL19} and \cite{Kos21}, is an open problem.

\om \label{om_p_q_negative}
The negative part of the answer to Srivastava's question consists of the following theorem from \cite{KL19}: Let $N=2^s, s\in\mathbb{N}$
, and $N^{-1/2}~\leq~ \e~\leq~ 1$. There exists a $2N\times N$ matrix $A=A(N, \e)$ such that estimate \eqref{cond} holds for any $\bx\in\mathbb{R}^N$ and $i_0\in\langle 2N\rangle$ but, nevertheless, the inequality
$$
\max(\|A(\Omega_1)\|_{(2,1)}, \|A(\Omega_2)\|_{(2,1)})\geq \frac{1}{\sqrt{2}}\left(\frac{1}{1+(\e N^{1/2})^{-1}}\right)\|A\|_{(2,1)}
$$
holds for any partition \eqref{partition} (with $M=2N$).
This result shows that \eqref{cond} for $q=1$ does not guarantee existence of a partition of $A$ of the form \eqref{partition} with a straightforward counterpart of property \eqref{MSS_obtained} for $q=1$. For its existence, it is necessary to assume that the parameter $\e$ tends to zero as $N\rightarrow \infty$. This example  was recently generalized by Limonova to the case of $1\leq q<4/3$.

\om Let $X$ be an $N$--dimensional normed space.
One can easily obtain the following corollary of {\bf{M.}\ref{pointwise_p_q}}
: Assume that, for an   $M\times N$ matrix, condition \eqref{cond} is satisfied for some $\e$ and $q$ with $0<\e\leq (\rank(A))^{-1/q}$ and $1\leq q<\infty$. Then there exists a partition of the form \eqref{partition} such that, for $k=1,2$
$$
\|A(\Omega_k)\|_{(X,q)}^q\leq \gamma \|A\|_{(X,q)}^q,
$$
where $\gamma$ is defined in {\bf{M.}\ref{pointwise_p_q}}.

\subsection{Algorithms}
All of the results in Section \ref{M} were initially established as existence results.  Construction of efficient algorithms for finding desirable submatrices is an important problem.  In this regard we note that the crucial difference between \cite{BSS} and \cite{MSS} is that in \cite{BSS} there is a deterministic polynomial-time algorithm. Many researchers apply ideas from \cite{BSS} to obtain algorithms for other problems. For example, Spielman and Srivastava suggested in \cite{SpSr} a short proof of the restricted invertibility theorem that contains a deterministic algorithm.

\section{Connections with other problems}
\label{connect}

 \subsection{Moments of marginals of
high-dimensional distributions}
\label{MM}

Recall that a random vector $\bu$ in $\mathbb{R}^N$
is a measurable function from some probability space
$(F, \mathcal{F}, {\mathbb P})$ with values in $\mathbb{R}^N$.
The distribution of the random vector $\bu$
is a Borel measure $\mu$ on $\mathbb{R}^N$
such that ${\mathbb P}(\bu \in A) = \mu(A)$ for each Borel set $A$.
This also means that for any reasonable function
(e.g. for any bounded continuous function) $\varphi\colon\mathbb{R}^N\to \mathbb{R}$ one has
$$
\mathbb{E}\varphi(\bu) = \int_{\mathbb{R}^N}\varphi(y)\, \mu(dy),
$$
where $\mathbb{E}$ denotes the expectation of a random variable, i.e.,
the Lebesgue integral with respect to the measure ${\mathbb P}$.

Consider a random vector $\bu$ in $\mathbb{R}^N$. The main problem considered in this section
is to
understand how well one can approximate one-dimensional marginals
of the distribution $\bu$ (i.e. one-dimensional images of this distribution) by sampling.
Consider $m$ independent copies $\bu^1, \ldots, \bu^m$
of the vector $\bu$ and for $q\in[1,\infty)$
consider the random variable
$$
V_q(K):= \sup\limits_{y\in K}
\Bigl|\frac{1}{m}\sum_{j=1}^{m}|\langle y, \bu^j\rangle|^q - \mathbb{E}|\langle y, \bu\rangle|^q\Bigr|,
$$
where $K\subset \mathcal{}\mathbb{R}^N$. When $K=B_2$, where
$B_2:=\{y\in \mathbb{R}^N\colon |y|\le1\}$ is the Euclidean ball,
i.e. $|y|:=\sqrt{\langle y, y\rangle}$,
we write $V_q$ in place of $V_q(B_2)$.
One now wants to
estimate the least possible number $m=m(\varepsilon, q, N)$ of sampling copies
of vector $\bu$ for which $V_q(K)\le \varepsilon$ with high probability.
This problem has been extensively studied for the last 20 years
(see \cite{ALPT-J}, \cite{B},
\cite{GM}, \cite{GR}, \cite{GLPT-J},
\cite{Rud2}, \cite{SV}, \cite{Tikh}, \cite{Versh},
\cite{Versh2} and citations therein).

On one hand, if
$$
K:=\{y\in \mathbb{R}^N\colon \mathbb{E}|\langle y, \bu\rangle|^q\le 1\},
$$
then the described problem of the
approximation of one--dimensional mar\-gi\-nals
of the distribution 
can be seen as the special case of the
sampling discretization problem,
when $\mu$ is the distribution of $\bu$ and the subspace
$X_N:=\{\langle y,\cdot\rangle, y\in \mathbb{R}^N\}\subset L_q(\mu)$.
On the other hand, for any $N$-dimensional subspace $X_N$
of $L_q(\Omega, \mu)$ 
one can take an orthonormal basis $u_1,\ldots, u_N$ in $X_N$
and consider the random vector $\bu:=(u_1(x), \ldots, u_N(x))$.
Thus, the general sampling discretization results
can be deduced from the probabilistic question
about
approximation of one-dimensional marginals
of the distribution.

{\bf MM.1.} The most studied case is $q=2$, in which the problem is equivalent
to estimating of the number of samples $\bu^1,\ldots, \bu^m$  sufficient
to approximate the covariance matrix
$\mathbb{E}\bu\otimes \bu$
of a random vector by a sample covariance matrix $\frac{1}{m}\sum\limits_{j=1}^{m}\bu^j\otimes \bu^j$.
This problem is closely related
to random matrix theory
 and the methods involved in studying it essentially depend on the properties
of the distribution.  Probably the first result of this type is the one  from \cite{BGN}
on bounds for the norm of a random matrix with i.i.d. Bernoulli random entries with values $\pm1$.
Similar results have already been discussed in Section \ref{M} where random submatrices of a given matrix
have been studied. The study of a sample covariance matrix was also partially motivated
by the paper \cite{KLS} where
a fast algorithm for calculating the volume of a convex body has been constructed (see also the book \cite{BGVV}).
For this algorithm one needs to calculate the covariance matrix of a random vector
uniformly distributed on a convex set. The authors of \cite{KLS} have shown
that it is enough to take $m=c\frac{N^2}{\varepsilon}$ points to make
$V_2\le \varepsilon$ with high probability.
 J.~Bourgain obtained in \cite{B}
that one can take $m\ge C(\varepsilon)N[\log N]^3$ for the same bound.
Later M.~Rudelson \cite{Rud2} improved  Bourgain's bound
to $m\ge C \frac{N}{\varepsilon^2}[\log\frac{N}{\varepsilon^2}]^2$
and Rudelson's bound was improved to $m\ge C(\varepsilon)N\log N$
in \cite{Pao}.
Actually, in \cite{Rud2} the following general theorem has been proved.

\begin{Theorem}\label{Rud}
There is a constant $C>0$ such that for any random vector $\bu$ in $\mathbb{R}^N$,
for which $|y|^2=\mathbb{E}|\langle y, \bu\rangle|^2$ for every $y\in \mathbb{R}^N$,
one has
$$
\mathbb{E}V_2
\le C\frac{\sqrt{\log N}}{\sqrt m}\bigl(\mathbb{E}\max\limits_{1\le j\le m}|\bu^j|^2\bigr)^{1/2}\cdot
\Bigl(\mathbb{E}\sup\limits_{|y|\le1}
\Bigl|\frac{1}{m}\sum_{j=1}^{m}|\langle y, \bu^j\rangle|^2\Bigr|\Bigr)^{1/2}
$$
where $\bu^1, \ldots, \bu^m$ are independent copies of the random vector $\bu$.
\end{Theorem}

One can easily verify that this theorem implies the following result
(see e.g. the end of the proof of Theorem~3 in \cite{GR}
or the end of the proof of Lemma~3.1 in \cite{Kos}).

\begin{Corollary}\label{CRud}
There is a constant $C>0$ such that for any random vector $\bu$ in $\mathbb{R}^N$,
for which $|y|^2=\mathbb{E}|\langle  y, \bu\rangle|^2$ for every $y\in \mathbb{R}^N$,
one has
$$
\mathbb{E}V_2
\le
C\Bigl(\frac{\log N}{m}\mathbb{E}\max\limits_{1\le j\le m}|\bu^j|^2+
\frac{\sqrt{\log N}}{\sqrt m}\bigl(\mathbb{E}\max\limits_{1\le j\le m}|\bu^j|^2\bigr)^{1/2}
\Bigr)
$$
where $\bu^1, \ldots, \bu^m$ are independent copies of the random vector $\bu$.
\end{Corollary}

For the general problem of sampling
discretization this
corollary means that for any $N$-dimensional subspace $X_N$ of $L_2(\Omega, \mu)$,
such that for every element $f\in X_N$ one has $\|f\|_\infty\le t N^{\frac{1}{2}}\|f\|_2$
(i.e. Condition E holds),
the $\cM(m, 2, \varepsilon)$ theorem holds for every
$m\ge c\frac{Nt^2}{\varepsilon^2}\log N$
for big enough constant $c$.
Indeed, consider the random vector
$\bu:=(u_1(x), \ldots, u_N(x))$, where $u_1, \ldots, u_N$ is an orthonormal in $L_2(\Omega, \mu)$
basis of $X_N$.
We need to estimate $\mathbb{E}\max\limits_{1\le j\le m}|\bu^j|^2$.
We have $|\bu|^2= |u_1(x)|^2+\ldots+|u_N(x)|^2\le t^2N$.
Thus, $\mathbb{E}\max\limits_{1\le j\le m}|\bu^j|^2\le t^2 N$, which implies
that there is a big enough constant $c$
such that for every $m\ge c\frac{Nt^2}{\varepsilon^2}\log N$
one has
$$
\int_{\Omega^m}\sup\limits_{\substack{\|f\|_2\le 1,\\ f\in X_N}}
\Bigl|\frac{1}{m}\sum_{j=1}^{m}|f(x^j)|^2 - \int_{\Omega}|f(x)|^2\, \mu(dx)\Bigr|\,
\mu(dx^1)\ldots \mu(dx^m)
\le \varepsilon.
$$
Thus, there are points $\xi^1,\ldots, \xi^m$ such that
$$
\sup\limits_{\substack{\|f\|_2\le 1,\\ f\in X_N}}
\Bigl|\frac{1}{m}\sum_{j=1}^{m}|f(\xi^j)|^2 - \int_{\Omega}|f(x)|^2\, \mu(dx)\Bigr|
\le \varepsilon,
$$
which implies the
$\cM(m, 2, \varepsilon)$ theorem for the space $X_N$ while $m\ge c\frac{Nt^2}{\varepsilon^2}\log N$.
As we have already mentioned, the case $q=2$
is connected with random matrix theory,
since in this case we approximate the covariance matrix of a random vector
by its random sample. In particular, one can obtain the same result as above
(see \cite{VT159}) using the recent development of random matrix theory
instead of Rudelson's
Theorem \ref{Rud}
(see the discussion in Subsection~\ref{greedy} below).
The further development of the sampling discretization
problem in the case $q=2$ is connected with the deep
result of A. Marcus, D.A. Spielman, and N. Srivastava from \cite{MSS}
and has already been discussed in the introduction and in Section \ref{tech}.

Returning to the problem of approximation by sample covariance matrix,
the initial question of \cite{KLS} about covariance  matrices of uniform distributions
on convex bodies was finally solved in terms of order in \cite{ALPT-J} (see also Chapter~10 in \cite{BGVV})
where the following two more general results were obtained.

\begin{Theorem}\label{log-conc1}
Let $\bu^1, \ldots, \bu^m$ be i.i.d. random vectors,
distributed uniformly on a symmetric convex body in $\mathbb{R}^N$
and assume that the distribution is isotropic, i.e. $|y|^2=\mathbb{E}|\langle y, \bu^1\rangle|^2$
for every $y\in \mathbb{R}^N$
(or vectors are distributed according to any isotropic logarithmically concave measure).
Then for every $q\ge 2$ and for every $\varepsilon\in(0, 1)$
there is a constant $C(\varepsilon, q)$
such that for any $m\ge C(\varepsilon, q) N^{\frac{q}{2}}$
one has $V_q\le \varepsilon$ with probability at least {$1-e^{-c_p\sqrt{N}}$}
(where the constant $c_q>0$ depends only on $q$).
\end{Theorem}

\begin{Theorem}\label{log-conc2}
Let $\bu^1, \ldots, \bu^m$ be i.i.d. random vectors,
distributed uniformly on a symmetric convex body in $\mathbb{R}^N$
and assume that the distribution is isotropic, i.e. $|y|^2=\mathbb{E}|\langle y, \bu^1\rangle|^2$
for every $y\in \mathbb{R}^N$
(or vectors are distributed according to any isotropic logarithmically concave measure).
Then for every $q\in [1,2)$ and for every $\varepsilon\in(0, 1)$
there is a constant $C(\varepsilon, q)$
such that for any $m\ge C(\varepsilon, q) N$, $m\le e^{\sqrt{N}}$
one has $V_q\le \varepsilon$ with probability at least {$1-e^{-c\sqrt{N}}$}
(where $c>0$ is an absolute constant).
\end{Theorem}

Let now $\Omega$ be a symmetric convex set endowed with the normalized Lebesgue
measure $\mu=\frac{1}{|\Omega|}\lambda_{\Omega}$.
There is a linear isomorphism $T$ such that the random vector
$\bu:=T\bx$ is isotropic, where $\bx=(x_1, \ldots, x_N)$.
Thus,
we can apply Theorems \ref{log-conc1} and \ref{log-conc2} with this vector $\bu$.
In particular, for $q>2$, by Theorem \ref{log-conc1}, there are $m\le C(\varepsilon, q)N^{q/2}$ points
$\xi^1,\ldots, \xi^m\in \Omega$ such that
\begin{equation}\label{eqb}
\Bigl|\int_{\Omega}|\langle \by, T\bx\rangle|^q\, \mu(d\bx)
- \frac{1}{m}\sum\limits_{j=1}^m|\langle \by, T\xi^j\rangle|^q\Bigr|
\le \varepsilon\Bigl(\int_{\Omega}|\langle \by, T\bx\rangle|^2\, \mu(d\bx)\Bigr)^{q/2}.
\end{equation}
This bound implies $\cM(m, q, \varepsilon)$ theorem
for the space $X_N$ of all linear functionals on $\mathbb{R}^N$ with
$m\le C(\varepsilon, q)N^{q/2}$ points for $q>2$,
since the $L_2$ norm is bounded by the $L_q$ norm.
Similarly, for $q\in [1, 2)$ we can apply Theorem \ref{log-conc2}
and get the same bound \eqref{eqb} with $m\le C(\varepsilon, q)N$
points $\xi^1,\ldots, \xi^m$.
It is known that all the $L_q$ norms are equivalent on a space of all linear
functions with respect to a logarithmically concave measure, thus
$$
\Bigl(\int_{\Omega}|\langle y, Tx\rangle|^2\, \mu(dx)\Bigr)^{q/2}
\le c(q)\int_{\Omega}|\langle y, Tx\rangle|^q\, \mu(dx)
$$
for some number $c(q)$, that depends only on $q$.
Thus, for $q\in[1,2)$ we get $\cM(m, q, \varepsilon)$ theorem
for the space of all linear functionals with
$m\le C(q)N$ points.
An interesting open question is whether similar results are true for algebraic polynomials
of  fixed degree on convex domains $\Omega$.

In the case $q=2$ the further study
(see \cite{Versh2}, \cite{GLPT-J}, \cite{SV},  \cite{MP-1}, \cite{MP-2}, \cite{Tikh})
of the problem of
approximation of the covariance matrix
by a sample covariance matrix
has continued for random vectors $\bu$ under $2+\varepsilon$ moment assumption,
i.e. it is assumed that
$$
\bigl(\mathbb{E}|\langle y, \bu\rangle|^p\bigr)^{1/p}
\le B\bigl(\mathbb{E}|\langle y, \bu\rangle|^2\bigr)^{1/2}\quad \forall y\in \mathbb{R}^N
$$\
for some $p>2$.

\vskip .1in

{\bf MM.2.} We now discuss the case $q\ne 2$
and start with the results from \cite{Versh}
for random vectors under the  assumption of higher order
integrability of one-dimensional marginals,
i.e. it is assumed that
\begin{equation}\label{ma}
\bigl(\mathbb{E}|\langle y, \bu\rangle|^p\bigr)^{1/p}
\le B\bigl(\mathbb{E}|\langle y, \bu\rangle|^2\bigr)^{1/2}\quad \forall y\in \mathbb{R}^N
\end{equation}
for some $p>\max\{q,2\}$.

For such random vectors the following two results were obtained in~\cite{Versh}.

\begin{Theorem}\label{T-Versh-1} Let $q>2$ and
let $\bu$ be a random vector in $\mathbb{R}^N$ such that
$|\bu|\le t\sqrt{N}$ a.s.,
$|y|^2=\mathbb{E}|\langle y, \bu\rangle|^2$
for every $y\in \mathbb{R}^N$,
and estimate \eqref{ma}
holds with $p=4q$. Then
for every $\varepsilon\in(0, 1)$ and for every $\delta\in(0, 1)$
there is a number
$C:=C(t, B, q, \varepsilon, \delta)$ such that
$V_q\le \varepsilon$ with probability at least $1-\delta$ for every $m\ge CN^{\frac{q}{2}}$.
\end{Theorem}

\begin{Theorem}\label{T-Versh-2} Let $q\in[1, 2)$ and
let $\bu$ be a random vector in $\mathbb{R}^N$ such that
$|\bu|\le t\sqrt{N}$ a.s.,
$|y|^2=\mathbb{E}|\langle y, \bu\rangle|^2$
for every $y\in \mathbb{R}^N$, and estimate \eqref{ma}
holds with $p\ge4q$, $p>4$. Then
for every $\varepsilon\in(0, 1)$ and for every $\delta\in(0, 1)$
there is a number
$C:=C(t, B, q, p, \varepsilon, \delta)$ such that
$V_q\le \varepsilon$ with probability at least $1-\delta$ for every $m\ge CN$.
\end{Theorem}

The former of these two theorems is Theorem 1.1 in \cite{Versh}
and the latter is formulated in \cite{Versh}
in the remark after Theorem 1.1.

We now reformulate these theorems
with regards to the general problem of sampling discretization
for a subspace $X_N$ of $L_q(\Omega, \mu)$.

\begin{Corollary}\label{C-ma-1}
Let $q>2$ and
let $X_N$ be an $N$-dimensional subspace of $L_q(\Omega, \mu)$.
Assume that $\|f\|_\infty\le t N^{\frac{1}{2}}\|f\|_2$ for every $f\in X_N$
and assume that
\be\label{sr1}
\|f\|_{4q}\le B\|f\|_2
\ee
for some constant $B>0$.
Then for every $\varepsilon\in(0,1)$
there is a constant $C:=C(t, B, q, \varepsilon)$
such that
$\cM(m, q, \varepsilon)$ theorem holds for every
$m\ge CN^{\frac{q}{2}}$.
\end{Corollary}

\begin{Corollary}\label{C-ma-2}
Let $q\in[1,2)$ and
let $X_N$ be an $N$-dimensional subspace of $L_q(\Omega, \mu)$.
Assume that $\|f\|_\infty\le t N^{\frac{1}{2}}\|f\|_2$ for every $f\in X_N$
and assume that
\be\label{sr2}
\|f\|_p\le B\|f\|_q
\ee
for some $p\ge 4q$, $p>4$, and for some $B>0$.
Then for every $\varepsilon\in(0,1)$
there is a constant $C:=C(t, B, q, p, \varepsilon)$
such that
$\cM(m, q, \varepsilon)$ theorem holds for every
$m\ge CN$.
\end{Corollary}

For example, conditions \eqref{sr1} and \eqref{sr2} are true for the space of all
polynomials of a fixed degree on a convex domain $\Omega$ endowed
with the normalized Lebesgue measure.

 In particular, Corollary \ref{C-ma-2} implies that,
with the extra assumption \eqref{sr2},
one does not need the logarithmic oversampling. Therefore, Corollary \ref{C-ma-2} allows us to slightly improve bound (\ref{A17}) from {\bf D.16}
by imposing a rather strict assumption \eqref{sr2}.

Both corollaries follow from Theorems \ref{T-Versh-1} and \ref{T-Versh-2}
by consideration of the random vector $\bu(x):=(u_1(x), \ldots, u_N(x))$,
where $u_1, \ldots, u_N$ is an orthonormal in $L_2(\Omega, \mu)$
basis in $X_N$.
Indeed, for any $f\in X_N$ there is $\by=(y_1, \ldots, y_N)$
such that $f=\sum\limits_{j=1}^Ny_ju_j$ and
$$
\|f\|_p^p = \int_{\Omega} |f(x)|^p\, \mu(dx)
=  \int_{\Omega} \bigl|\sum_{j=1}^{N}y_ju_j(x)\bigr|^p\, \mu(dx)
= \mathbb{E} |\langle \by, \bu\rangle|^p.
$$
Thus, for $q\in[1, 2)$ by Theorem \ref{T-Versh-2}
for every $\varepsilon>0$ there is a positive number
$C:=C(t, B, q, \varepsilon)$ such that for any $m\ge CN$
there are points $\xi^1, \ldots, \xi^m$ for which
$$
\sup\limits_{\substack{\|f\|_2\le 1,\\ f\in X_N}}
\Bigl|\frac{1}{m}\sum_{j=1}^{m}|f(\xi^j)|^q - \|f\|_q^q\Bigr|\le \frac{1}{B^q}\varepsilon.
$$
Since $p>2$
one has $\|f\|_2^q\le \|f\|_p^q\le B^q\|f\|_q^q$ and
$X_N\in \cM(m, q, \varepsilon)$.

For $q>2$ by Theorem \ref{T-Versh-1}
for every $\varepsilon>0$ there is a positive number
$C:=C(t, B, q, \varepsilon)$ such that for any $m\ge CN^{\frac{q}{2}}$
there are points $\xi^1, \ldots, \xi^m$ for which
$$
\sup\limits_{\substack{\|f\|_2\le 1,\\ f\in X_N}}
\Bigl|\frac{1}{m}\sum_{j=1}^{m}|f(\xi^j)|^q - \|f\|_q^q\Bigr|\le \varepsilon.
$$
Since $q>2$
one has $\|f\|_2^q\le \|f\|_q^q$ and
$X_N\in \cM(m, q, \varepsilon)$.

{
We point out
that using a different type of approximation method one could improve the
bound for the number of points in Theorem \ref{T-Versh-1}
(see {the} recent paper \cite{Mend}).
}

\vskip .1in

{\bf MM.3.} The {last} result of this subsection is
due to O.~Guedon and M.~Rudelson (see \cite[Theorem 3]{GR}).

Let $K$ be a symmetric convex body. Then
$\|\cdot\|_K$ denote the norm in which $K$ is the unit ball, i.e.
$\|z\|_K:=\inf\{s\ge0\colon s^{-1}z\in K\}$.
Let $p\ge 2$. We recall that the symmetric convex body $K$
in $\mathbb{R}^N$
is called $p$-convex with constant $\eta$ if
$$
\Bigl\|\frac{z_1+z_2}{2}\Bigr\|_K\le 1-\eta\|z_1-z_2\|_K^p
$$
for all vectors $z_1, z_2\in K$.

\begin{Theorem}\label{T-GR}
Let $K\subset \mathbb{R}^N$ be a symmetric convex
body of radius $D$.
Assume that $K$ is $p$-convex with constant $\eta$
for some $p\ge 2$. Assume that $q\ge p$.
Then for any random vector $\bu$ one has
$$\mathbb{E}V_q(K)\le A^2+A\sqrt{B},$$
where
$$
A:=C^q\eta^q\frac{(\log m)^{1-1/p}}{\sqrt{m}}D^{q/2}
\bigl(\mathbb{E}\max\limits_{1\le j\le m}|\bu^j|^q \bigr)^{1/2},\quad
B:= \sup\limits_{y\in K}\mathbb{E}|\langle y,\bu\rangle|^q.
$$
\end{Theorem}

This theorem implies the following result for the general problem of sampling discretization.

\begin{Corollary}\label{C-GR}
Let $q\ge 2$ and let $X_N$ be an $N$-dimensional
subspace of $L_q(\Omega, \mu)$
such that
$$
\|f\|_\infty\le t\sqrt{N}\|f\|_2\quad  \forall f\in X_N.
$$
Then for every $\varepsilon\in(0, 1)$
there is a number $C:=C(t, q, \varepsilon)$
such that for each $m\ge CN^{\frac{q}{2}}[\log N]^{2-\frac{2}{q}}$
one has $X_N\in \cM(m, q, \varepsilon)$.
\end{Corollary}

Indeed, we take $K=B_q:=\{f\in X_N\colon \|f\|_q\le 1\}$.
It is known that $B_q$
is $q$-convex with some constant $\eta=\eta(q)$ for $q\ge 2$.
We again consider the random vector $\bu:=(u_1(x), \ldots, u_N(x))$,
where $u_1, \ldots, u_N$ is an $L_2(\Omega, \mu)$-orthonormal   
basis in $X_N$,
and again $|\bu|\le t\sqrt{N}$.
Moreover,  $\|f\|_2\le\|f\|_q$, i.e. $D=1$.
Thus,
$$
A\le C(q)t^{q/2}\frac{(\log m)^{1-1/q}}{\sqrt{m}}\bigl(\sqrt{N}\bigr)^{q/2},\quad
B=1.
$$
The discretization result now follows by the standard argument
when we take $m\ge C(t, q, \varepsilon) N^{\frac{q}{2}}[\log N]^{2-\frac{2}{q}}$
for large enough constant $C(t, q, \varepsilon)$
such that $A\le \frac{1}{4}\varepsilon$.

Lewis's change of density theorem (see \cite[Theorem 2.1]{SZ}),
which we discuss further, combined with Corollary \ref{C-GR},
gives the following statement about weighted sampling discretization.

\begin{Corollary}\label{C-W-GR}
Let $q\ge 2$ and let $X_N$ be an $N$-dimensional
subspace of $L_q(\Omega, \mu)\cap C(\Omega)$.
Then for every $\varepsilon\in(0, 1)$
there is a number $C:=C(q, \varepsilon)$
such that for each $m\ge CN^{\frac{q}{2}}[\log N]^{2-\frac{2}{q}}$
one has $X_N\in \cM^w(m, q, \varepsilon)$.
\end{Corollary}

Actually, one always has $X_N\in \cM^w(m, q, \varepsilon)$
for each $m\ge CN^{\frac{q}{2}}[\log N]$.
For example, this follows from the results on  {the} embedding problem
from Subsection~\ref{Em} or from the
slightly improved version of
Theorem~\ref{T-GR} from \cite{Kos} (see Corollaries 4.4 and 4.5 and Remark 4.6 there).

\subsection{Embedding of finite-dimensional subspaces of $L_q$ into $\ell_q^m$}
\label{Em}

Let $q\in[1, \infty)$ and $\varepsilon>0$ be fixed.
The main problem of this subsection
is to understand what is the smallest possible
integer
$m:=m_q(\varepsilon, N)$ such that for every
$N$-dimensional subspace $X_N$ of $L_q([0, 1])$, $1\le q<\infty$
there is an $N$-dimensional subspace $Y_N$ of $\ell_q^m$
with $d(X_N, Y_N)\le 1+\varepsilon$.
Here $\ell_q^m$ is $\mathbb{R}^m$ endowed with
the norm $\|y\|_q:=\bigl(\sum\limits_{j=1}^m|y_j|^q\bigr)^{1/q}$
and $d(X_N, Y_N)$ is the Banach–Mazur distance
between $X_N$ and $Y_N$,
which
is the infimum of $\|T\|\cdot\|T^{-1}\|$
over all linear isomorphisms $T$ between $X_N$ and $Y_N$,
where $\|\cdot\|$ is an operator norm
(actually, the logarithm of $d$ is a distance).

This important and long standing problem has been extensively studied from the early $1980$-s
to the present by many mathematicians
(see \cite{BLM}, \cite{Sche3}, \cite{Sche4}, \cite{Ta2}, \cite{Ta3}, \cite{Ta1},  \cite{SZ}, \cite{Zvav} and citations therein).
There is also a very nice expository paper by W.B.~Johnson and G.~Schechtman \cite{JS}
where one could find a profound discussion of the problem
along with some historical comments and more literature on the subject.

The main goal of this subsection is to present main ideas and methods of the theory behind
the stated problem
which {can also be} applied for various other problems. The important connection of this
problem with the discretization problem from Section \ref{A} is through the  
empirical method, which proposes to study the problem of embedding through discretization.
We point out that this part of our paper is based
on the above mentioned article \cite{JS}
and two books' sections \cite[Section 16.8]{Ta1}
and \cite[Section 15.5]{LedTal}.

{\bf Em.1.} First of all we formulate the known bounds for $m_q(\varepsilon, N)$.
For $q>2$ one has $m_q(\varepsilon, N)\le C(q, \varepsilon)N^{\frac{q}{2}}\log N$
(see \cite{BLM}).
For $q\in(1, 2)$ one has $m_q(\varepsilon, N)\le C(\varepsilon)N[\log N][\log\log N]^2$
(see \cite{Ta3}).
For $q=1$ one has $m_1(\varepsilon, N)\le C(\varepsilon)N\log N$ (see \cite{Ta2}).
The problem when $q\in(0, 1)$ has also been studied (see \cite{SZ} and \cite{Zvav})
and the known bound is the same as for $q\in(1, 2)$ (see \cite{Zvav}).

{\bf Em.2.} We now discuss the main technique, which was initiated in \cite{Sche3},
then developed in \cite{BLM} and \cite{Ta2} and since then
has been often applied in the study of the embedding problem.
We point out that presented in this paper ideas do not provide the above best known results on their own
and one needs some other technical steps to obtain the best known bounds
for $m_q(\varepsilon, N)$.

First of all, one can always assume that
the subspace $X_N$ of $L_q([0, 1])$ is at a distance less than $1+\varepsilon$
from some $N$-dimensional subspace $X_N^0$ of $L_q(\Omega_M, \mu)$ for some large number $M$,
where $\Omega_M=\{x_1, \ldots, x_M\}$ is a discrete point set with cardinality $M$
and $\mu$ is some probability measure on $\Omega_M$. We will write $\mu_j$
in place of $\mu(\{x_j\})$.
The first idea is to find a good replacement for the measure $\mu$,
i.e. to find a measure $\nu_j=\varrho(x_j)\mu_j$, such that
the subspace $(X_N^0, \|\cdot\|_{L_q(\Omega_M, \mu)})$
is isometric to some subspace $(Y_N^0, \|\cdot\|_{L_q(\Omega_M, \nu)})$
and this new subspace $Y_N^0$ possesses some good properties.
This can be done via the Lewis change
of density theorem (see \cite{Lew}).
We formulate this theorem as it was stated in \cite{SZ}, where the result of Lewis
was extended to the case $q\in(0, 1)$.

\begin{Theorem}\label{C-D-Lew}
Let $X_N$ be an  $N$-dimensional subspace of $L_q(\Omega, \mu)$, where $0<q<\infty$.
Then there is a basis $f_1, \ldots, f_N$ in $X_N$ such that
$$
\int_{\Omega}\Bigl(\sum_{k=1}^{N}f_k^2\Bigr)^{\frac{q-2}{2}} f_if_j\, d\mu = \delta_{i, j}\quad
\forall i, j\in\{1, \ldots, N\}.
$$
In particular, there is a probability measure $\nu$ on $\Omega$ and a space
$Y_N$ in $L_q(\Omega, \nu)$  which is
isometric to $X_N$ and which
has {a basis 
    $u_1, \ldots, u_N$ that is orthonormal in $L_2(\Omega, \nu)$ such that
    $\sum\limits_{k=1}^{N}u_k^2(x)=N\quad \forall x\in \Omega.$}
\end{Theorem}

{\bf Em.3.}  Lewis's theorem implies that
there is a new probability measure $\nu=\varrho\mu$ on $\Omega_M$
such that $(X_N^0, \|\cdot\|_{L_q(\Omega_M, \mu)})$
is isometric to the subspace {
    $$Y_N^0:=\{f'\colon \Omega_M\to \mathbb{R}\colon f'(x_j)=\varrho(x_j)^{-1/q}f(x_j), f\in X_N^0\}$$
    }
of
$L_q(\Omega_M, \nu)$
and  {there is a basis
$u_1, \ldots, u_N$ of $Y_N^0$ that is orthonormal in  $L_2(\Omega_M, \nu)$}  such that
$\sum\limits_{k=1}^{N}u_k^2(x)=N$ for each $x\in \Omega_M$.
Moreover, one can split the atoms $\nu_i$, which are greater than $\frac{2}{M}$,
into $[\nu_i\frac{M}{2}]+1$ equal pieces.
The cardinality $M'$ of
the new set $\Omega'_{M'}=\{x_1', \ldots, x'_{M'}\}$
after the splitting is {not greater} than $\frac{3M}{2}$.
Let $\mu'$ be the new measure on $\Omega'_{M'}$
which appears from $\nu$ after
the splitting of atoms, $\mu'_j\le \frac{2}{M}$.
Let $X_N'$ be a subset of $L_q(\Omega'_{M'}, \mu')$
which appears from $Y_N^0$ after the splitting of atoms.

{\bf Em.4.} The next idea is to
``split'' an $L_q(\Omega'_{M'}, \mu')$ norm on $X_N'$ into two almost equal parts,
i.e. one wants to find a splitting of the set $\Omega'_{M'}$
into two disjoint parts $W_1$ and $W_2$
such that
$$
\bigl(\tfrac{1}{2}-
\varepsilon_{N, M'}\bigr)\sum_{j=1}^{M'}\mu'_j|f'(x'_j)|^q\le
\sum_{x_j\in W_i}\mu'_j|f'(x'_j)|^q\le
\bigl(\tfrac{1}{2}+\varepsilon_{N, M'}\bigr)\sum_{j=1}^{M'}\mu'_j|f'(x'_j)|^q
$$
for each $f'\in X_N'$ and $i\in\{1, 2\}$. {We point out that exactly this type of  splitting
is one of the main results of A. Marcus, D.A. Spielman, and N.~Srivastava
from paper \cite{MSS}.}

Clearly,
one of the sets $W_1$ or $W_2$ is of cardinality {not greater}
than $\frac{1}{2}M'$, i.e. of cardinality bounded by $\frac{3}{4}M$
and exactly this set we denote $\Omega^1_{M_1}$, $M_1\le \frac{3}{4}M$.
If we return to the initial measure $\mu$ and the subspace $X_N$,
the above bound actually reads as
$$
\bigl(\tfrac{1}{2}-\varepsilon_{N, M'}\bigr)\sum_{j=1}^{M}\mu_j|f(x_j)|^q
\le
\sum_{y_j\in \Omega^1_{M_1}}\mu^1_j|f(y_j)|^q\le
\bigl(\tfrac{1}{2}+\varepsilon_{N, M'}\bigr)\sum_{j=1}^{M}\mu_j|f(x_j)|^q
$$
for all $f\in X_N$,
where $\{\mu^1_j\}$ are new weights (new measure) on the new domain
$\Omega^1_{M_1}$ with cardinality $M_1\le \lambda M$, $\lambda=\frac{3}{4}$.
Then one iterates the described procedure and, while
$\prod\limits_{k=0}^{k_0}
\Bigl(\frac{1+\varepsilon_{N, \lambda^k M}}{1-\varepsilon_{N, \lambda^k M}}\Bigr)^{1/q}\le 1+\varepsilon$,
the initial subspace $X_N$ is at a distance at most $1+\varepsilon$
from an $N$-dimensional subspace in $\ell_q^{m}$ with $m\le \lambda^{-k_0-1}M$.

The  above iteration procedure was applied to the embedding problem
in \cite{Ta2} { but a similar idea of iteration
has also appeared earlier in the work of Bourgain, Lindenstrauss, and Milman 
\cite{BLM}
and in the work of Lunin \cite{Lun}
where he studied the problem of finding a submatrix with
the smallest possible norm in a given matrix (see more on this in Section \ref{M}).
}

{\bf Em.5.} The first bound in {\bf Em.4.} is equivalent to the bound
$$
\Bigl|\sum_{x_j\in W_1}\mu'_j|f'(x'_j)|^q-\sum_{x_j\in \Omega'_{M'}\setminus W_1}\mu'_j|f'(x'_j)|^q\Bigr|
\le 2\varepsilon_{N, M'}\sum_{j=1}^{M'}\mu'_j|f'(x'_j)|^q\quad \forall f'\in X_N'.
$$
This inequality means that there is a choice of signs
$\{\varepsilon_1^0, \ldots, \varepsilon_{M'}^0\}$
such that
$$
\Bigl|\sum_{j=1}^{M'}\varepsilon_j^0\mu'_j|f(x'_j)|^q\Bigr|
\le 2\varepsilon_{N, M'}\sum_{j=1}^{M'}\mu'_j|f'(x'_j)|^q\quad \forall f'\in X_N'
$$
or, equivalently,
$$
\sup\limits_{\substack{\|f'\|_{L_q(\Omega'_{M'}, \mu')}\le 1\\ f'\in X_N'}}
\Bigl|\sum_{j=1}^{M'}\varepsilon_j^0\mu'_j|f'(x'_j)|^q\Bigr|\le
2\varepsilon_{N, M'}.
$$
Clearly, there is a choice of signs $\{\varepsilon_1^0, \ldots, \varepsilon_{M'}^0\}$
such that
$$
\sup\limits_{\substack{\|f'\|_{L_q(\Omega'_{M'}, \mu')}\le 1\\ f'\in X_N'}}
\Bigl|\sum_{j=1}^{M'}\varepsilon_j^0\mu'_j|f'(x'_j)|^q\Bigr|\le
\mathbb{E}\sup\limits_{\substack{\|f'\|_{L_q(\Omega'_{M'}, \mu')}\le 1\\ f'\in X_N'}}
\Bigl|\sum_{j=1}^{M'}\varepsilon_j\mu'_j|f'(x'_j)|^q\Bigr|,
$$
where $\mathbb{E}$ stands for the average over all possible choice of signs
$\{\varepsilon_1, \ldots, \varepsilon_{M'}\}$. Thus,
one interested in the bounds for the expectation
$$
\mathbb{E}\sup\limits_{\substack{\|f'\|_{L_q(\Omega'_{M'}, \mu')}\le 1\\ f'\in X_N'}}
\Bigl|\sum_{j=1}^{M'}\varepsilon_j\mu'_j|f'(x'_j)|^q\Bigr|
$$
and the embedding problem reduces to finding good bounds
for the expectation of the supremum of a random process over a convex set.
 These bounds are similar to the one discussed in Subsection \ref{MM}
and can be studied by the same techniques
(e.g. by the chaining technique).

In the case $p\in[1, 2)$ the bounds for the expectation above
have been obtained by Talagrand
(see \cite[Theorem 15.12]{LedTal} and \cite[Theorem 16.8.2]{Ta1} or \cite{Ta2} and \cite{Ta3}).

We point out that in the case $p>2$ the paper \cite{BLM}
uses a different approach.

{\bf Em.6.}  As it was   by G. Schechtman brought to our attention that these
 methods of proofs  
may actually provide the same bounds for the number of points
in the corresponding problems of sampling discretization with weights.

\subsection{Sparse approximation}
\label{greedy}

We have explained in Section \ref{Int} (see (\ref{Gb})) how sampling discretization of the $L_2$ norm is connected with the $m$-term approximation of the identity matrix $I$ with respect to the system $\{G(x)\}_{x\in \Omega}$ of rank one matrices.
Namely, for a set of points $\xi^k\in \Omega$, $k=1,\dots,m$, and $f=\sum_{i=1}^N b_iu_i$ we have
$$
\frac{1}{m}\sum_{k=1}^m f(\xi^k)^2 - \int_\Omega f(x)^2 d\mu = {\mathbf b}^T\left(\frac{1}{m}\sum_{k=1}^m G(\xi^k)-I\right){\mathbf b},
$$
where ${\mathbf b} = (b_1,\dots,b_N)^T$ is the column vector. Therefore,
\be\label{5.8'}
\left|\frac{1}{m}\sum_{k=1}^m f(\xi^k)^2 - \int_\Omega f(x)^2 d\mu \right| \le
\left\|\frac{1}{m}\sum_{k=1}^m G(\xi^k)-I\right\|\|{\mathbf b}\|_2^2,
\ee
where $\|\cdot\|$ is the operator norm (spectral norm) of a matrix. Then the following deep result (see \cite{Tro12}, Theorem 1.1) on random matrices can be used.

\begin{Theorem}\label{T5.3} Consider a finite sequence $\{T_k\}_{k=1}^m$ of independent, random, self-adjoint matrices with dimension $N$. Assume that
each random matrix is positive semi-definite and satisfies
$$
\lambda_{\max}(T_k) \le R\quad  \text{almost surely}.
$$
{Note that $\lambda_{\max}$ and $\lambda_{\min}$ denote the maximal and the minimal eigenvalues. }
Define
$$
s_{\min} := \lambda_{\min}\left(\sum_{k=1}^m \bE(T_k)\right) \quad \text{and}\quad
s_{\max} := \lambda_{\max}\left(\sum_{k=1}^m \bE(T_k)\right).
$$
Then
$$
\bP\left\{\lambda_{\min}\left(\sum_{k=1}^m T_k\right) \le (1-\eta)s_{\min}\right\} \le
N\left(\frac{e^{-\eta}}{(1-\eta)^{1-\eta}}\right)^{s_{\min}/R}
$$
for $\eta\in[0,1)$ and for $\eta\ge 0$
$$
\bP\left\{\lambda_{\max}\left(\sum_{k=1}^m T_k\right) \ge (1+\eta)s_{\max}\right\} \le
N\left(\frac{e^{\eta}}{(1+\eta)^{1+\eta}}\right)^{s_{\max}/R}.
$$
\end{Theorem}

It is shown in \cite{VT159} that Theorem \ref{T5.3} implies the following sampling discretization result: Let $\{u_i\}_{i=1}^N$ be an orthonormal system, satisfying Condition~E. Then
\begin{equation*}
X_N=\sp\{u_1,\dots,u_N\} \in \cM(m,2,\e)\quad \text{provided}\quad m\ge C\frac{t^2}{\e^2}N\log N .
\end{equation*}
 We only formulated one result -- Theorem \ref{T5.3} -- from the theory of random matrices which is useful in sampling discretization. We refer the reader for other results in this direction to the papers
\cite{Rud2}, \cite{Tro12}, \cite{Oliv}, \cite{Rau}, \cite{GPR}, and \cite{MPa}. Also, we mention a very recent paper
\cite{MUl} where results from the theory of random matrices were applied to sampling discretization and recovery.

Inequality (\ref{5.8'}), shows that the Marcinkiewicz-type
discretization theorem in $L_2$ is closely related to approximation of the identity matrix $I$ by an $m$-term approximant of the form $\frac{1}{m}\sum_{k=1}^m G(\xi^k)$ in the operator norm from $\ell^N_2$ to $\ell^N_2$ (spectral norm).
Therefore, we can consider the following sparse approximation problem:
Assume that the system $\{u_i(x)\}_{i=1}^N$ satisfies Condition~E and consider
the dictionary
$$
\D^u := \{g_x\}_{x\in\Omega},\quad g_x:= G(x)(Nt^2)^{-1},\quad G(x):=[u_i(x)u_j(x)]_{i,j=1}^N.
$$
 Let the Hilbert space
$H$ be the closure in the Fr{\"o}benius norm of  the space $\sp\{g_x, x\in\Omega\}$ with the inner product generated by the Fr{\"o}benius norm: for $A=[a_{i,j}]_{i,j=1}^N$ and
$B=[b_{i,j}]_{i,j=1}^N$
$$
\<A,B\> = \sum_{i,j=1}^N a_{i,j}b_{i,j}
$$
in case of real matrices (with standard modification in case of complex matrices).

By known results from greedy approximation (see, for instance, \cite{VTbook}, p.90, Theorem 2.15) we obtain that, for any $m\in \N$, we can constructively find (by the Relaxed Greedy Algorithm)
points $\xi^1,\dots,\xi^m$ such that
\begin{equation}\label{5.12}
\left\|\frac{1}{m}\sum_{k=1}^m G(\xi^k)-I\right\|_F \le 2Nt^2 m^{-1/2}.
\end{equation}
Taking into account the inequality $\|A\|\le \|A\|_F$   and  (\ref{5.8'}) we get from (\ref{5.12})
the following proposition.

\begin{Proposition}\label{P5.1} Let $\{u_i\}_{i=1}^N$ be an orthonormal system, satisfying Condition E. Then there exists a constructive set $\{\xi^j\}_{j=1}^m \subset \Omega$ with $m\le C(t)N^2$
such that for any $f=\sum_{i=1}^N c_iu_i$ we have
$$
\frac{1}{2}\|f\|_2^2 \le \frac{1}{m} \sum_{j=1}^m f(\xi^j)^2 \le \frac{3}{2}\|f\|_2^2.
$$
\end{Proposition}

We now comment on a  breakthrough result by J. Batson, D.A. Spielman, and N. Srivastava \cite{BSS} on the weighted sampling discretization. We formulate their result in our notation. Let as above $\Omega_M=\{x^j\}_{j=1}^M$ be a discrete set with the probability measure $\mu(x^j)=1/M$, $j=1,\dots,M$. Assume that
$\{u_i(x)\}_{i=1}^N$ is a real orthonormal on $\Omega_M$ system. Then for any
number $d>1$ there exist a set of weights $w_j\ge 0$ such that $|\{j: w_j\neq 0\}| \le dN$ so that for any $f\in \sp\{u_1,\dots,u_N\}$ we have
$$
\|f\|_2^2 \le \sum_{j=1}^M w_jf(x^j)^2 \le \frac{d+1+2\sqrt{d}}{d+1-2\sqrt{d}}\|f\|_2^2.
$$
The proof of this result is based on a delicate study of the $m$-term approximation of the identity matrix $I$ with respect to the system
$\D := \{G(x)\}_{x\in \Omega}$, $G(x):=[u_i(x)u_j(x)]_{i,j=1}^N$ in the spectral norm. The authors of \cite{BSS} control the change of the maximal and minimal eigenvalues of a matrix, when they add a rank one matrix of the form $wG(x)$.
Their proof provides an algorithm for construction of the weights $\{w_j\}$.
In particular, this implies that
$$
X_N(\Omega_M) \in \cM^w(m,2,\e)\quad \text{provided} \quad m \ge CN\e^{-2}
$$
with large enough $C$.

 \subsection{Supervised learning theory}
\label{learn}

In Section \ref{Int} we have discussed one of the settings of learning theory.  In this subsection we present some of the well known results from supervised learning theory, comment on { important techniques},  and discuss connection with the sampling discretization. The chaining technique and concentration of measure inequalities are the foundation blocks of both learning theory and sampling discretization. We begin with a detailed discussion of the interplay between
the use of expectation and the probability distribution function in the error analysis.
  For the reader's convenience we begin our discussion with an application of probabilistic methods in numerical integration.   We recall the classical Monte Carlo method from numerical integration. Let $\Omega$ be a bounded domain of $\R^d$.
Consider a real function $f\in L_2(\Omega):=L_2(\Omega,\mu)$ with respect to a probability measure $\mu$. Define $\Omega^m:= \Omega\times\dots\times\Omega$ and $\mu^m:=\mu \times\dots\times \mu$. For $\bx^j\in \Omega$ denote $\bz:= (\bx^1,\dots,\bx^m)\in \Omega^m$ and for $g\in L_1(\Omega^m, \mu^m)$
$$
\bbE(g):= \int_{\Omega^m}g(\bz)d\mu^m.
$$
Then it is known that for $f\in L_2(\Omega,\mu)$ we have
\be\label{MC}
\bbE\left(\left(\int_\Omega f d\mu - \frac{1}{m}\sum_{j=1}^m f(\bx^j)\right)^2\right) \le \frac{\|f\|_2^2}{m}.
\ee
In particular, inequality (\ref{MC}) guarantees existence of a cubature formula $Q_m(\cdot,\xi)$, with $\xi:=\{\xi^j\}_{j=1}^m$   such that
\be\label{B1}
\left|\int_\Omega fd\mu -Q_m(f,\xi)\right| \le m^{-1/2}\|f\|_2.
\ee
Certainly, the existence of a good cubature formula for an individual function is a simple fact. {For a  function $f$ that is continuous on $\Omega$} we can always find a point $\xi^1$ and a weight $\la_1$ such that $\int_\Omega fd\mu =\la_1f(\xi^1)$. However, the above probabilistic argument allows us to guarantee a bound similar to (\ref{B1}) with high probability. This observation, in turn, allows us to find a  cubature formula, which is simultaneously good for several functions.
However, the use of expectation does not provide a good enough bound on  {the} probability to guarantee a tight error bound like (\ref{B1}) for many functions. The concentration of measure inequalities, which we formulate momentarily, provide a  very good bound on the probability under some extra assumptions on $f$. Under the condition $\|f\|_\infty \le M$ the Hoeffding inequality (see, for instance, \cite{VTbook}, p.197) gives
\be\label{B2}
\mu^m\left\{\bz: \left|\int_\Omega fd\mu - \frac{1}{m}\sum_{j=1}^m f(\bx^j)\right|\ge \eta\right\}\le 2\exp\left(-\frac{m\eta^2}{8M^2}\right).
\ee
  Bernstein's inequality (see, for instance, \cite{VTbook}, p.198) gives the following bound
under conditions $\|f\|_\infty\le M_\infty$ and $\|f\|_2 \le M_2$
\be\label{B3}
\mu^m\left\{\bz: \left|\int_\Omega fd\mu - \frac{1}{m}\sum_{j=1}^m f(\bx^j)\right|\ge \eta\right\}\le 2\exp\left(-\frac{m\eta^2}{2(M_2^2+2M_\infty\eta/3)}\right).
\ee

The above inequalities (\ref{B2}) and (\ref{B3}) can be used directly for proving the existence of good cubature formulas for function classes containing a finite number of elements. Denote by $|\bW|$  {the} cardinality of a set $\bW$. Assume that for all $f\in \bW$ we have $\|f\|_\infty \le M$. Then, the Hoeffding inequality (\ref{B2}) gives
\be\label{B4}
\mu^m\left\{\bz: \sup_{f\in \bW}\left|\int_\Omega fd\mu - \frac{1}{m}\sum_{j=1}^m f(\bx^j)\right|\le \eta\right\}\ge 1- 2|\bW|\exp\left(-\frac{m\eta^2}{8M^2}\right).
\ee
Thus, when the right hand side of (\ref{B4}) is positive, inequality (\ref{B4}) guarantees the existence of a good cubature formula for the whole class $\bW$.  However, it is clear from relation (\ref{B4}) that such an argument does not work for a class $\bW$, which contains
infinitely many functions. It is known that chaining techniques allow us to overcome this problem.

We now concentrate on   results on  the probability distribution function of the defect function introduced in Section \ref{Int}. We use notation from Section~\ref{Int}.
 An important consideration in finding an optimal $f_\bz$ is how to describe the class $\Theta$ of priors? In other words, what characteristics of $\Theta$ govern, say, the optimal rate of decay of $\bbE(\|f_\rho-f_\bz\|^2_{L_2(\rho_\Omega)})$ for $f_\rho\in \Theta$?
Previous and recent works in statistics and learning theory (see, for instance, \cite{DKPT}, \cite{KoTe}, and \cite{VTbook}, Ch.4) indicate that the compactness characteristics of $\Theta$ play a fundamental role. It is convenient for us to express compactness of $\Theta$ in terms of entropy numbers. We discuss the classical concept of entropy (see Subsection~\ref{Ag}). We note that some other concepts of entropy, for instance, entropy with bracketing, prove to be useful in the theory of empirical processes and nonparametric statistics (see \cite{VG},   \cite{Va}). There is a concept of $VC$ dimension that plays a fundamental role in the problem of pattern recognition and classification \cite{Va}. This concept is also useful in describing compactness characteristics of sets.

In   this subsection we  assume that the measure $\rho$ is concentrated on a set that is bounded with respect to $y$, i.e. the set $Z$ satisfies the condition $Z\subset \Omega\times [-M,M]$ (or a little weaker $|y|\le M$ a.e. with respect to $\rho_\Omega$, i.e the $\rho_\Omega$-measure of those $\bx$, for which there exists a $y$ such that $(\bx,y)\in Z$ and $|y|>M$ is equal to zero) with some fixed $M$. Then it is clear that for $f_\rho$ we have $|f_\rho(\bx)|\le M$ for all $\bx$ (for almost all $\bx$). Therefore, it is natural to assume that a class $\Theta$ of priors where $f_\rho$ belongs is embedded into the $\cC(X)$-ball ($L_\infty$-ball) of radius $M$.

We begin with the case of
 $\cC(\Omega)$, the space of functions continuous on a compact subset $\Omega$ of $\R^d$ with the norm
$$
\|f\|_\infty:= \sup_{\bx\in \Omega}|f(\bx)|.
$$
  We use the abbreviated notation
$$
  \e_n(\bF):= \e_n(\bF,\cC).
$$
The following well known theorem (see, for instance, \cite{VTbook}, p.211) shows how compactness characteristics of $\bF$ can be used in estimating the defect function.
 \begin{Theorem}\label{BT1} Let $\bF$ be a compact subset of $\cC(\Omega)$. Assume that there exists  $M>0$ such that  for all $f\in \bF$ and for all $(\bx,y)\in Z$ we have $|f(\bx)-y| \le M$.  Then, for all $\e>0$
\begin{equation}\label{B5}
\rho^m\{\bz:\sup_{f\in \bF}|L_\bz(f)|\le\e\} \ge1- 2N\left(\bF,\frac{\e}{8M}\right)\exp\left(-\frac{m\e^2}{8(\sigma^2+M^2\e/6)}\right).
\end{equation}
Here $\sigma^2:=  \sup_{f\in \bF}\sigma^2((f(\bx)-y)^2)$ and $\sigma^2(g)$ is the variance of a random variable $g$.
\end{Theorem}
\begin{Remark}\label{BR1} In general we cannot guarantee that the set \newline $\{\bz:\sup_{f\in \bF}|L_\bz(f)|\ge\eta\}$ is $\rho^m$-measurable. In such a case the relation (\ref{B5}) and further relations of this type are understood in the sense of outer measure associated with the $\rho^m$. For instance, for (\ref{B5}) this means that there exists $\rho^m$-measurable set $G$ such that $\{\bz:\sup_{f\in \bF}|L_\bz(f)|\ge\eta\}\subset G$ and (\ref{B5}) holds for $G$.
\end{Remark}

We note that Theorem \ref{BT1} is related to the concept of
 the Glivenko-Cantelli sample complexity of a class $\Phi$ with accuracy $\e$ and confidence $\de$:
$$
S_\Phi(\e,\de):= \min\{n:\quad\text{ for \, all}\quad m\ge n, \quad \text{for \, all}\quad \rho
$$
$$
 \rho^m\{\bz=(z_1,\dots,z_m):\sup_{\phi\in \Phi}|\int_Z \phi d\rho-\frac{1}{m}\sum_{i=1}^m\phi(z_i)|\ge\e\} \le \de\}.
$$
In order to see that, we define points $z_i:=(\bx^i,y_i)$, $i=1,\dots,m$; functions $\phi(\bx,y):=(f(\bx)-y)^2$; and, finally, the class $\Phi:=\{(f(\bx)-y)^2, f\in \bF\}$. One can find a survey of results on the Glivenko-Cantelli sample complexity in \cite{Me} and find results and the corresponding historical remarks related to Theorem \ref{BT1}  in \cite{GKKW}.

We now formulate two theorems, which provide somewhat more delicate estimates for the defect function (see \cite{VTbook}, pp. 213--217). We assume that $\rho$ and $\bF$ satisfy the following condition.
\begin{equation}\label{B6}
\text{For all}\quad f\in \bF,\quad f:\Omega\to Y\quad \text{and any}\, (\bx,y)\in Z,\quad |f(\bx)-y| \le M.
\end{equation}
The following Theorem \ref{BT2} is from \cite{VT98} (see also \cite{VTbook}, Section 4.3.3, p.213).
 We point out that the proofs of both Theorems \ref{BT2} and \ref{BT3} are based on the chaining technique and on the concentration measure inequalities. Some of the ideas here are similar {to the previously discussed ideas}, which were used for analyzing Monte Carlo
cubature formulas.

\begin{Theorem}\label{BT2} Assume that $\rho$, $\bF$ satisfy (\ref{B6}) and $\bF$ is such that
\begin{equation*}
\sum_{n=1}^\infty n^{-1/2}\e_n(\bF) <\infty.
\end{equation*}
Then for $m\eta^2\ge 1$ we have
$$
\rho^m\{\bz:\sup_{f\in \bF}|L_\bz(f)|\ge \eta\} \le C(M,\e(\bF))\exp(-c(M)m\eta^2)
$$
with $C(M,\e(\bF))$ that may depend on $M$ and $\e(\bF):=\{\e_n(\bF,\cC)\}$; $c(M)$ may depend only on $M$.
\end{Theorem}

\begin{Theorem}\label{BT3} Assume that $\rho$, $\bF$ satisfy (\ref{B6}) and $\bF$ is such that
$$
\sum_{n=1}^\infty n^{-1/2}\e_n(\bF) =\infty.
$$
For $\eta>0$ define $J:=J(\eta/M)$ as the minimal $j$ satisfying $\e_{2^j}\le \eta/(8M)$ and
$$
S_J:= \sum_{j=1}^J2^{(j+1)/2}\e_{2^{j-1}}.
$$
Then for $m$, $\eta$ satisfying $m(\eta/S_J)^2 \ge 480M^2$ we have
$$
\rho^m\{\bz:\sup_{f\in \bF} |L_\bz(f)|\ge \eta\} \le C(M,\e(\bF))\exp(-c(M)m(\eta/S_J)^2).
$$
\end{Theorem}

\begin{Corollary}\label{BC1} Assume $\rho$, $\bF$ satisfy (\ref{B6}) and $\e_n(\bF)\le Dn^{-1/2}$.
Then for $m$, $\eta$ satisfying $m(\eta/(1+\log(M/\eta)))^2 \ge C_1(M,D)$ we have
$$
\rho^m\{\bz:\sup_{f\in \bF} |L_\bz(f)|\ge \eta\} \le C(M,D)\exp(-c(M,D)m(\eta/(1+\log (M/\eta)))^2).
$$
\end{Corollary}
\begin{Corollary}\label{BC2} Assume $\rho$, $\bF$ satisfy (\ref{B6}) and
$\e_n(\bF)\le Dn^{-r}$, $r\in (0,1/2)$.
Then for $m$, $\eta$ satisfying $m \eta^{1/r} \ge C_1(M,D,r)$ we have
$$
\rho^m\{\bz:\sup_{f\in \bF} |L_\bz(f)|\ge \eta\} \le C(M,D,r)\exp(-c(M,D,r)m\eta^{1/r} ).
$$
\end{Corollary}

We now demonstrate how the above results from learning theory can be used for sampling discretization of the $L_2$ norm. We use notation from Section \ref{Int}.
 Condition (\ref{B6}) is satisfied with $M$ such that for all $f\in \bF$ we have $\|f\|_\infty \le M$. The argument from the Introduction shows that we can derive results on discretization of the $L_2$ norm directly from the corresponding results from learning theory. We assume that $\bF$ satisfies the following condition:
\be\label{B8}
f\in \bF\quad \Rightarrow \quad \|f\|_\infty \le M.
\ee

Theorem \ref{BT2} implies the following result (see \cite{VT171}).

\begin{Theorem}\label{BT4} Assume that $\bF$ satisfies (\ref{B8}) and the condition
$$
\sum_{n=1}^\infty n^{-1/2} \e_n(\bF) < \infty.
$$
Then there exists a constant $K$ such that for any $m$ there is a set of points $\xi=\{\xi^1,\dots,\xi^m\}$ such that
\begin{equation*}
er_m(\bF,L_2) \le Km^{-1/2}.
\end{equation*}
In particular, if
$$
\e_n(\bF) \le C_1 n^{-r}, \quad r>1/2,
$$
 then
 $$
 er_m(\bF,L_2) \le C(r,C_1) m^{-1/2}.
 $$
\end{Theorem}
 The reader can find further results in this direction in the recent paper \cite{MUl}.

Corollary \ref{BC2} implies the following result (see \cite{VT171}).

\begin{Theorem}\label{BT5} Assume that $\bF$ satisfies (\ref{B8}) and the condition
$$
  \e_n(\bF) \le C_1 n^{-r},\qquad r\in (0,1/2).
$$
Then there exists a constant $C(r,C_1)$ such that for any $m$ there is a set of points $\xi=\{\xi^1,\dots,\xi^m\}$ such that
\begin{equation*}
er_m(\bF,L_2) \le C(r,C_1)m^{-r}.
\end{equation*}
\end{Theorem}

{\bf Acknowledgements.} The authors are grateful to Feng Dai, Tino Ullrich, Mario Ullrich, and the referees for their useful comments and suggestions. {The second author is a Young Russian Mathematics award winner and would like to thank its sponsors and jury.}

The work was supported by the Russian Federation Government Grant N{\textsuperscript{\underline{o}}}14.W03.31.0031.

\end{document}